\documentclass[11pt]{amsart}

\usepackage{framed}
\usepackage{braket}
\usepackage{amsmath,amssymb,amsthm,stmaryrd}
\usepackage{color}
\usepackage{bm}
\usepackage[all]{xy}
\usepackage{amscd}
\usepackage{graphicx}
\usepackage{ascmac}
\usepackage{mathrsfs}
\usepackage[top=30truemm,bottom=30truemm,left=25truemm,right=25truemm]{geometry}
\usepackage{mathtools}
\mathtoolsset{showonlyrefs=true}
\usepackage{tikz}
\usetikzlibrary{tikzmark,calc}
\usetikzlibrary{cd}
\usepackage[
  colorlinks=true,
  linkcolor=red,
  citecolor=blue,
  urlcolor=blue
]{hyperref}

\makeatletter
\@addtoreset{equation}{section}

\makeatother

\newcommand{\Z}{\mathbb{Z}}

\newcommand{\Q}{\mathbb{Q}}
\newcommand{\R}{\mathbb{R}}
\newcommand{\C}{\mathbb{C}}
\newcommand{\F}{\mathbb{F}}

\newcommand{\fM}{\mathfrak{M}}

\newcommand{\fS}{\mathfrak{S}}

\newcommand{\fm}{\mathfrak{m}}

\newcommand{\fp}{\mathfrak{p}}

\newcommand{\fw}{\mathfrak{w}}

\newcommand{\cA}{\mathcal{A}}
\newcommand{\cB}{\mathcal{B}}
\newcommand{\cC}{\mathcal{C}}
\newcommand{\cD}{\mathcal{D}}

\newcommand{\cF}{\mathcal{F}}
\newcommand{\cG}{\mathcal{G}}
\newcommand{\cH}{\mathcal{H}}

\newcommand{\cJ}{\mathcal{J}}

\newcommand{\cM}{\mathcal{M}}

\newcommand{\cO}{\mathcal{O}}

\newcommand{\cR}{\mathcal{R}}

\newcommand{\cT}{\mathcal{T}}

\newcommand{\cX}{\mathcal{X}}

\newcommand{\what}{\widehat}

\newcommand{\parenth}[1]{\left( #1 \right)}

\DeclareMathOperator{\gal}{Gal}
\DeclareMathOperator{\Coker}{Coker}

\DeclareMathOperator{\Ker}{Ker}
\DeclareMathOperator{\Imag}{Im}
\DeclareMathOperator{\core}{cor}
\DeclareMathOperator{\rest}{res}

\DeclareMathOperator{\Cl}{Cl}
\DeclareMathOperator{\ord}{ord}
\DeclareMathOperator{\Hom}{Hom}

\DeclareMathOperator{\rank}{rank}
\DeclareMathOperator{\corank}{corank}

\DeclareMathOperator{\cyc}{cyc}
\DeclareMathOperator{\cha}{char}

\DeclareMathOperator{\Tr}{Tr}
\DeclareMathOperator{\an}{an}
\DeclareMathOperator{\alg}{alg}
\DeclareMathOperator{\Frob}{Frob}

\DeclareMathOperator{\sgn}{sgn}

\DeclareMathOperator{\tr}{Tr}
\DeclareMathOperator{\col}{Col}
\DeclareMathOperator{\sel}{Sel}
\DeclareMathOperator{\iw}{Iw}
\DeclareMathOperator{\augm}{aug}
\DeclareMathOperator{\loc}{loc}

\usepackage[OT2,T1]{fontenc}
\DeclareSymbolFont{cyrletters}{OT2}{wncyr}{m}{n}
\DeclareMathSymbol{\Sha}{\mathalpha}{cyrletters}{"58}

\let\oldenumerate\enumerate
\renewcommand{\enumerate}{
   \oldenumerate
   \setlength{\itemsep}{1pt}
   \setlength{\parskip}{0pt}
   \setlength{\parsep}{0pt}
}
\let\olditemize\itemize
\renewcommand{\itemize}{
   \olditemize
   \setlength{\itemsep}{1pt}
   \setlength{\parskip}{0pt}
   \setlength{\parsep}{0pt}
}

\theoremstyle{plain}
\newtheorem{thm}{Theorem}[section]
\newtheorem{lem}[thm]{Lemma}
\newtheorem{conj}[thm]{Conjecture}
\newtheorem{prop}[thm]{Proposition}
\newtheorem{cor}[thm]{Corollary}

\theoremstyle{definition}
\newtheorem{dfn}[thm]{Definition}

\newtheorem{rem}[thm]{Remark}
\newtheorem{eg}[thm]{Example}

\usepackage{diagbox}
\usepackage{tikz}
\usetikzlibrary{shapes.geometric}
\usepackage{mathtools}
\usepackage{float}
\usepackage{booktabs}

\DeclareMathOperator{\aug}{aug}

\makeatletter
\@namedef{subjclassname@2020}{\textup{2020} Mathematics Subject Classification}
\makeatother

\title
[A 2-Adic Kida-Type Formula and Quadratic Twists]
{A 2-Adic Kida-Type Formula for Chromatic Selmer Groups of Elliptic Curves and Applications to Quadratic Twists}
\author[T.~Adachi]{Taiga Adachi}
\address{Joint Graduate School of Mathematics for Innovation, Kyushu University,
Motooka 744, Nishi-ku, Fukuoka 819-0395, Japan}
\email{t.adachi1729@gmail.com}
\date{\today}

\begin{document}

\maketitle

\begin{abstract}
Let $E$ be an elliptic curve over $\Q$ with good
supersingular reduction at $2$.
We study the chromatic ($\sharp/\flat$) Selmer groups
of $E$ over the cyclotomic $\Z_2$-extension of $\Q$,
defined using Sprung's chromatic Coleman maps.
We show that at least one of these Selmer groups
is $\Lambda$-cotorsion and, assuming the vanishing
of the corresponding $\mu$-invariant, establish a
$2$-adic Kida-type formula describing the change
in the Iwasawa $\lambda$-invariants under base change
from $\Q$ to a real quadratic field $F$ in which
$2$ is unramified.
When $E$ has square-free conductor, applying this
formula to $E$ and its quadratic twist $E^D$
yields an explicit formula describing the variation
of the chromatic $\lambda$-invariants under twists
by positive square-free integers $D\equiv1\pmod4$.
As an application, under suitable hypotheses on
the root number and chromatic Iwasawa invariants
of $E$, we construct explicit sets of primes $\ell$
such that the $2^\infty$-Selmer group of $E^\ell$
over $\Q$ has $\Z_2$-corank one, and determine
their positive natural densities.
\end{abstract}
%
%
\section{Introduction}\label{sec:intro}
\noindent
Let $E$ be an elliptic curve over $\Q$ with good supersingular reduction at $2$, and let $\Q^{\cyc}$ be the cyclotomic $\Z_2$-extension of $\Q$. For a square-free integer $D$, we denote by $E^D$ the quadratic twist of $E$ associated with $\Q(\sqrt{D})/\Q$. Put $\Lambda:=\Z_2\llbracket \gal(\Q^{\cyc}/\Q)\rrbracket\cong \Z_2\llbracket T\rrbracket$.
Let $a_p=a_p(E)$ denote the trace of Frobenius at a good prime $p$ of $E$. For an odd prime $p$ of good supersingular reduction with $a_p=0$, Kobayashi \cite{Kob03} introduced the plus/minus Selmer groups, whose Pontryagin duals are torsion over $\Lambda$ and which provide algebraic counterparts to Pollack's plus/minus $p$-adic $L$-functions \cite{Pol03}. Sprung \cite{Spr12} removed the assumption $a_p=0$ by constructing chromatic ($\sharp/\flat$) Coleman maps and the corresponding $p$-adic $L$-functions, and introduced chromatic Selmer groups using these Coleman maps. While his construction of the Coleman maps and $p$-adic $L$-functions includes the case $p=2$, his treatment of chromatic Selmer groups in \cite[\S 7]{Spr12} assumes that $p$ is odd. 

In this paper, we study the chromatic Selmer groups of $E$ over $\Q^{\cyc}$ and their behavior under base change and quadratic twists. In \S \ref{sec:multi_Selmer}, we construct these Selmer groups at $p=2$ and prove that at least one of them is $\Lambda$-cotorsion and that both are if $a_2=0$ (Theorem \ref{Thm:Lambda_torsion}). Our first main result is a $2$-adic Kida-type formula (Theorem \ref{Main_Kida}) for base change to a real quadratic field $F$ in which $2$ is unramified. We show that the vanishing of the corresponding $\mu$-invariant is preserved under this base change and express the $\lambda$-invariant over the cyclotomic $\Z_2$-extension $F^{\cyc}$ of $F$ in terms of that over $\Q^{\cyc}$ and certain local correction terms. When $E$ has square-free conductor, applying this formula to $E$ and $E^D$ over $\Q(\sqrt{D})$ yields an explicit formula describing the variation of the chromatic $\lambda$-invariants under twists by positive square-free integers $D\equiv1\pmod4$ (Theorem \ref{Main_Thm1}). As applications, we obtain distributional results for the chromatic $\lambda$-invariants in quadratic twist families (Corollary \ref{Cor:density_explicit} and Theorem \ref{Main_Lambda_Distribution}). Furthermore, combining the difference formula with the corank bound in Proposition \ref{prop:chromatic_rank_bound} and the $2$-parity theorem, we construct, under suitable hypotheses on $E$, explicit sets of primes $\ell$ such that $\sel_{2^\infty}(E^\ell/\Q)$ has $\Z_2$-corank one, and determine their positive natural densities (Theorems \ref{Main_Theorem2} and \ref{Main_density_root_1}).

To prove the Kida-type formula, we compute Herbrand quotients, following the approach of Hachimori--Matsuno \cite{HM99} and Matsuno \cite{Mat08}. Unlike in the ordinary case, the local conditions at $2$ are defined through the kernels of the chromatic Coleman maps, and their compatibility with restriction under base change requires a direct comparison of these maps under corestriction. Moreover, the presence of real places makes the $2$-cohomological dimension of the relevant global Galois groups infinite. In particular, the treatment of real local conditions is essential for controlling the $\mu$-invariants when the minimal discriminant of $E$ is positive. A further difficulty in deriving the difference formula under quadratic twists arises when $D\equiv5\pmod8$, in which case $a_2(E^D)=-a_2(E)$. We show that the isomorphism between $E$ and $E^D$ over $\Q(\sqrt{D})$ identifies their chromatic local conditions at $2$ with the same sign ($\sharp/\flat$).

We introduce some notation to state our main results. We fix algebraic closures $\overline{\Q},\overline{\Q_2}$, an embedding $\overline{\Q}\hookrightarrow \overline{\Q_2}$, and let $v_2$ be the corresponding $2$-adic valuation normalized by $v_2(2)=1$. For an elliptic curve $A$ defined over $\Q$ with good supersingular reduction at $2$, let $\sel^{\sharp/\flat}(A/\Q^{\cyc})$ be the sharp/flat (chromatic) Selmer groups. For $\star\in\{\sharp,\flat\}$, whenever $\sel^{\star}(A/\Q^{\cyc})$ is $\Lambda$-cotorsion, we denote its Iwasawa $\mu$- and $\lambda$-invariants by $\mu^{\star}(A/\Q)$ and $\lambda^{\star}(A/\Q)$, respectively. Since $2$ is unramified in $F$, restriction induces a natural isomorphism $\gal(F^{\cyc}/F)\cong\gal(\Q^{\cyc}/\Q)$, through which we identify the corresponding Iwasawa algebras. We also consider the doubly-signed Selmer group $\sel^{\vec{s}_{\star}}(A/F^{\mathrm{cyc}})$ for $\star\in\{\sharp,\flat\}$, where $\vec{s}_{\star}$ specifies the $\star$-local condition at every prime of $F$ above $2$. When this Selmer group is $\Lambda$-cotorsion, we denote its Iwasawa $\mu$- and $\lambda$-invariants by $\mu^{\vec{s}_{\star}}(A/F)$ and $\lambda^{\vec{s}_{\star}}(A/F)$, respectively.

\begin{thm}[Theorem \ref{Thm:Matsuno-type_formula}]\label{Main_Kida}
Let $E$ be an elliptic curve defined over $\Q$ with good supersingular reduction at $2$. Let $\star\in\{\sharp,\flat\}$ and let $F$ be a real quadratic field in which $2$ is unramified. For each finite prime $\fw\nmid 2$ of $F^{\cyc}$, put
\begin{align}\label{eq:modified_ramification_index_intro}
e'_E(\fw):=
\begin{cases}
0&\text{if $E$ has additive reduction at $v$},\\
e_{F^{\cyc}/\Q^{\cyc}}(\fw)&\text{otherwise},
\end{cases}
\end{align}
where $e_{F^{\cyc}/\Q^{\cyc}}(\fw)$ is the ramification index of $\fw$ in $F^{\cyc}/\Q^{\cyc}$ and $v$ is the prime of $\Q^{\cyc}$ below $\fw$. If $\sel^{\star}(E/\Q^{\cyc})$ is $\Lambda$-cotorsion and $\mu^{\star}(E/\Q)=0$, then $\sel^{\vec{s}_{\star}}(E/F^{\cyc})$ is also $\Lambda$-cotorsion and $\mu^{\vec{s}_{\star}}(E/F)=0$. Moreover, we have
\begin{align}\label{Kida_p=2_chromatic_algebraic}
\lambda^{\vec{s}_\star}(E/F)=2\lambda^{\star}(E/\Q)+\sum_{\fw\in P_1}(e'_E(\fw)-1)+2\sum_{\fw\in P_2}(e'_E(\fw)-1),
\end{align}
where $P_1$ and $P_2$ are the sets of primes of $F^{\cyc}$ defined by
\begin{align}
&P_1:=\{\fw\ |\ \fw\nmid 2,\ E\text{ has split multiplicative reduction at }\fw \},\\
&P_2:=\{\fw\ |\ \fw\nmid 2,\ E\text{ has good reduction at }\fw,\ E[2]((F^{\cyc})_{\fw})\neq\{O\}\}.
\end{align}
\end{thm}

Only primes $\fw$ above odd primes ramified in $F/\Q$ contribute to the sums in \eqref{Kida_p=2_chromatic_algebraic}; hence these sums are finite.
Let $\ell$ be the odd rational prime below $\fw$.
Since $(F^{\cyc})_{\fw}$ contains the unramified quadratic extension of $\Q_\ell$, multiplicative reduction at $\fw$ is automatically split.
The modified index $e'_E(\fw)$ in \eqref{eq:modified_ramification_index_intro} accounts for the case where $E$ has additive reduction at $v$ but multiplicative or good reduction at $\fw$.

The above Kida-type formula is a supersingular analogue of Matsuno's result \cite[Theorem 3.1]{Mat08}, which extends the work of Hachimori--Matsuno \cite{HM99} to the case $p=2$ in the ordinary case.
For odd supersingular primes, Pollack--Weston
\cite{PW06} and Hatley--Lei \cite{HL19}
established Kida-type formulae for signed Selmer
groups. More recently, building on Lim's work \cite{Lim22}, Forr\'{a}s--M\"{u}ller
\cite{FM25} extended such results to certain
weakly ramified base fields for $p\geq5$. Kataoka \cite{Kat24+} developed a general framework
for proving Kida-type formulae via Selmer complexes,
encompassing a wide range of settings, including
the supersingular case. On the analytic side, corresponding Kida-type
formulae for signed $p$-adic $L$-functions
were obtained by Pollack--Weston \cite{PW06}
and, more recently, by Pratap--Ray \cite{PR26}. 

\begin{rem}
As noted in \cite[Remark 4.8]{Kat24+}, Kataoka's approach via Selmer complexes can, in principle, also be applied when $p=2$ and the base field is totally imaginary. To the best of the author's knowledge, however, the existing results do not directly cover our setting, which involves the prime $p=2$, real quadratic extensions, and chromatic local conditions for elliptic curves with good supersingular reduction at $2$. In contrast to the totally imaginary case,
the presence of real places must be taken
into account in the global cohomological arguments.
\end{rem}

For a square-free integer $D>1$ with $D\equiv 1\pmod{4}$, the quadratic extension
$\Q(\sqrt{D})/\Q$ is unramified at $2$, and
$E^D$ also has good supersingular reduction at $2$.
Thus, applying Theorem \ref{Main_Kida} to both $E$ and $E^D$
over $\Q(\sqrt{D})$, we obtain the following difference formula.

\begin{thm}[Theorem \ref{chromatic_Matsuno-type_formula}]\label{Main_Thm1}
Let $E$ be an elliptic curve defined over $\Q$ with good supersingular reduction at $2$ and let $D>0$ be a square-free integer such that $D\equiv 1\pmod 4$. Fix $\star\in\{\sharp,\flat\}$. We put, for an odd prime $\ell$,
\begin{align}
n_{\ell}:=v_2\left(\frac{\ell^2-1}{8}\right).
\end{align}
Assume that the conductor $N_E$ of $E$ is square-free and $\sel^{\star}(E/\Q^{\cyc})$ is $\Lambda$-cotorsion. If $\mu^{\star}(E/\Q)=0$, then $\mu^{\star}(E^D/\Q)=0$ and we have
\begin{align}\label{eq:difference_formula}
\lambda^{\star}(E^D/\Q)=\lambda^{\star}(E/\Q)+\sum_{\substack{\ell\mid D\\ \ell\mid N_E}}2^{n_{\ell}}+\sum_{\substack{\ell\mid D\\ \ell\nmid N_E, 2\mid \sharp \widetilde{E}_{\ell}(\F_{\ell})}}2^{n_{\ell}+1}.
\end{align}
Here, $\widetilde{E}_{\ell}$ is the reduction of $E$ modulo $\ell$.
\end{thm}

This is a supersingular analogue of Matsuno's formula \cite[Theorem 5.1]{Mat08} in the ordinary case. The author \cite[Theorem 1.3]{Ada26+} obtained the analytic counterpart of Theorem \ref{Main_Thm1}.
The arguments in \cite[\S\S 3--4]{Ada26+}
also apply to the algebraic difference formula
\eqref{eq:difference_formula}, yielding the
following algebraic counterparts of the analytic
distribution results obtained there.

\begin{cor}[Corollary \ref{Cor:Adachi26+_1.4}]\label{Cor:density_explicit}
Let $E$ be an elliptic curve defined over $\Q$ with good supersingular reduction at $2$ and square-free conductor, and let $\star\in \{\sharp,\flat\}$. Assume that $\sel^{\star}(E/\Q^{\cyc})$ is $\Lambda$-cotorsion. If $\mu^{\star}(E/\Q)=0$, then the set
\begin{align}
\{\lambda^{\star}(E^{\ell}/\Q)\ |\ \ell\text{ is a rational prime and }\ell\equiv 1\pmod 4\}
\end{align}
is infinite. More precisely, for each nonnegative integer $k$, the natural density (in the set of all primes) of
\begin{align}
\{\ell\text{ is a rational prime and }\ell\equiv 1\pmod 4\ |\ \lambda^{\star}(E^{\ell}/\Q)=\lambda^{\star}(E/\Q)+2^{k+1}\}
\end{align}
is $\frac{1}{3\cdot 2^{k+1}}$.
\end{cor} 
Corollary \ref{Cor:density_explicit} gives a supersingular analogue of
\cite[Proposition 6.1]{Mat08}, strengthening the
infinitude result by determining the exact natural
density of primes yielding each prescribed
$\lambda$-invariant.
The author obtained the corresponding analytic result
for chromatic $2$-adic $L$-functions in
\cite[Corollary 1.4]{Ada26+}.

Hatley--Ray \cite{HR24+} studied quadratic twists of elliptic curves with good ordinary
reduction at $2$ using Matsuno's formula.
They obtained results on the distribution of Iwasawa
$\lambda$-invariants \cite[Theorem 4.1]{HR24+} and on prime twists with
$2^\infty$-Selmer corank one \cite[Theorems 3.13 and 3.15]{HR24+}.
They also suggested extending their approach 
to the technically more delicate supersingular case
at $p=2$, highlighting the role
of Sprung's chromatic ($\sharp/\flat$) Selmer groups
\cite[pp.~4--5]{HR24+}. Our results provide a concrete realization
of this direction. The counting argument in \cite[\S 4]{Ada26+}
is based on a refinement of the Selberg--Delange
argument used by Hatley--Ray \cite{HR24+}. Theorem \ref{Main_Lambda_Distribution}
provides a supersingular counterpart of
\cite[Theorem 4.1]{HR24+}, giving lower bounds
for the number of quadratic twists with a
prescribed chromatic $\lambda$-invariant.
Furthermore, Theorem \ref{Main_Theorem2}
provides a supersingular counterpart of
\cite[Theorem 3.13]{HR24+}, yielding a
positive-density set of prime twists whose
$2^\infty$-Selmer groups have $\Z_2$-corank one.

Throughout this paper, the notation $f(X)\gg g(X)$ means that
there exists a constant $C>0$ such that
$f(X)\ge Cg(X)$ for all sufficiently large $X$.
The constant $C$ may depend on the fixed parameters,
but not on $X$.

\begin{thm}[Theorem \ref{Thm:Ada26+_1.5}]\label{Main_Lambda_Distribution}
Let $E$ be an elliptic curve defined over $\Q$ with good supersingular reduction at $2$ and square-free conductor, and let $\star\in \{\sharp,\flat\}$. Assume that $\sel^{\star}(E/\Q^{\cyc})$ is $\Lambda$-cotorsion and $\mu^{\star}(E/\Q)=0$. Let $N^{\star}$ be a nonnegative integer such that $N^{\star}\geq\lambda^{\star}(E/\Q)$ and $\lambda^{\star}(E/\Q)\equiv N^{\star}\pmod 2$. Let
\begin{align}
m^{\star}_{E,N^{\star}}(X):=\sharp\{D\leq X\ |\ D\in\Z_{>0}\text{ is square-free, }D\equiv1\pmod 4,\text{ and }\lambda^{\star}(E^D/\Q)=N^{\star}\}.
\end{align}
Then we have
\begin{align}
m^{\star}_{E,N^{\star}}(X)\gg \frac{X}{(\log X)^{\frac{2}{3}}}.
\end{align}
\end{thm}

Using the difference formula \eqref{eq:difference_formula} for the
$\lambda$-invariants, we obtain the following result on the behavior of the $\Z_2$-corank of $2^{\infty}$-Selmer groups over $\Q$ under quadratic twists. Let $\omega(E)$ be the root number of $E$, and let $\Delta_E$ denote the minimal discriminant of $E$.

\begin{thm}[Theorem \ref{thm:prime_twists_corank_one}]\label{Main_Theorem2}
Let $E/\Q$ be an elliptic curve with good supersingular reduction at $2$ and square-free conductor $N_E$, and let $\star\in\{\sharp,\flat\}$. Assume that $\sel^{\star}(E/\Q^{\cyc})$ is $\Lambda$-cotorsion and
\begin{align}
\omega(E)=-1,\qquad
\mu^\star(E/\Q)=0,\qquad
\lambda^\star(E/\Q)\leq2.
\notag\end{align}
Let $\cM_E^-$ be the set of rational primes $\ell$ satisfying
\begin{itemize}
\item $\ell\nmid 2N_E$,
\item $\widetilde E_{\ell}(\F_\ell)[2]=\{\widetilde{O}\}$,
\item $\ell$ splits completely in $\Q(i,\sqrt{-N_E})$.
\end{itemize}
Then, for every $\ell\in\cM_E^-$, we have
\begin{align}
\corank_{\Z_2}\sel_{2^\infty}(E^\ell/\Q)=1.
\notag\end{align}
Moreover, $\cM_E^-$ has natural density
\begin{align}\label{eq:density_prime_twists_intro}
\operatorname{dens}(\cM_E^-)=
\begin{cases}
\frac{1}{6}&\text{if }\Q(\sqrt{\Delta_E})\subseteq \Q(i,\sqrt{-N_E}),\\[4pt]
\frac{1}{12}&\text{otherwise}.
\end{cases}
\end{align}
\end{thm}

The set $\cM_E^-$ in Theorem \ref{Main_Theorem2} is closely related to the
prime-twist families studied by Barrera Salazar--Pacetti--Tornar\'{i}a
\cite[Theorems 3.3 and 3.7]{BPT21}. Under the additional local hypotheses
of \cite[Hypotheses 2.1]{BPT21}, their results imply that
$\sel_2(E^\ell/\Q)=\sel_2(E/\Q)$ for every $\ell\in\cM_E^-$.
Consequently, if $\omega(E)=-1$ and
$\dim_{\F_2}\sel_2(E/\Q)=1$, the $2$-parity theorem
\cite[Theorem 1.4]{DD10} gives
$\corank_{\Z_2}\sel_{2^\infty}(E^\ell/\Q)=1$
for all $\ell\in\cM_E^-$. In contrast, Theorem \ref{Main_Theorem2}
requires neither the local hypotheses of \cite{BPT21} nor the assumption
$\dim_{\F_2}\sel_2(E/\Q)=1$, and hence is not subsumed by their results.
Moreover, we determine the exact natural density of $\cM_E^-$ in
\eqref{eq:density_prime_twists_intro}.

Theorem \ref{Main_Theorem2} also extends to square-free products of primes in $\cM_E^-$.
This yields a supersingular analogue of \cite[Theorem 3.12]{HR24+}, with a sharper counting lower bound in certain cases.

\begin{cor}[Corollary \ref{cor:composite_twists_corank_one}]
\label{cor:composite_twists_corank_one_intro}
Assume the hypotheses of Theorem \ref{Main_Theorem2}.
Let $\cD_E$ be the set of square-free integers $D>1$
whose prime factors all belong to $\cM_E^-$.
Then, for $D\in \cD_E$, we have
\begin{align}
\corank_{\Z_2}\sel_{2^\infty}(E^D/\Q)=1.
\notag\end{align}
Moreover,
\begin{align}
\sharp\{D\in\cD_E\mid D\le X\}\gg
\begin{cases}
\dfrac{X}{(\log X)^{5/6}}
&\text{if }\Q(\sqrt{\Delta_E})\subseteq\Q(i,\sqrt{-N_E}),\\[8pt]
\dfrac{X}{(\log X)^{11/12}}
&\text{otherwise}.
\end{cases}
\notag\end{align}
\end{cor}

The following theorem provides a different explicit positive-density set
$\cM_E^+$ of prime twists, which is not covered by the inert-prime families
of \cite{BPT21}, since $\widetilde E_{\ell}(\F_\ell)[2]\neq\{\widetilde{O}\}$ for every
$\ell\in\cM_E^+$. 
It also treats the case $\omega(E)=+1$ by choosing
primes $\ell$ for which the local correction term
in \eqref{eq:difference_formula} equals $2$
and $\omega(E^\ell)=-1$.
Under the assumptions $\mu^\star(E/\Q)=\lambda^\star(E/\Q)=0$,
the difference formula and the $2$-parity theorem,
together with the corank bound, then imply that
$\sel_{2^\infty}(E^\ell/\Q)$ has $\Z_2$-corank one.
Hatley--Ray \cite[Theorem 3.15]{HR24+} obtained a root-number $+1$ result
for elliptic curves with good ordinary reduction at $2$, assuming
$E(\Q)[2]\neq\{O\}$. Our result treats instead good supersingular reduction
at $2$, where $E(\Q)[2]=\{O\}$, and gives the exact density of the resulting
prime-twist family.

\begin{thm}[Theorem {\ref{thm:prime_twists_corank_one_plus}}]\label{Main_density_root_1}
Let $E/\Q$ be an elliptic curve with good supersingular reduction at $2$ and square-free conductor $N_E$, and let $\star\in\{\sharp,\flat\}$. Assume that $\sel^{\star}(E/\Q^{\cyc})$ is $\Lambda$-cotorsion and
\begin{align}
\omega(E)=+1,\qquad
\mu^\star(E/\Q)=0,\qquad
\lambda^\star(E/\Q)=0.
\notag\end{align}
Let $\cM_E^+$ be the set of rational primes $\ell$ satisfying
\begin{itemize}
\item $\ell\nmid2N_E$,
\item $\ell\equiv5\pmod8$,
\item $\widetilde E_{\ell}(\F_\ell)[2]\neq\{\widetilde{O}\}$,
\item $\displaystyle\left(\frac{-N_E}{\ell}\right)=-1$.
\end{itemize}
Then, for every $\ell\in\cM_E^+$, we have
\begin{align}\label{eq:prime_twists_plus_conclusion_intro}
\mu^\star(E^\ell/\Q)=0,\qquad
\lambda^\star(E^\ell/\Q)=2,\qquad
\corank_{\Z_2}\sel_{2^\infty}(E^\ell/\Q)=1.
\end{align}
Moreover, $\cM_E^+$ has natural density
\begin{align}\label{eq:density_prime_twists_plus_intro}
\operatorname{dens}(\cM_E^+)=
\begin{cases}
\frac{1}{8}&\text{if }\Q(\sqrt{\Delta_E})\subseteq \Q(i,\sqrt{-N_E}),\\[4pt]
\frac{1}{12}&\text{otherwise}.
\end{cases}
\end{align}
\end{thm}

Theorems \ref{Main_Theorem2} and \ref{Main_density_root_1} are also related to Goldfeld's
conjecture \cite{Gol79} on the distribution
of ranks of quadratic twists. Recently,
Smith \cite{Smi25+} proved the $50$--$50$
distribution of $2^\infty$-Selmer coranks
in the full family of quadratic twists of
any elliptic curve over $\Q$. 
Our approach instead gives explicit families of prime
twists and their square-free products.

\subsection*{Organization of the paper}
In \S \ref{sec:multi_Selmer}, we construct the $2$-adic chromatic Coleman maps and doubly-signed Selmer groups. In \S \ref{sec:Kida}, we prove a $2$-adic Kida-type formula for chromatic Selmer groups over real quadratic extensions using Herbrand quotients, following \cite{HM99,Mat08}. We then derive a difference formula for the $\lambda$-invariants under quadratic twists. Combining this result with the arguments in \cite{HR24+} and \cite[\S\S 3,4]{Ada26+}, we also obtain several density results. 
In \S \ref{sec:MW-rank_twists}, as an application of the difference formula obtained in \S \ref{sec:Kida}, we establish density results for prime twists with $2^\infty$-Selmer corank one, using the method of \cite{HR24+}.
We also obtain a counting lower bound for square-free products of primes in $\cM_E^-$ whose twists have $2^\infty$-Selmer corank one.
In \S \ref{sec:examples}, we use data from LMFDB
\cite{lmfdb} and Pollack's SageMath code
\cite{PolCode} to give examples of
Corollaries \ref{cor:prime_twists_analytic} and
\ref{cor:prime_twists_plus_analytic}, which give criteria in terms of analytic Iwasawa
invariants for the Selmer-corank conclusions of
Theorems \ref{Main_Theorem2} and
\ref{Main_density_root_1}, respectively.

\subsection*{Acknowledgments}
The author is grateful to his supervisor, Shinichi Kobayashi, for helpful comments, encouragement, and insights from the paper \cite{Kob03}. The author would also like to express sincere gratitude to Takenori Kataoka for answering many questions about his papers \cite{Kat21, Kat24+}, for insightful comments, and for continued encouragement. The author would also like to thank Florian Sprung for answering questions about his papers \cite{Spr12,Spr24}. The author is also grateful to Ben Forr\'{a}s for helpful correspondence
concerning his joint work \cite{FM25} with Katharina M\"{u}ller. Part of this work was carried out during the author's visit to Seoul National
University from December 2025 to January 2026. The author is deeply grateful
to Dohyeong Kim for kindly hosting the visit and for many insightful comments. Part of this work was also carried out during the author's visit to the
University of Ottawa in October 2024. The author is deeply grateful to
Antonio Lei for his helpful comments and generous support during the visit. He is also thankful to Ashay Burungale and Shoma Sueyoshi for providing useful information on related work. The author thanks Chan-Ho Kim, Mahiro Atsuta, Naoto Dainobu, and Ryota Shii for helpful discussions. The author is also grateful to Keiichiro Nomoto and Ryota Shii, whose joint work \cite{ANS26}
with the author provided the initial impetus for the research project
encompassing \cite{Ada26+} and the present work. The author was supported by the WISE program (MEXT) at Kyushu University.

\subsection*{Use of AI}
After completing the first draft, the author used
ChatGPT-6 Sol to check parts of the arguments.
It pointed out that the compatibility of the Coleman map
in \S~\ref{subsec:Coleman} with Sprung's rational Coleman map
\cite{Spr12} and with the $2$-adic $L$-functions
of \cite{Spr17} had to be verified
(Lemmas \ref{lem:Col_Kato_comparison}
and \ref{lem:Col_analytic_comparison}),
and it drew the author's attention to the works
discussed in Remark \ref{Rem:relation_previous_works}.
The author has checked all arguments and takes
full responsibility for the content of this paper.

\section{$2$-adic doubly-signed Selmer groups}\label{sec:multi_Selmer}
In this section, we construct the chromatic local condition at $p=2$ based on the method of \cite{Spr12}, \cite[\S\S3,4]{Kat21} and define the doubly-signed Selmer group. The constructions and results discussed in this section originate in the work of Kobayashi \cite{Kob03} and have subsequently been generalized and further developed in many papers; see, for example, \cite{IV06}, \cite{Kim07}, \cite{Kim11}, \cite{Kim13}, and \cite{KO18}. 

Let $F/\Q$ be a quadratic extension in which $2$ is unramified. For a prime $\fp\mid 2$ of $F$, let $F_{\fp}$ denote the completion of $F$ at $\fp$. We also fix a generator $(\zeta_{2^n})_{n\geq 0}$ of $\Z_2(1)$, namely $\zeta_{2^n}$ is a primitive $2^n$-th root of unity, and $\zeta_{2^{n+1}}^2=\zeta_{2^n}$ and $\zeta_{2^0}=1$. For $n\geq-2$, put $k_n:=F_{\fp}(\zeta_{2^{n+2}})$, let $\cO_n$ be its ring of integers, and let $\fm_n$ be the maximal ideal of $\cO_n$. Let $k_{\infty}:=\bigcup_{n\geq-2}k_n$. For simplicity, put $\cO:=\cO_{-2}$ and $\fm:=\fm_{-2}$. Let $\varphi\in\gal(F_{\fp}/\Q_2)$ be the Frobenius automorphism, i.e., $\varphi(x)\equiv x^2\pmod {2\cO}$ for $x\in\cO$. We put, for $n\geq-2$, $G_n:=\gal(F_{\fp}(\zeta_{2^{n+2}})/\Q_2)$, and $\Gamma:=\gal(k_{\infty}/k_0)\cong \Z_2$. Fix a topological generator $\gamma$ of $\Gamma$. We identify $\Lambda$ with $\Z_2\llbracket T\rrbracket$
via $\gamma\mapsto 1+T$.

\begin{rem}
When $F_{\fp}=\Q_2$, our construction agrees with
that of Sprung \cite{Spr12}.
When $F_{\fp}/\Q_2$ is the unramified quadratic extension,
it gives a $p=2$ analogue of the local construction
of Kataoka \cite[\S\S3,4]{Kat21}.
\end{rem}

\subsection{A system of local points}\label{subsec:local_points}
Let $\zeta$ be a primitive third root of unity if
$[F_{\fp}:\Q_2]=2$, and let $\zeta=1$ if $F_{\fp}=\Q_2$. Let $E$ be an elliptic curve defined over $\Q_2$ with good supersingular reduction, and let $a_2:=3-\sharp\widetilde{E}(\F_2)$ be the trace of Frobenius, where $\widetilde{E}$ is the reduction of $E$ modulo $2$. Let $\widehat{E}$ be the formal group law over $\Z_2$ associated with a fixed minimal Weierstrass model of $E$. For a formal group $\cF$ over $\cO$, let $\exp_{\cF}$ and $\log_{\cF}$ be the formal exponential and logarithm, respectively.

\begin{dfn}[cf. {\cite[Definition 3.19]{Kat21}}]\label{Def:eta_system}
For $n\geq-2$, we define
\begin{align}
\eta_n:= \frac{2}{9-a_2^2}
\left(
3\zeta^{\varphi^{-n}}
+a_2\zeta^{\varphi^{-(n+1)}}
\right).
\end{align}
\end{dfn}

\begin{rem}\label{Rem:phi_relation_eta_n}
Since $\varphi^2=1$, we have
\begin{align}
\eta_{n+1}=\eta_n^\varphi,\qquad
\eta_{n+2}=\eta_n
\end{align}
for every $n\geq-2$. Moreover,
\begin{align}
(3-a_2\varphi)\eta_n=2\zeta^{\varphi^{-n}}.
\end{align}
When $F_{\fp}=\Q_2$, we have
\begin{align}
\eta_n=\frac{2}{3-a_2},
\end{align}
which coincides with the value
$\log_{\mathcal F_{\mathrm{ss}}}(\varepsilon)$
used in \cite[Proof of Theorem 2.2]{Spr12}.
\end{rem}

\begin{dfn}\label{Def:sequence}
We define a sequence $\{x_j\}_{j\geq-1}$ in $\Q_2$ by $x_{-1}:=0$, $x_0:=1$, and $x_j:=2^{-1}(a_2x_{j-1}-x_{j-2})$ for $j\geq1$.
\end{dfn}
\begin{dfn}\label{Def_power}
We define a power series
\begin{align}
f(T):=\sum_{j\geq0}x_j((T+\zeta)^{2^j}-\zeta^{2^j})\in F_{\fp}\llbracket T\rrbracket.
\end{align}
\end{dfn}

\begin{lem}
For $n\geq-2$, $f^{\varphi^{-n}}$ is of the Honda type $t^2-a_2 t+2$.
\end{lem}
\begin{proof}
See \cite[Proof of Proposition 3.22]{Kat21}, for example.
\end{proof}

By Honda theory \cite{Hon70}, for each $n\geq-2$,
there exists a unique formal group $\cG_n$ over $\cO$
such that $\log_{\cG_n}=f^{\varphi^{-n}}$.

\begin{lem}\label{lem:formal_log_unramified}
For every $n\geq-2$, the formal logarithm induces
an isomorphism
\begin{align}
\log_{\cG_n}:\cG_n(\fm)\xrightarrow{\sim}\fm.
\end{align}
\end{lem}

\begin{proof}
Since $F_{\fp}/\Q_2$ is unramified, we have
$\fm=2\cO$. By
\cite[Chapter IV, Theorem 6.4(b)]{Sil09},
$\log_{\cG_n}$ is an isomorphism on $4\cO$.

Put $u_n:=\log_{\cG_n}'(0)$. By Definition
\ref{Def:sequence}, we have $u_n\equiv1\pmod{2\cO}$,
hence $u_n\in\cO^\times$.
Since $\cG_n$ has height $2$, we have, for $x\in2\cO$,
\begin{align}
[2]_{\cG_n}(x)\equiv2x\pmod{8\cO}.
\end{align}
Moreover, for $y\in4\cO$, the integrality of
the invariant differential gives
$\log_{\cG_n}(y)\equiv u_ny\pmod{8\cO}$.
Consequently, for $x\in2\cO$,
\begin{align}
2\log_{\cG_n}(x)
&=\log_{\cG_n}([2]_{\cG_n}(x))\notag\\
&\equiv2u_nx\pmod{8\cO}.
\end{align}
Thus
\begin{align}
\log_{\cG_n}(x)\equiv u_nx\pmod{4\cO}.
\end{align}
Since $u_n$ is a unit, the induced map on
$\cG_n(2\cO)/\cG_n(4\cO)$ is an isomorphism
onto $2\cO/4\cO$. Together with the
isomorphism on $4\cO$, this proves the assertion.
\end{proof}

By Lemma \ref{lem:formal_log_unramified},
there exists a unique element
$\varepsilon_n\in\cG_n(\fm)$ such that
$\log_{\cG_n}(\varepsilon_n)=\eta_n$.

\begin{dfn}
For $n\geq-2$, we define $\pi_n\in \fm_n$ by
\begin{align}
\pi_n:=\zeta^{\varphi^{-n}}(\zeta_{2^{n+2}}-1),
\end{align}
and define $d_n\in \what{E}(\fm_n)$ by
\begin{align}
d_n:=\exp_{\what{E}}\circ \log_{\cG_n}(\varepsilon_n[+]_{\cG_n} \pi_n).
\end{align}
\end{dfn}
\begin{lem}\label{injective_log}
For every $n\geq -2$, $\what{E}(\fm_n)$ is
torsion-free. In particular,
\begin{align}
\log_{\what{E}}\colon \what{E}(\fm_n)\longrightarrow k_n
\end{align}
is injective.
\end{lem}

\begin{proof}
The argument is the same as in the odd-prime case
(cf.~\cite[Proposition 3.1]{KO18}), with a slight modification at
$p=2$. Since $\what{E}$ has height $2$, every non-zero point of
$\what{E}[2]$ generates a totally ramified extension of
$F_{\fp}$ of degree $2^2-1=3$. On the other hand, the degrees $[k_n:F_{\fp}]$ are powers of $2$. Hence $\what{E}(\fm_n)[2]=0$. 
Since $\widehat E(\mathfrak m_n)$ is a pro-$2$ group,
it follows that it is torsion-free.
The last assertion follows
since the kernel of the formal logarithm consists of torsion
points.
\end{proof}

\begin{prop}\label{Prop:trace_relation}
The set of local points $\{d_n\}_{n\geq-2}$ forms a system satisfying
\begin{align}
\begin{cases}\tr_{n/n-1}(d_n)=
a_2d_{n-1}-d_{n-2}&(n\geq0),\\
d_{-1}=(a_2-2\varphi)d_{-2}.
\end{cases}
\end{align}
\end{prop}
\begin{proof}
Since $\log_{\widehat{E}}$ is injective by Lemma \ref{injective_log}, it suffices to verify the
assertions after applying $\log_{\widehat{E}}$. 
By the definitions of $d_n$ and $f$, we have, for $n\geq-2$,
\begin{align}\label{eq:log_dn}
\log_{\what E}(d_n)
&=\eta_n+f^{\varphi^{-n}}(\pi_n)\notag\\
&=\eta_n+\sum_{j\geq0}x_j
\bigl(\zeta^{\varphi^{-n}}\bigr)^{2^j}
\bigl(\zeta_{2^{n+2}}^{2^j}-1\bigr).
\end{align}
For $n\geq 1$, the trace relation follows from a computation similar to that in \cite[Proof of Proposition 3.23]{Kat21}. We therefore only treat the bottom layers $n=0$ and $n=-1$. 
By \eqref{eq:log_dn}, we have
\begin{align}
&\log_{\what E}(d_0)=\eta_0+\pi_0+x_1\pi_{-1},\\
&\log_{\what E}(d_{-1})=\eta_{-1}+\pi_{-1},\\
&\log_{\what E}(d_{-2})=\eta_{-2}.
\end{align}
Note that
\begin{align}
  \pi_{-1}
  =\zeta^{\varphi^{-1}}(\zeta_2-1)
  =-2\zeta^{\varphi^{-1}}\neq 0.
\end{align}

For $n=0$, since $\eta_0,\pi_{-1}\in k_{-1}$ and
$\Tr_{0/-1}(\pi_0)=-2\zeta$, we obtain
\begin{align}
\Tr_{0/-1}\bigl(\log_{\what E}(d_0)\bigr)
&=2\eta_0-2\zeta+2x_1\pi_{-1}\\
&=a_2\eta_0^\varphi-\eta_0+a_2\pi_{-1}\\
&=a_2(\eta_{-1}+\pi_{-1})-\eta_{-2}\\
&=a_2\log_{\what E}(d_{-1})
-\log_{\what E}(d_{-2}),
\end{align}
where we used $2x_1=a_2$,
$(3-a_2\varphi)\eta_0=2\zeta$,
$\eta_0^\varphi=\eta_{-1}$, and
$\eta_0=\eta_{-2}$. Hence
\begin{align}
\Tr_{0/-1}(d_0)=a_2d_{-1}-d_{-2}.
\end{align}
Since $(3-a_2\varphi)\eta_{-1}
=2\zeta^{\varphi^{-1}}$ and
$\eta_{-1}=\eta_{-2}^\varphi$, we have
\begin{align}
\log_{\what E}(d_{-1})
&=\eta_{-1}+\pi_{-1}\\
&=(a_2\varphi-2)\eta_{-1}\\
&=(a_2-2\varphi)\eta_{-2}\\
&=(a_2-2\varphi)\log_{\what E}(d_{-2}).
\end{align}
Hence
\begin{align}
d_{-1}=(a_2-2\varphi)d_{-2}.
\end{align}
\end{proof}

\begin{prop}\label{Prop:generate_formal}
We have
\begin{align}
\what{E}(\fm_n)=\begin{cases}
\langle d_n,d_{n-1} \rangle_{\Z_2[G_n]} &(n\geq1),\\
\langle d_0\rangle_{\Z_2[G_0]}&(n=0),\\
\langle d_{-2} \rangle_{\Z_2[G_{-2}]}&(n=-1,-2).
\end{cases}
\end{align}
\end{prop}
\begin{proof}

We first justify the argument
of \cite[Lemma 3.9, Proposition 3.10]{KO18}
and \cite[Lemma 3.18, Proposition 3.24]{Kat21}
in our setting at $p=2$.
Since $\cO_n=\cO[\zeta_{2^{n+2}}]$ and
$\cO=\langle\zeta\rangle_{\Z_2[\gal(F_{\fp}/\Q_2)]}$,
the argument of \cite[Lemma 3.9]{KO18} gives
\begin{align}
\fm_n/\fm_{n-1}
=\langle\pi_n\rangle_{\Z_2[G_n]}.
\end{align}
For $x\in\fm_n$, write
$x=a+b(\zeta_{2^{n+2}}-1)$ with
$a\in\fm_{n-1}$ and $b\in\cO_{n-1}$.
Since
\begin{align}
(\zeta_{2^{n+2}}-1)^2
=(\zeta_{2^{n+2}}^2-1)-2(\zeta_{2^{n+2}}-1),
\end{align}
we have
\begin{align}
x^{2^r}\in 2^r\fm_n+\fm_{n-1}
\qquad(r\geq1).
\end{align}
As $\log_{\what E}'(T)\in\Z_2\llbracket T\rrbracket$,
the argument of \cite[Lemma 3.18]{Kat21}
therefore gives
\begin{align}
\log_{\what E}(\what E(\fm_n))
\subseteq\fm_n+k_{n-1}.
\end{align}
Since $\log_{\what E}$ is injective
by Lemma \ref{injective_log}
and commutes with the action of
$\gal(k_n/k_{n-1})$, it induces
an injection
\begin{align}
\what E(\fm_n)/\what E(\fm_{n-1})
&\hookrightarrow
(\fm_n+k_{n-1})/k_{n-1}\notag\\
&\cong\fm_n/\fm_{n-1}.
\end{align}
By \eqref{eq:log_dn}, the image of $d_n$
is the class of $\pi_n$, which generates
$\fm_n/\fm_{n-1}$.
Hence the above injection is surjective, and
\begin{align}\label{eq:quotient_gen}
\what{E}(\fm_n)/\what{E}(\fm_{n-1})=\langle d_n\rangle_{\Z_2[G_n]}.
\end{align}
Since $(3-a_2\varphi)\eta_{-2}=2\zeta$, $3-a_2\varphi\in\Z_2[G_{-2}]^\times$, and
$\fm_{-2}=\fm_{-1}=\fm=\langle2\zeta\rangle_{\Z_2[G_{-2}]}$, we obtain
\begin{align}
\what{E}(\fm_{-1})=\what{E}(\fm_{-2})=\langle d_{-2}\rangle_{\Z_2[G_{-2}]}.
\end{align}
By Proposition \ref{Prop:trace_relation}, we see that
\begin{align}
\tr_{0/-1}(d_0)=\bigl(a_2(a_2-2\varphi)-1\bigr)d_{-2}.
\end{align}
Hence $d_{-2}\in\langle d_0\rangle_{\Z_2[G_0]}$, and \eqref{eq:quotient_gen} yields
\begin{align}
\what{E}(\fm_0)=\langle d_0\rangle_{\Z_2[G_0]}.
\end{align}
The case $n\geq1$ now follows from \eqref{eq:quotient_gen} and
\begin{align}
\tr_{n/n-1}(d_n)=a_2d_{n-1}-d_{n-2}
\end{align}
by induction, as in \cite[Corollary 3.11]{KO18}.
\end{proof}

\subsection{The chromatic Coleman map}\label{subsec:Coleman}

For simplicity, put $R_n:=\Z_2[G_n]$ for each $n\geq-2$ and $R:=\Z_2\llbracket \gal(k_{\infty}/\Q_2)\rrbracket$. Since $F_{\fp}/\Q_2$ is unramified and $\Q_2(\zeta_{2^{n+2}})/\Q_2$ is totally ramified, we have a natural decomposition 
\begin{align}
G_n\cong G_0\times (\Gamma/\Gamma^{2^n}).
\end{align}
For each $n\geq0$, we identify $R_0$ with its image in $R_n$ via the natural embedding
\begin{align}\label{R_0-R_n}
R_0\hookrightarrow R_0[\Gamma/\Gamma^{2^n}]\cong R_n
\end{align}
and use the natural homomorphism
\begin{align}\label{Lambda_R_n}
\Lambda=\Z_2\llbracket \Gamma \rrbracket \twoheadrightarrow \Z_2[\Gamma/\Gamma^{2^n}]\hookrightarrow R_n.
\end{align}
Following the convention of \cite[\S4]{Kat21}, we regard elements of $R_n^{\oplus2}$ as row vectors, and matrices act on  them from the right. For a matrix $M$, we denote its transpose by $M^{\mathsf T}$. For $n\geq1$, we put
\begin{align}
N_n:=1+\gamma^{2^{n-1}},\quad A_n:=\begin{pmatrix}0&-N_n\\ 1&a_2\end{pmatrix}\in M_2(\Lambda),\quad \cB_n:=A_n\cdots A_1,
\end{align}
and let $\cB_0$ be the identity matrix in $M_2(\Lambda)$.
Proposition \ref{Prop:trace_relation} implies that, for $n\geq1$,
\begin{align}
(d_n,d_{n-1})A_n=(d_{n-1},d_{n-2}).
\end{align}
\begin{dfn}
We define, for $n\geq0$,
\begin{align}
\cM_n:=R_n^{\oplus2}/R_n^{\oplus2}(\gamma-1)\cB_n^{\mathsf T},
\end{align}
and $t_n:\cM_n\to \what{E}(\fm_n)$ by $t_n(\overline{x,y}):=xd_n+yd_{n-1}$. Moreover, we define $s_n:R_0\to \cM_n$ by 
\begin{align}
s_n(z):=\overline{z(-w,1)\cB_n^{\mathsf T}}
\end{align} 
for $z\in R_0$. Here, 
\begin{align}
w:=(a_2-2\varphi)(a_2(a_2-2\varphi)-1)^{-1}\sum_{\sigma\in\gal(k_0/F_{\fp})}\sigma\in R_0.
\end{align}
In the definition of $s_n$, the product $z(-w,1)\cB_n^{\mathsf{T}}$ is computed in $R_n^{\oplus 2}$ via \eqref{R_0-R_n} and \eqref{Lambda_R_n}, and its image in the quotient $\cM_n$ is taken.
\end{dfn}

\begin{lem}
For each $n\geq0$, the sequence
\begin{align}\label{eq:fundamental_seq}
\xymatrix{
0\ar[r] & R_0\ar^{s_n}[r] & \cM_n\ar^{t_n\ }[r] & \what{E}(\fm_n)\ar[r] & 0
}
\end{align}
is exact.
\end{lem}
\begin{proof}
By Proposition \ref{Prop:generate_formal}, the map $t_n$ is surjective. 
Moreover, iterating
\begin{align}
(d_j,d_{j-1})A_j=(d_{j-1},d_{j-2})
\end{align}
for $1\leq j\leq n$ gives
\begin{align}
(d_n,d_{n-1})\cB_n=(d_0,d_{-1}).
\end{align}
Since $d_{-1}=wd_0$, it follows that, for $z\in R_0$,
\begin{align}
(t_n\circ s_n)(z)
&=z(d_n,d_{n-1})\cB_n
\begin{pmatrix}-w\\1\end{pmatrix}\notag\\
&=z(-wd_0+d_{-1})=0.
\end{align}
By the argument of \cite[Lemma 4.9 (1)]{Kat21},
which uses
\begin{align}
\det((\gamma-1)\cB_n^{\mathsf T})
=(\gamma-1)(\gamma^{2^n}-1),
\end{align}
the module $\cM_n$ is free over $\Z_2$, with
\begin{align}
\operatorname{rank}_{\Z_2}\cM_n
=[k_0:\Q_2](2^n+1).
\end{align}

We next show that $s_n$ is injective with
$\Z_2$-torsion-free cokernel.
It suffices to show that the reduction
\begin{align}
\overline{s}_n:R_0/2R_0
\longrightarrow\cM_n/2\cM_n
\end{align}
is injective.
For $n=0$, we have $\cM_0=R_0^{\oplus2}$
and $s_0(z)=z(-w,1)$, so the assertion is clear.
Suppose $n\geq1$ and put $T:=\gamma-1$.
Since $a_2\in2\Z_2$, we have $w\in2R_0$.
Moreover,
\begin{align}
R_n/2R_n
\cong(R_0/2R_0)[T]/(T^{2^n}),\quad 
A_j^{\mathsf T}
\equiv
\begin{pmatrix}
0&1\\
T^{2^{j-1}}&0
\end{pmatrix}
\pmod2.
\label{eq:An_mod2}
\end{align}
Since $w\equiv0\pmod2$, we have
\begin{align}
s_n(z)\equiv z(0,1)\cB_n^{\mathsf T}\pmod2.
\end{align}
Moreover, for $j\geq1$,
\begin{align}
A_{2j-1}^{\mathsf T}A_{2j}^{\mathsf T}
\equiv
\begin{pmatrix}
T^{2^{2j-1}}&0\\
0&T^{2^{2j-2}}
\end{pmatrix}
\pmod2.
\label{eq:An_pair_mod2}
\end{align}

Put $\overline{R_n}:=R_n/2R_n$. By \eqref{eq:An_mod2} and \eqref{eq:An_pair_mod2}, for some nonnegative integers $u,e$ with $e<2^n$, we have
\begin{align}
\cB_n^{\mathsf T}\equiv
\begin{cases}
\begin{pmatrix}
0&T^u\\
T^e&0
\end{pmatrix}
&\text{if $n$ is odd},\\[6pt]
\begin{pmatrix}
T^u&0\\
0&T^e
\end{pmatrix}
&\text{if $n$ is even}
\end{cases}
\pmod 2.
\end{align}
Since
\begin{align}
\cM_n/2\cM_n
\cong \overline{R_n}^{\oplus2}/
\overline{R_n}^{\oplus2}T\overline{\cB_n^{\mathsf T}},
\end{align}
we obtain
\begin{align}
\cM_n/2\cM_n\cong
\begin{cases}
\overline{R_n}/(T^{e+1})\oplus \overline{R_n}/(T^{u+1})
&\text{if $n$ is odd},\\
\overline{R_n}/(T^{u+1})\oplus \overline{R_n}/(T^{e+1})
&\text{if $n$ is even}.
\end{cases}
\end{align}
Under these identifications, the reduction
of $s_n$ is given by
\begin{align}
\overline{s}_n(z)=
\begin{cases}
(zT^e,0)&\text{if $n$ is odd},\\
(0,zT^e)&\text{if $n$ is even}.
\end{cases}
\end{align}
Since $zT^e=0$ in $\overline{R_n}/(T^{e+1})$
implies $z=0$, the map
$\overline{s}_n$ is injective.
Since $R_0$ and $\cM_n$ are free
$\Z_2$-modules, it follows that $s_n$
is injective with $\Z_2$-torsion-free cokernel.

Choose $r$ sufficiently large so that the formal logarithm induces an isomorphism
\begin{align}
\log_{\what{E}}:\what{E}(\fm_n^r)\overset{\sim}{\longrightarrow} \fm_n^r.
\end{align}
Since $\what{E}(\fm_n^r)$ has finite index in $\what{E}(\fm_n)$, the latter is a finitely generated $\Z_2$-module and 
\begin{align}
\rank_{\Z_2}\what{E}(\fm_n)=\rank_{\Z_2}\what{E}(\fm_n^r)=[k_n:\Q_2].
\end{align}
By Lemma \ref{injective_log}, $\what{E}(\fm_n)$ is $\Z_2$-torsion-free, and hence it is $\Z_2$-free. Since $[k_n:k_0]=2^n$ and 
\begin{align}
\rank_{\Z_2}R_0=[k_0:\Q_2],
\end{align}
we obtain
\begin{align}
\rank_{\Z_2}\cM_n=\rank_{\Z_2}R_0+\rank_{\Z_2}\what{E}(\fm_n).
\end{align}
Hence all the hypotheses of \cite[Lemma 4.15]{Kat21} are satisfied, and the sequence is exact.
\end{proof}

For a compact or discrete finitely generated $\Z_2$-module $X$, we put
\begin{align}
X^*:=\Hom_{\Z_2}(X,\Z_2).
\end{align}
Let $\iota$ denote the involution of $R_n$ induced by
$\sigma\mapsto\sigma^{-1}$, and let $\widetilde{\cB}_n$
denote the adjugate matrix of $\cB_n$. By \cite[Lemma 4.2 (2)]{Kat21}, we have
\begin{align}
\cM_n^*
\cong
\ker\left(
R_n^{\oplus2}
\xrightarrow{\times((\gamma-1)\cB_n)^\iota}
R_n^{\oplus2}
\right).
\end{align}
Moreover, since
\begin{align}
(\gamma-1)\cB_n\widetilde{\cB}_n
=(\gamma^{2^n}-1)I_2,
\end{align}
the argument of \cite[Lemma 4.9 (2)]{Kat21}
gives an isomorphism
\begin{align}
\cM_n^*
\cong
R_n^{\oplus2}/
R_n^{\oplus2}((\gamma-1)\cB_n)^\iota.
\end{align}
Thus, taking the $\Z_2$-linear dual of $t_n$,
we obtain an $R_n$-homomorphism
\begin{align}
\col'_n:
\what{E}(\fm_n)^*
\longrightarrow
R_n^{\oplus2}/
R_n^{\oplus2}((\gamma-1)\cB_n)^\iota.
\end{align}
By \cite[Lemma 4.11]{Kat21}, these maps are compatible
with the natural transition maps. Moreover, as in
\cite[Definition 4.12]{Kat21}, we have
\begin{align}
\varprojlim_n
R_n^{\oplus2}/
R_n^{\oplus2}((\gamma-1)\cB_n)^\iota
\cong R^{\oplus2}.
\end{align}
Therefore, taking the projective limit, we obtain
an $R$-homomorphism
\begin{align}
\col':
\varprojlim_n\what{E}(\fm_n)^*
\longrightarrow R^{\oplus2}.
\end{align}

\begin{dfn}[cf. {\cite[Definition 4.12]{Kat21}}]\label{Definition: Coleman_map}
We define the \emph{Coleman map}
\begin{align}
\col:=(\col^\sharp,\col^\flat):=\col'U:\varprojlim_n \what{E}(\fm_n)^*
\longrightarrow R^{\oplus2},
\end{align}
where $U:=\begin{pmatrix}-a_2&-1\\1&0\end{pmatrix}$. The components
\begin{align}
\col^\sharp,\col^\flat:
\varprojlim_n \what{E}(\fm_n)^*
\longrightarrow R
\end{align}
are called the \emph{chromatic Coleman maps}.
\end{dfn}

Let $\cT:=T_2E$ be the $2$-adic Tate module of $E$.
For each $n\geq 0$, let
\begin{align}
(-,-)_n:
\what{E}(\fm_n)\times H^1(k_n,\cT)\longrightarrow \Z_2
\end{align}
be the local Tate pairing. It induces a natural surjective homomorphism
\begin{align}
H^1(k_n,\cT)\twoheadrightarrow \what{E}(\fm_n)^*.
\end{align}
By the compatibility of the local Tate pairing
with restriction and corestriction, these surjections
are compatible with the natural transition maps.
By compactness, the induced map on projective limits
is also surjective. Thus, we also regard the chromatic
Coleman maps as
\begin{align}
\col^\sharp,\col^\flat:
H^1_{\mathrm{Iw}}(k_\infty,\cT)\longrightarrow R,
\end{align}
where
\begin{align}
H^1_{\mathrm{Iw}}(k_\infty,\cT)
:=
\varprojlim_n H^1(k_n,\cT).
\end{align}
For $M=\Q_2$ or the unramified quadratic extension of $\Q_2$,
put $M_n:=M(\zeta_{2^{n+2}})$ for $n\geq-2$ and $M_{\infty}:=\bigcup_{n\geq-2}M_n$, and write
\begin{align}
H^1_{\iw,M}(\cT)
:=
\varprojlim_n H^1(M_n,\cT).
\end{align}
In particular,
$H^1_{\iw,F_{\fp}}(\cT)=H^1_{\mathrm{Iw}}(k_\infty,\cT)$.

By local Tate duality and the Weil pairing, we have a perfect pairing
\begin{align}
H^1(k_\infty,E[2^{\infty}])\times
H^1_{\mathrm{Iw}}(k_\infty,\cT)
\longrightarrow \Q_2/\Z_2.
\end{align}

\begin{dfn}\label{Definition: Local_condition}
For $\star\in\{\sharp,\flat\}$, we define
\begin{align}
E^\star_\infty
:=
\parenth{\Ker\col^\star}^{\perp}
\subset H^1(k_\infty,E[2^{\infty}]),
\end{align}
where the orthogonal complement is taken with respect to the above
local Tate pairing.
\end{dfn}
By Pontryagin duality, we have a natural isomorphism
\begin{align}
(E^\star_\infty)^\vee
\cong
H^1_{\mathrm{Iw}}(k_\infty,\cT)/\Ker\col^\star
\cong
\Imag\col^\star.
\end{align}

\subsection{Doubly-signed Selmer groups}\label{subsec:doubly_Selmer}

In this subsection, we define $2$-adic doubly-signed Selmer groups following \cite{Spr21}, \cite{LL22}, \cite{Lim22}, and \cite{FM25}. Let $F$ be a real quadratic field in which $2$ is unramified and let $E$ be an elliptic curve defined over $\Q$ with good supersingular reduction at 2.

For $n\geq0$, put $F_n:=F(\zeta_{2^{n+2}})$ and
$F_\infty:=\bigcup_{n\geq0}F_n$, and let
$F^{\cyc}$ and $\Q^{\cyc}$ be the cyclotomic
$\Z_2$-extensions of $F$ and $\Q$, respectively.
 Let $\Lambda:=\Z_2\llbracket \gal(\Q^{\cyc}/\Q)\rrbracket$. Since $F_{\fp}/\Q_2$ is unramified,
restriction induces a natural isomorphism
\begin{align}
\Gamma=\gal(k_\infty/k_0)
\cong \gal(\Q^{\cyc}/\Q).
\end{align}
We identify these two Galois groups throughout the paper.

For each prime $\fp\mid2$ of $F$, let $F_{\infty,\fp}$ denote the union of the completions of $F_n$ at the unique primes above $\fp$.
Then
\begin{align}
F_{\infty,\fp}=F_\fp(\zeta_{2^\infty})=k_\infty.
\end{align}
\begin{lem}\label{lem:2-primary_torsion}
For every prime $\fp\mid2$ of $F$ and every $n\geq-2$, we have
\begin{align}
E(k_n)[2^\infty]=\{O\}.
\end{align}
Consequently,
\begin{align}
E(F_\infty)[2^\infty]=\{O\}.
\end{align}
\end{lem}

\begin{proof}
Let $P\in E(k_n)[2^\infty]$, and let $\kappa_n$
be the residue field of $k_n$.
Then the reduction $\widetilde P$ belongs to
$\widetilde E(\kappa_n)[2^\infty]$.
Since $\widetilde E$ is supersingular, we have
\begin{align}
\widetilde E(\kappa_n)[2^\infty]=\{\widetilde{O}\}.
\end{align}
Thus $\widetilde{P}=\widetilde{O}$, so
$P\in\widehat E(\fm_n)[2^\infty]$.
Lemma \ref{injective_log} gives $P=O$.
Therefore $E(k_n)[2^\infty]=\{O\}$.
\end{proof}
Applying the construction in Definition \ref{Definition: Local_condition} to each $\fp\mid2$, for
$\star\in\{\sharp,\flat\}$, we write
\begin{align}
E^\star_{\infty,\fp}
:=
\parenth{\Ker\col^\star_\fp}^{\perp}
\subset
H^1(F_{\infty,\fp},E[2^\infty]).
\end{align}
Thus, by Pontryagin duality,
\begin{align}
\parenth{E^\star_{\infty,\fp}}^\vee
\cong
H^1_{\mathrm{Iw}}(F_{\infty,\fp},T)/
\Ker\col^\star_\fp
\cong
\Imag\col^\star_\fp.
\end{align}

Let $\Sigma_\Q$ be a finite set of places of $\Q$
containing the archimedean place, the prime $2$,
the primes of bad reduction of $E$, and the primes
ramified in $F/\Q$.
Let $\Sigma_F$ be the set of places of $F$ lying
above $\Sigma_\Q$.
For $M=\Q$ or $F$, let $M_\Sigma$ denote the maximal
algebraic extension of $M$ unramified outside $\Sigma_M$. For an algebraic extension $M'/M$, we denote by
$\Sigma(M')$ the set of places of $M'$ lying above
$\Sigma_M$, and put
\begin{align}
&\Sigma_2(M') := \{w\in\Sigma(M') \mid w\mid 2\},\\
&\Sigma_\infty(M') := \{w\in\Sigma(M') \mid
w \text{ is archimedean}\}.
\end{align}
For $w\in\Sigma(M')$, let $M'_w$ denote the completion
of $M'$ at $w$ when $M'/M$ is finite.
If $M'/M$ is infinite, we set
$M'_w:=\bigcup_L L_{w_L}$, where $L$ runs through
the finite intermediate fields of $M'/M$ and $w_L$
is the restriction of $w$ to $L$.

For $M=\Q$ or $F$, put
$M_\infty:=M(\zeta_{2^\infty})$ and
$\Delta_M:=\gal(M_\infty/M^{\cyc})$.
Restriction induces an isomorphism
\begin{align}
\Delta_M\cong\gal(M(\zeta_4)/M)
\cong\Z/2\Z.
\end{align}
We identify $\Delta_\Q$ and $\Delta_F$ via
restriction and simply write $\Delta$.

\begin{dfn}
Let $M=\Q$ or $F$. For
\begin{align}
\vec{s}
=
(s_\fp)_{\fp\in\Sigma_2(M)}
\in
\prod_{\fp\in\Sigma_2(M)}
\{\sharp,\flat\},
\end{align}
we define the \emph{(at most) doubly-signed Selmer group} by
\begin{align}
\begin{aligned}
\sel^{\vec{s}}(E/M_\infty)
:={}&
\Ker\Biggl(
H^1(M_\Sigma/M_\infty,E[2^\infty])
\\
&\longrightarrow
\prod_{\fw\in\Sigma_2(M_\infty)}
\frac{H^1(M_{\infty,\fw},E[2^\infty])}
     {E^{s_\fw}_{\infty,\fw}}
\times
\prod_{\fw\in\Sigma(M_\infty)\setminus\Sigma_2(M_\infty)}
H^1(M_{\infty,\fw},E[2^\infty])
\Biggr).
\end{aligned}
\end{align}
Here, for $\fw\in \Sigma_2(M_{\infty})$ lying above $\fp\in \Sigma_2(M)$, $s_{\fw}:=s_{\fp}$. As in the classical case, the above definition is independent of the choice of
$\Sigma$.
\end{dfn}

For $M=\Q$ or $F$ and $n\geq 0$, put
$M_n^{\cyc}:=M(\zeta_{2^{n+2}})\cap M^{\cyc}$.

\begin{dfn}\label{Definition_Selmer_chromatic_cyclotomic}
Let $M=\Q$ or $F$.
We define
\begin{align}
\mathscr H_\infty(M,E)
:=
\varinjlim_{n}
\prod_{v\in\Sigma_\infty(M_n^{\cyc})}
H^1((M_n^{\cyc})_v,E[2^\infty]).
\end{align}
We then define the chromatic Selmer group over
$M^{\cyc}$ by
\begin{align}
\sel^{\vec{s}}(E/M^{\cyc})
&:=
\Ker\left(
\sel^{\vec{s}}(E/M_\infty)^{\Delta_M}
\xrightarrow{\operatorname{loc}_\infty}
\mathscr H_\infty(M,E)
\right)\\
&=
\Ker\Biggl(
H^1(M_\Sigma/M^{\cyc},E[2^\infty])
\longrightarrow
\left(\prod_{\fw\in\Sigma_2(M_\infty)}
\frac{H^1(M_{\infty,\fw},E[2^\infty])}
     {E^{s_\fw}_{\infty,\fw}}\right)^{\Delta_M}\\
&\times
\left(\prod_{\fw\in\Sigma(M_\infty)\setminus\Sigma_2(M_\infty)}
H^1(M_{\infty,\fw},E[2^\infty])
\right)^{\Delta_M}
\times \mathscr{H}_\infty(M,E)\Biggr),
\end{align}
where $\operatorname{loc}_\infty$ is induced by
the inflation--restriction isomorphism and
localization at the archimedean places. 
\end{dfn}

We denote the Pontryagin dual of a module $X$ by $X^{\vee}$, and put 
\begin{align}
\cX^{\vec{s}}(E/M^{\cyc}):=\sel^{\vec{s}}(E/M^{\cyc})^{\vee}
\end{align}
for $M=\Q$ or $F$.

\begin{dfn}
Let $M=\Q$ or $F$. For comparison, we define the chromatic
Selmer group with relaxed local conditions at the archimedean
places by
\begin{align}
\sel_{\mathrm{rel}}^{\vec{s}}(E/M^{\cyc}):=\sel^{\vec{s}}(E/M_\infty)^{\Delta_M},\quad \cX_{\mathrm{rel}}^{\vec{s}}(E/M^{\cyc}):=\sel_{\mathrm{rel}}^{\vec{s}}(E/M^{\cyc})^{\vee}.
\end{align}
Since 
\begin{align}
H^1(\R,E[2^{\infty}])\cong
\begin{cases}
\Z/2\Z&\text{if }\Delta_E>0,\\
0&\text{otherwise},
\end{cases}
\end{align}
$\sel_{\mathrm{rel}}^{\vec{s}}(E/M^{\cyc})\cong \sel^{\vec{s}}(E/M^{\cyc})$ if $\Delta_E<0$.
\end{dfn}

\begin{lem}\label{Lem:archimedean_comparison}
There is an exact sequence of $\Lambda$-modules
\begin{align}\label{eq:archimedean_comparison}
0\longrightarrow
\Imag(\loc_\infty)^\vee
\longrightarrow \cX_{\mathrm{rel}}^{\vec{s}}(E/M^{\cyc})
\longrightarrow \cX^{\vec{s}}(E/M^{\cyc})
\longrightarrow 0.
\end{align}
Moreover, $\Imag(\loc_\infty)^\vee$ is annihilated by $2$.
If the two modules on the right of \eqref{eq:archimedean_comparison}
are $\Lambda$-torsion, they have the same $\lambda$-invariant.
Inverting $2$ in the exact sequence gives a natural isomorphism
\begin{align}\label{eq:archimedean_rational_comparison}
\cX_{\mathrm{rel}}^{\vec{s}}(E/M^{\cyc})[1/2]
\cong \cX^{\vec{s}}(E/M^{\cyc})[1/2].
\end{align}
The corresponding $\mu$-invariants need not coincide.
\end{lem}

\begin{proof}
The exact sequence \eqref{eq:archimedean_comparison} follows from Definition \ref{Definition_Selmer_chromatic_cyclotomic} by Pontryagin duality.
Since $H^1(\R,E[2^\infty])$ is annihilated by $2$, the same holds for $\Imag(\loc_\infty)^\vee$.
Moreover, $\Imag(\loc_\infty)^\vee$ is finitely generated over $\Lambda$, since it is a quotient of $\mathscr{H}_\infty(M,E)^\vee$.
Thus its characteristic ideal is a power of $(2)$, and hence its $\lambda$-invariant is zero.
The equality of the $\lambda$-invariants follows from the multiplicativity of characteristic ideals in \eqref{eq:archimedean_comparison}.
Therefore inverting $2$ in \eqref{eq:archimedean_comparison} gives the isomorphism \eqref{eq:archimedean_rational_comparison}.
\end{proof}

The following conjecture extends
\cite[Conjecture 7.15]{Spr12} to $p=2$
and real quadratic fields.

\begin{conj}\label{Conj:Lambda_cotorsion}
Let $M=\Q$ or $F$. Then $\cX^{\vec{s}}(E/M^{\cyc})$ is $\Lambda$-torsion.
\end{conj}

In the arguments below concerning $E/\Q$, we take $F_{\fp}=\Q_2$
in the constructions of Subsections \ref{subsec:local_points} and
\ref{subsec:Coleman}. In particular,
$R=\Z_2\llbracket \gal(\Q_\infty/\Q)\rrbracket$ via the natural
identification of the local and global Galois groups.
We write $H^1_{\mathrm{Iw}}(\cT):=H^1_{\iw,\Q_2}(\cT)$.
For $\star\in\{\sharp,\flat\}$, put
$\cX^\star:=\sel^\star(E/\Q_\infty)^\vee$.
Let $\delta$ denote the non-trivial element of $\Delta$ and put
$\Lambda_{\Q_2}:=\Lambda[1/2]$. Since $2$ is invertible in
$R[1/2]$, we have the idempotent
\begin{align}
e_\Delta:=\frac{1+\delta}{2},
\end{align}
and we identify $e_\Delta R[1/2]$ with $\Lambda_{\Q_2}$
via the natural projection.

Choose the generator $\gamma$ of $\Gamma$ so that
$\chi_{\mathrm{cyc}}(\gamma)=5^{-1}$, where
$\chi_{\mathrm{cyc}}$ is the cyclotomic character, and put
$T:=\gamma-1$. In the comparison with Sprung's convention,
we use the variable $X:=\gamma^{-1}-1=(1+T)^{-1}-1$.

For $n\geq0$, put $L_n:=\Q(\zeta_{2^{n+2}})$, so that
$\Q_\infty=\bigcup_{n\geq0}L_n$. Let
$j_{L_n}:\operatorname{Spec}L_n
\to\operatorname{Spec}\cO_{L_n}[1/2]$
be the canonical morphism induced by the inclusion
$\cO_{L_n}[1/2]\hookrightarrow L_n$.
We distinguish the global Iwasawa cohomology groups
\begin{align}
\mathbb H^i_\Sigma(\cT):=\varprojlim_n H^i(\Q_\Sigma/L_n,\cT),\quad \mathbb H^i_{\mathrm{Kato}}(\cT):=\varprojlim_n H^i(\operatorname{Spec}\cO_{L_n}[1/2],j_{L_n,*}\cT),
\end{align}
where the transition maps are corestrictions.
The Leray spectral sequence for the open immersion
$\operatorname{Spec}\cO_{L_n}[1/\Sigma]\hookrightarrow
\operatorname{Spec}\cO_{L_n}[1/2]$
gives the following exact sequence
(cf. \cite[\S6]{Kur02}):
\[
0\to H^1(\operatorname{Spec}\cO_{L_n}[1/2],
j_{L_n,*}\cT)
\to H^1(\Q_\Sigma/L_n,\cT)
\to \bigoplus_{\substack{v\in\Sigma(L_n)\\v\nmid 2,\infty}}
H^0(\kappa(v),H^1(L_{n,v}^{\mathrm{nr}},\cT)).
\]
Here, $\kappa(v)$ denotes the residue field of $L_{n,v}$,
and $L_{n,v}^{\mathrm{nr}}$ the maximal unramified
extension of $L_{n,v}$.
As in \cite[\S 6, p.~217]{Kur02}, the inverse limit of
the rightmost term vanishes. The same argument applies
at $p=2$, since the decomposition groups at odd primes
are open in $\gal(\Q_\infty/\Q)$.
Consequently,
$\mathbb H^1_{\mathrm{Kato}}(\cT)
\cong\mathbb H^1_\Sigma(\cT)$.
Thus we write
\begin{align}
\mathbb H^1(\cT)
:=\mathbb H^1_\Sigma(\cT)
\cong\mathbb H^1_{\mathrm{Kato}}(\cT),
\qquad
\mathbb H^2(\cT):=\mathbb H^2_{\mathrm{Kato}}(\cT).
\end{align}
For $\star\in\{\sharp,\flat\}$, let
$L_{2,\mathrm{Kato}}^\star(E,T)\in\Lambda[1/2]$
denote the chromatic $2$-adic $L$-function defined by
Sprung in \cite[Definition 6.1]{Spr12} using Kato's zeta
element and the rational Coleman map.
Let $z=z^+\in\mathbb H^1(\cT)[1/2]$ be Kato's rational
zeta element associated with the trivial character of $\Delta$
as in \cite[Definition 6.1]{Spr12} and
\cite[Theorem 12.5 (1)]{Kat04}.
We denote its localization at $2$ by
$\loc(z)\in H^1_{\mathrm{Iw}}(\cT)[1/2]$.
For $\star\in\{\sharp,\flat\}$, put
\begin{align}
g_\star(T):=e_\Delta\col^\star(\loc(z))
\in\Lambda_{\Q_2}.
\end{align}

We first identify the Coleman map $\col$ constructed in
Subsection \ref{subsec:Coleman} following
\cite[Definition 4.12]{Kat21} with Sprung's rational
Coleman map used to define
$L_{2,\mathrm{Kato}}^\star(E,T)$
\cite[Definitions 5.9 and 6.1]{Spr12}.
The change of basis in Definition \ref{Definition: Coleman_map}
from $\col'$ to $\col$ is precisely the one chosen in
\cite[Remark 4.13]{Kat21} to agree with Sprung's sharp/flat
convention \cite[Definition 3.8]{Spr12}.
It remains to make this identification explicit for our
local points and our choice of the cyclotomic generator.

\begin{lem}\label{lem:Col_Kato_comparison}
For $\star\in\{\sharp,\flat\}$, we have
\begin{align}
g_\star(T)=L_{2,\mathrm{Kato}}^\star
\left(E,(1+T)^{-1}-1\right).
\end{align}
\end{lem}

\begin{proof}
We follow the logarithm-matrix comparison in
\cite[Proof of Theorem 6.9]{Kat21}, using the
$p=2$ results of \cite{Spr12}.

By \eqref{eq:log_dn}, our local points for
$F_{\fp}=\Q_2$ have the same normalization as
those of \cite[Theorem 2.2]{Spr12}.
In Sprung's construction, the finite-level Coleman
map is obtained by lifting the pair of group-ring-valued
local Tate pairings associated with $d_n$ and $d_{n-1}$
under right multiplication by the matrix $H_n(X)$
(cf. \cite[Propositions 3.9 and 5.3]{Spr12}).
In our construction, the corresponding pair of
local Tate pairings is obtained from the dual map
$t_n^*$ in Subsection \ref{subsec:Coleman}.
By \cite[Notation 3.7]{Spr12}, we have
\[
H_n(X)
=U^{-1}\prod_{j=1}^n
\begin{pmatrix}
a_2&\Phi_{2^j}(1+X)\\
-1&0
\end{pmatrix},
\]
where $U$ is the matrix appearing in
Definition \ref{Definition: Coleman_map}
and $\Phi_{2^j}$ is the $2^j$-th cyclotomic
polynomial.
Since
$\Phi_{2^j}(\gamma^{-1})
=1+\gamma^{-2^{j-1}}$,
we obtain
\begin{align}\label{relation_H_n-U}
H_n((1+T)^{-1}-1)
=U^{-1}\widetilde{\cB}_n^\iota.
\end{align}
Together with \cite[(4.7), Lemma 4.9 (2),
Definition 4.12 and Remark 4.13]{Kat21},
\eqref{relation_H_n-U} identifies the normalization of our
Coleman map with that of
\cite[Definition 5.9]{Spr12}.
In particular, applying the two Coleman maps
to Kato's zeta element gives, for $\star\in\{\sharp,\flat\}$,
\begin{align}
g_\star(T)
=L_{2,\mathrm{Kato}}^\star
\left(E,(1+T)^{-1}-1\right),
\end{align}
as required.
\end{proof}

This provides a $p=2$ counterpart of
\cite[Theorem 7.14]{Spr12}.

\begin{thm}\label{Thm:Lambda_torsion}
At least one of
\begin{align}
\cX^{\sharp}(E/\Q^{\cyc})
\quad\text{and}\quad
\cX^{\flat}(E/\Q^{\cyc})
\end{align}
is $\Lambda$-torsion. Moreover, if $a_2=0$,
then both of them are $\Lambda$-torsion.
\end{thm}

\begin{proof}
We modify the proof of \cite[Theorem 7.3 ii)]{Kob03} and
\cite[Theorem 7.14]{Spr12} by inverting $2$.
By Lemma \ref{Lem:archimedean_comparison},
\begin{align}
\cX^{\star}(E/\Q^{\cyc})[1/2]
\cong \cX^{\star}_{\mathrm{rel}}(E/\Q^{\cyc})[1/2]
\cong (\cX^\star)_\Delta[1/2].
\end{align}
Moreover,
\begin{align}\label{eq:delta_component}
(\cX^\star)_\Delta[1/2]=\cX^\star[1/2]/(\delta-1)\cX^\star[1/2]
\cong
e_\Delta\cX^\star[1/2].
\end{align}
Indeed,
\begin{align}
\cX^\star[1/2]
=
e_\Delta\cX^\star[1/2]
\oplus
(1-e_\Delta)\cX^\star[1/2],
\end{align}
and $\delta-1$ acts as $0$ on the first summand and as $-2$
on the second one.

By the Cassels--Poitou--Tate exact sequence,
we have the exact sequence
\begin{align}\label{seq:PT}
\xymatrix{
\mathbb H^1(\cT)(=\mathbb H^1_{\Sigma}(\cT))
\ar[r] &
\frac{H^1_{\mathrm{Iw}}(\cT)}{\Ker\col^\star}
\ar[r] &
\cX^{\star}\ar[r] &
\cX^0\ar[r] & 0,
}
\end{align}
where $\cX^0$ denotes the Pontryagin dual of the
fine Selmer group defined in
\cite[Definition 4.1]{Kur02}.
Tensoring \eqref{seq:PT} with $\Q_2$ and applying
$e_\Delta$, which is an exact functor over
$\Q_2[\Delta]$, gives
\begin{align}\label{seq_CPT}
\xymatrix{
e_\Delta\mathbb H^1(\cT)[1/2]
\ar[r]^{\iota^\star}
&
e_\Delta
\left(
\frac{H^1_{\mathrm{Iw}}(\cT)}
{\Ker\col^\star}
\right)[1/2]
\ar[r]
&
e_\Delta\cX^\star[1/2]
\ar[r]
&
e_\Delta\cX^0[1/2]
\longrightarrow0.
}
\end{align}

We justify the comparison of $\cX^0$ with Kato's
second Iwasawa cohomology. For each rational prime
$\ell\in\Sigma_\Q$, put
\begin{align}
\mathscr L_\ell^2(\cT)
:=\varprojlim_n\bigoplus_{v\mid\ell}
H^2((L_n)_v,\cT),
\end{align}
and put
$\Sigma_\Q^{(2)}:=\{\ell\in\Sigma_\Q\mid
\ell\neq2\text{ is a finite prime}\}$.
By the localization sequence and Poitou--Tate duality
(cf. \cite[\S 6]{Kur02}), we obtain
\begin{align}\label{eq:Kato_fine_localization}
0&\longrightarrow \mathbb H^2(\cT)
\longrightarrow \mathbb H^2_\Sigma(\cT)
\longrightarrow
\bigoplus_{\ell\in\Sigma_\Q^{(2)}}\mathscr L_\ell^2(\cT),
\\
0&\longrightarrow \cX^0
\longrightarrow \mathbb H^2_\Sigma(\cT)
\longrightarrow
\bigoplus_{\ell\in\Sigma_\Q^{(2)}\cup\{2\}}
\mathscr L_\ell^2(\cT).
\notag
\end{align}
Here, the non-$2$-adic local terms need not vanish;
the key point is that they occur in both maps in
\eqref{eq:Kato_fine_localization}.
By Lemma \ref{lem:2-primary_torsion} and local Tate
duality, we have
\begin{align}\label{2-component_vanish}
\mathscr L_2^2(\cT)
\cong E(\Q_2(\zeta_{2^\infty}))[2^\infty]^\vee=0.
\end{align}
Since the two maps in \eqref{eq:Kato_fine_localization}
have the same components at primes not dividing $2$,
their kernels coincide by \eqref{2-component_vanish}.
Consequently,
\begin{align}
e_\Delta\cX^0[1/2]
\cong e_\Delta\mathbb H^2(\cT)[1/2].
\end{align}
Since $\mathbb H^2(\cT)$ is $R$-torsion by
\cite[Theorem 12.4 (1)]{Kat04},
$e_\Delta\cX^0[1/2]$ is a torsion
$\Lambda_{\Q_2}$-module.
By \cite[Theorem 12.4 (2)]{Kat04},
$e_\Delta\mathbb H^1(\cT)[1/2]$ is a free
$\Lambda_{\Q_2}$-module of rank one.
Moreover,
\begin{align}\label{eq:ker_quotient_imge}
e_\Delta
\left(
\frac{H^1_{\mathrm{Iw}}(\cT)}
{\Ker\col^\star}
\right)[1/2]
\cong
e_\Delta\Imag(\col^\star)[1/2]
\subseteq
\Lambda_{\Q_2}.
\end{align}

By \cite[Proposition 6.14]{Spr12}, which relies on
Rohrlich's theorem \cite[p.~409, Theorem]{Roh84a}
on the non-vanishing of central Hasse--Weil
$L$-values $L(E,\chi,1)$ for all but finitely many
Dirichlet characters $\chi$ of $2$-power conductor,
at least one of $L_{2,\mathrm{Kato}}^\sharp(E,T)$ and
$L_{2,\mathrm{Kato}}^\flat(E,T)$ is non-zero.
Suppose that $a_2=0$.
By \cite[Remark 6.16]{Spr12}, we have
\begin{align}
\parenth{L_{2,\mathrm{Kato}}^\sharp(E,T),
L_{2,\mathrm{Kato}}^\flat(E,T)}
=
\parenth{-L_2^-(E,T),L_2^+(E,T)}
\end{align}
in Pollack's notation. Moreover, in the case $a_2=0$,
both $L_2^-(E,T)$ and $L_2^+(E,T)$ are non-zero by
\cite[Corollary 5.11]{Pol03}.

Fix $\star\in\{\sharp,\flat\}$ such that
$L_{2,\mathrm{Kato}}^\star(E,T)\neq0$.
By Lemma \ref{lem:Col_Kato_comparison},
we have $g_\star(T)\neq0$. Since
\[
\iota^\star(e_\Delta z)
=[e_\Delta\loc(z)],
\qquad
\col^\star(e_\Delta\loc(z))=g_\star(T)\neq0,
\]
the map $\iota^\star$ in \eqref{seq_CPT}
is non-zero.
In particular, the module on the right-hand side
of \eqref{eq:ker_quotient_imge} is non-zero,
and hence is a torsion-free $\Lambda_{\Q_2}$-module
of rank one.

Thus $\iota^\star$ is a non-zero homomorphism between
torsion-free $\Lambda_{\Q_2}$-modules of rank one.
It is therefore injective and its cokernel is a torsion
$\Lambda_{\Q_2}$-module.
Since $e_\Delta\cX^0[1/2]$ is a torsion
$\Lambda_{\Q_2}$-module, \eqref{seq_CPT} implies that
$e_\Delta\cX^\star[1/2]$ is also torsion over
$\Lambda_{\Q_2}$.
By Lemma \ref{Lem:archimedean_comparison} and
\eqref{eq:delta_component}, we have
\begin{align}\label{isom:cyc_Delta_component}
\cX^\star(E/\Q^{\cyc})[1/2]
&\cong
\cX_{\mathrm{rel}}^\star(E/\Q^{\cyc})[1/2]\\
&\cong
(\cX^\star)_\Delta[1/2]\\
&\cong
e_\Delta\cX^\star[1/2].
\end{align}
Consequently,
\begin{align}
\rank_\Lambda\cX^\star(E/\Q^{\cyc})
=
\rank_{\Lambda_{\Q_2}}
\parenth{e_\Delta\cX^\star[1/2]}=0.
\end{align}
Therefore, $\cX^\star(E/\Q^{\cyc})$
is $\Lambda$-torsion.
The first assertion follows from the above non-vanishing
results for $L_{2,\mathrm{Kato}}^\star(E,T)$.
If $a_2=0$, the above non-vanishing holds for both
$\star=\sharp$ and $\star=\flat$, and the same argument
proves that both chromatic Selmer groups are
$\Lambda$-cotorsion
(cf. \cite[Theorem 7.3 ii)]{Kob03}).
\end{proof}

In addition to the construction of
$L_{2,\mathrm{Kato}}^\star(E,T)$ described above,
Sprung constructed chromatic $2$-adic $L$-functions
by decomposing the classical $p$-adic $L$-functions of
Amice--V\'elu--Vi\v{s}ik \cite{AV75,Vis76}
using Mazur--Tate elements \cite[Theorem 2.14]{Spr17}.
In this paper, $L_2^\star(E,T)$ denotes the latter
construction, normalized using the Shimura periods of
\cite[Convention 2.1]{Ada26+}.
The variable $X$ in $L_2^\star(E,X)$ in
\cite[\S2]{Ada26+} corresponds to the generator whose image under the cyclotomic character is $5$.
Thus the two variables are related by
$X=(1+T)^{-1}-1$, as above.

\begin{dfn}
Let $E$ be an elliptic curve defined over $\Q$
with good supersingular reduction at $2$.
When $\cX^{\star}(E/\Q^{\cyc})$ is $\Lambda$-torsion
for $\star\in\{\sharp,\flat\}$, we can define its
Iwasawa $\mu$- and $\lambda$-invariants and denote
them by $\mu^{\star}(E/\Q)$ and
$\lambda^{\star}(E/\Q)$, respectively.
Similarly, for
$\vec{s}:=(s_{\fp})_{\fp\in\Sigma_2(F)}
\in \prod_{\fp\in \Sigma_2(F)}\{\sharp,\flat\}$,
if $\sel^{\vec{s}}(E/F^{\cyc})$ is $\Lambda$-cotorsion,
we define
\begin{align}
\mu^{\vec{s}}(E/F):=\mu\left(\sel^{\vec{s}}(E/F^{\cyc})\right),\quad \lambda^{\vec{s}}(E/F):=\lambda\left(\sel^{\vec{s}}(E/F^{\cyc})\right).
\end{align}
For $\star\in\{\sharp,\flat\}$ such that
$L_2^\star(E,T)\neq0$, we put
\begin{align}
\mu_{\an}^{\star}(E/\Q)
:=\mu\parenth{L_2^\star(E,T)},
\qquad
\lambda_{\an}^{\star}(E/\Q)
:=\lambda\parenth{L_2^\star(E,T)}.
\end{align}
When $\cX^\star(E/\Q^{\cyc})$ is $\Lambda$-torsion,
we also put
\begin{align}
\mu_{\alg}^{\star}(E/\Q):=\mu^\star(E/\Q),
\qquad
\lambda_{\alg}^{\star}(E/\Q):=\lambda^\star(E/\Q).
\end{align}
\end{dfn}
Throughout this paper, the Iwasawa invariants
without the subscripts $\an$ or $\alg$ are
understood to be algebraic.
When both analytic and algebraic invariants
are considered, we explicitly use the
subscripts $\an$ and $\alg$ to distinguish them.

\begin{lem}\label{lem:Col_analytic_comparison}
For $\star\in\{\sharp,\flat\}$, if
$L_2^\star(E,T)\neq0$, then $g_\star(T)\neq0$ and
\begin{align}
\lambda(g_\star)=\lambda(L_2^\star(E,T))
=\lambda_{\an}^\star(E/\Q).
\end{align}
\end{lem}

\begin{proof}
Let $\alpha,\beta$ be the roots of
$X^2-a_2X+2$, and write
\[
\mathcal L_{\alpha,\beta}(X):=
\begin{pmatrix}
\log_\alpha^\sharp(1+X)&
\log_\beta^\sharp(1+X)\\
\log_\alpha^\flat(1+X)&
\log_\beta^\flat(1+X)
\end{pmatrix}
\]
for the logarithm matrix of
\cite[Definition 6.8]{Spr12}.
For $\nu\in\{\alpha,\beta\}$, let
$L_2(E,\nu,X)$ denote the Amice--V\'elu--Vi\v{s}ik
$2$-adic $L$-function in the trivial
$\Delta$-component, normalized by the
N\'eron period $\Omega_E^+$, as in
\cite[Theorem 6.2]{Spr12}.
Here, $\Omega_E^+=\int_{E(\R)}|\omega_E|$,
where $\omega_E$ is a N\'eron differential of $E$
(cf. \cite[\S 3.1]{Ada26+}).
Combining Lemma \ref{lem:Col_Kato_comparison} with
Kato's explicit reciprocity law in the form of
\cite[Theorem 6.4]{Spr12} and
\cite[Main Theorem 6.12]{Spr12}, we obtain
\begin{align}
&\bigl(g_\sharp(T),g_\flat(T)\bigr)
\mathcal L_{\alpha,\beta}((1+T)^{-1}-1)\\
&\qquad=
\bigl(
L_2(E,\alpha,(1+T)^{-1}-1),
L_2(E,\beta,(1+T)^{-1}-1)
\bigr).
\end{align}
Let $f_E$ be the normalized newform of weight $2$
and level $N_E$ associated with $E$, and let
$\Omega_{f_E}^+$ denote its positive Shimura period,
chosen as in \cite[Convention 2.1]{Ada26+}.
On the other hand, \cite[Theorem 2.14 and
Corollaries 5.16--5.17]{Spr17}, together with
\cite[Theorem 2.11]{Ada26+}, gives
\begin{align}
&\bigl(
L_2(E,\alpha,(1+T)^{-1}-1),
L_2(E,\beta,(1+T)^{-1}-1)
\bigr)\\
&\qquad=
\frac{\Omega_{f_E}^+}{\Omega_E^+}
\bigl(
L_2^\sharp(E,(1+T)^{-1}-1),
L_2^\flat(E,(1+T)^{-1}-1)
\bigr)
\mathcal L_{\alpha,\beta}((1+T)^{-1}-1).
\end{align}
Since
$\det\mathcal L_{\alpha,\beta}\neq0$
by \cite[Lemma 6.9]{Spr12},
we obtain
\begin{align}\label{eq:Col_analytic_Kato}
g_\star(T)
=\frac{\Omega_{f_E}^+}{\Omega_E^+}
L_2^\star(E,(1+T)^{-1}-1),
\qquad \star\in\{\sharp,\flat\},
\end{align}
up to a common sign arising from the choice of periods.
Multiplication by a non-zero scalar and the substitution
\[
T\longmapsto (1+T)^{-1}-1
\]
preserve the $\lambda$-invariant.
In particular, if $L_2^\star(E,T)\neq0$, then
$g_\star(T)\neq0$ and \eqref{eq:Col_analytic_Kato}
gives the asserted equality of $\lambda$-invariants.
\end{proof}

\begin{prop}\label{prop:alg_an_lambda}
Let $E/\Q$ be an elliptic curve with good supersingular
reduction at $2$, and let $\star\in\{\sharp,\flat\}$.
Assume that $L_2^\star(E,T)\neq0$.
Then $\cX^\star(E/\Q^{\cyc})$ is $\Lambda$-torsion and
\begin{align}
\lambda_{\alg}^{\star}(E/\Q)
\leq\lambda_{\an}^{\star}(E/\Q).
\end{align}
\end{prop}

\begin{proof}
By Lemma \ref{lem:Col_analytic_comparison},
$g_\star(T)\neq0$.
Hence the argument in the proof of
Theorem \ref{Thm:Lambda_torsion} shows that
$\cX^\star(E/\Q^{\cyc})$ is $\Lambda$-torsion
and that the map $\iota^\star$ in
\eqref{seq_CPT} is injective.

Put
\begin{align}
A:=e_\Delta\mathbb H^1(\cT)[1/2],
\qquad
B:=e_\Delta
\left(
\frac{H^1_{\mathrm{Iw}}(\cT)}{\Ker\col^\star}
\right)[1/2],
\end{align}
\begin{align}
C:=e_\Delta\cX^\star[1/2],
\qquad
D:=e_\Delta\cX^0[1/2].
\end{align}
Then \eqref{seq_CPT} yields an exact sequence
\begin{align}\label{eq:ABCD}
0\longrightarrow B/A
\longrightarrow C
\longrightarrow D
\longrightarrow0.
\end{align}
By \eqref{isom:cyc_Delta_component}, we have
$C\cong \cX^\star(E/\Q^{\cyc})[1/2]$.
Moreover, as in the proof of
Theorem \ref{Thm:Lambda_torsion}, we have
\begin{align}
D\cong e_\Delta\mathbb H^2(\cT)[1/2].
\end{align}
For the element $z$ defined above, put
\begin{align}
Z:=\Lambda_{\Q_2}e_\Delta z\subseteq A.
\end{align}
Under the embedding
$B\subseteq\Lambda_{\Q_2}$
induced by $\col^\star$, the image of $Z$
is $g_\star(T)\Lambda_{\Q_2}$.
Hence we have an exact sequence
\begin{align}\label{eq:AB_L}
0\longrightarrow A/Z
\longrightarrow B/g_\star(T)\Lambda_{\Q_2}
\longrightarrow B/A
\longrightarrow0.
\end{align}
Let $\fp$ be a height-one prime of
$\Lambda$ with $2\notin \fp$.
Since $E$ has good reduction at $2$, the local
correction term in
\cite[Theorem 12.5 (3)]{Kat04} vanishes.
Applying this theorem to the trivial
$\Delta$-component, we obtain
\begin{align}\label{eq:Kato_rational_bound}
\operatorname{length}_{\Lambda_\fp}D_\fp
\leq\operatorname{length}_{\Lambda_\fp}(A/Z)_\fp.
\end{align}
Localizing \eqref{eq:ABCD} and \eqref{eq:AB_L}
at $\fp$ and using \eqref{eq:Kato_rational_bound},
we obtain
\begin{align}
\operatorname{length}_{\Lambda_\fp}C_\fp
&=\operatorname{length}_{\Lambda_\fp}(B/A)_\fp
+\operatorname{length}_{\Lambda_\fp}D_\fp\\
&\leq\operatorname{length}_{\Lambda_\fp}(B/A)_\fp
+\operatorname{length}_{\Lambda_\fp}(A/Z)_\fp\\
&=\operatorname{length}_{\Lambda_\fp}
\parenth{B/g_\star(T)\Lambda_{\Q_2}}_\fp.
\end{align}
Since $B_\fp$ is a torsion-free $\Lambda_\fp$-module
of rank one contained in $\Lambda_\fp$, we have
\begin{align}
\operatorname{length}_{\Lambda_\fp}
\parenth{B/g_\star(T)\Lambda_{\Q_2}}_\fp
\leq\ord_\fp\parenth{g_\star(T)}.
\end{align}
By the isomorphism
$C\cong\cX^\star(E/\Q^{\cyc})[1/2]$
established above, we have
\begin{align}
\operatorname{length}_{\Lambda_\fp}C_\fp
=
\operatorname{length}_{\Lambda_\fp}
\cX^\star(E/\Q^{\cyc})_\fp.
\end{align}
Consequently, by \eqref{eq:Col_analytic_Kato},
for every height-one prime $\fp\neq(2)$ of $\Lambda$,
we have
\begin{align}
\ord_\fp\parenth{
\cha_\Lambda\cX^\star(E/\Q^{\cyc})
}
\leq \ord_\fp(g_\star(T))
=\ord_{\iota(\fp)}\parenth{L_2^\star(E,T)},
\end{align}
where $\iota(\fp)$ is the image of $\fp$
under the involution $\iota$.
Since $\iota$ permutes the height-one
primes of $\Lambda$ other than $(2)$,
Lemma \ref{lem:Col_analytic_comparison} gives
\begin{align}
\lambda_{\alg}^\star(E/\Q)
\leq \lambda(g_\star(T))
=\lambda(L_2^\star(E,T))
=\lambda_{\an}^\star(E/\Q).
\end{align}
\end{proof}

\subsection{Structure of chromatic local conditions}\label{subsec:Herbrand_chromatic}
In this subsection, we determine the structure of the chromatic local conditions at primes lying above $2$ and compute their $\Lambda$-ranks for use in Section \ref{sec:Kida}. Recall the definition of $w$ as
\begin{align}
w:=(a_2-2\varphi)(a_2(a_2-2\varphi)-1)^{-1}\sum_{\sigma\in\gal(k_0/F_{\fp})}\sigma\in R_0.
\end{align}
Although our setting is non-equivariant, the following
two propositions can be regarded as $p=2$ analogues
of the exact sequences in \cite[(4.14) and (4.15)]{Kat21}.
\begin{prop}\label{Prop:flat_condition}
We have $(E^{\flat}_{\infty})^{\vee}\cong R$.
\end{prop}
\begin{proof}
In computing $s_n^*$, we essentially
follow the argument in the proof of
\cite[Theorem 4.16]{Kat21}. For the reader's convenience, we spell out
the relevant maps and identifications. Taking the $\Z_2$-linear dual of \eqref{eq:fundamental_seq}, we obtain
\begin{align}
0\longrightarrow\what{E}(\fm_n)^*
\xrightarrow{t_n^*}\cM_n^*
\xrightarrow{s_n^*}R_0^*\longrightarrow0.
\end{align}
We first describe $s_n^*$ explicitly using the identification of
$\Z_2$-linear duals with group rings in
\cite[(4.1), Lemma 4.2 (2)]{Kat21}.
For $d=1,2$ and $1\leq i\leq d$, let $e_i$ denote the $i$-th standard row vector
of $R_n^{\oplus d}$. For a $\Z_2$-linear map
$f:R_n^{\oplus d}\to\Z_2$, define
\begin{align}
\Theta_{n,d}(f):=
\left(\sum_{g\in G_n}f(ge_i)g\right)_{1\leq i\leq d}
\in R_n^{\oplus d}.
\end{align}
Each $\Theta_{n,d}$ is a $\Z_2$-linear isomorphism. Indeed, for $v=(v_1,\ldots,v_d)\in R_n^{\oplus d}$,
write $v_i=\sum_{g\in G_n}v_{i,g}g$
with $v_{i,g}\in\Z_2$. Then
$\Theta_{n,d}^{-1}(v)$ is the $\Z_2$-linear map
given by
\begin{align}
\sum_{i=1}^d\sum_{g\in G_n}a_{i,g}ge_i
\longmapsto
\sum_{i=1}^d\sum_{g\in G_n}a_{i,g}v_{i,g}.
\end{align}
If $A=(a_{ij})\in M_2(R_n)$ and
$m_A(u):=uA$, then
\begin{align}\label{eq:1}
\Theta_{n,2}(f\circ m_A)
=\Theta_{n,2}(f)(A^{\mathsf T})^\iota.
\end{align}
Indeed, writing $a_{ij}=\sum_{h\in G_n}a_{ij,h}h$,
the $i$-th component of the left-hand side is
\begin{align}
\begin{aligned}
\sum_{g\in G_n}f(ge_iA)g
&=\sum_{g,j,h}a_{ij,h}f(gh e_j)g\\
&=\sum_j\left(\sum_{k\in G_n}f(ke_j)k\right)
       \left(\sum_{h\in G_n}a_{ij,h}h^{-1}\right),
\end{aligned}
\end{align}
which is the $i$-th component of the right-hand side.
Let $\mathrm{pr}_n:R_n^{\oplus2}\twoheadrightarrow\cM_n$ be the
quotient map, and set $D_n:=((\gamma-1)\cB_n)^\iota$.
By the description of $\cM_n^*$ in Subsection
\ref{subsec:Coleman} and
\cite[Lemma 4.9 (2)]{Kat21},
right multiplication by $\widetilde{\cB}_n^\iota$ induces an
isomorphism
\begin{align}
R_n^{\oplus2}/R_n^{\oplus2}D_n
\xrightarrow{\ \sim\ }
\Ker\left(R_n^{\oplus2}\xrightarrow{\times D_n}
R_n^{\oplus2}\right),
\qquad
\overline{(x,y)}\longmapsto
(x,y)\widetilde{\cB}_n^\iota.
\end{align}
Consequently, an element $F\in\cM_n^*$ represented by
$\overline{(x,y)}$ satisfies
\begin{align}\label{eq:2}
\Theta_{n,2}(F\circ\mathrm{pr}_n)
=(x,y)\widetilde{\cB}_n^\iota.
\end{align}
This identity specifies the above identification of
$\cM_n^*$; in particular, $(x,y)$ denotes a representative
of a class in the displayed quotient.

Next, let $P_n:R_n\to R_0$ be the $\Z_2$-linear map
\begin{align}
P_n\left(\sum_{j=0}^{2^n-1}a_j\gamma^j\right):=a_0,
\qquad a_j\in R_0.
\end{align}
Write $\aug_\Gamma:R_n\to R_0$ for the augmentation along
$\Gamma$, so that
\begin{align}
\aug_\Gamma\left(\sum_{j=0}^{2^n-1}
 a_j\gamma^j\right)=\sum_{j=0}^{2^n-1}a_j.
\end{align}
We define the $\Z_2$-linear map $j_n:R_0\longrightarrow R_n^{\oplus2}$ by $j_n(z):=z(-w,1)$.
Write $w=\sum_{\tau\in G_0}w_\tau\tau$ with
$w_\tau\in\Z_2$.
For every $\Z_2$-linear map $f:R_n^{\oplus 2}\to\Z_2$
and $\sigma\in G_0$, we have
\begin{align*}
(f\circ j_n)(\sigma)
=-\sum_{\tau\in G_0}w_\tau f(\sigma\tau e_1)
+f(\sigma e_2).
\end{align*}
Hence
\begin{align}
\Theta_{0,1}(f\circ j_n)
&=-\sum_{\sigma,\tau\in G_0}
w_\tau f(\sigma\tau e_1)\sigma
+\sum_{\sigma\in G_0}f(\sigma e_2)\sigma
\notag\\
&=-\sum_{\rho,\tau\in G_0}
w_\tau f(\rho e_1)\rho\tau^{-1}
+\sum_{\sigma\in G_0}f(\sigma e_2)\sigma
\notag\\
&=-w^\iota P_n\bigl(\Theta_{n,2}(f)_1\bigr)
+P_n\bigl(\Theta_{n,2}(f)_2\bigr).
\label{eq:3}
\end{align}
By the definition of $s_n$, we have
$s_n=\mathrm{pr}_n\circ m_{\cB_n^{\mathsf T}}\circ j_n$.
Put
\begin{align}
\nu_n:=\prod_{j=1}^n N_j
=\sum_{j=0}^{2^n-1}\gamma^j.
\end{align}
Since $\widetilde{\cB}_n\cB_n=\nu_nI_2$ and
$\nu_n^\iota=\nu_n$ in $R_n$, \eqref{eq:1} yields
\begin{align}\label{eq:4}
\Theta_{n,2}(F\circ\mathrm{pr}_n\circ m_{\cB_n^{\mathsf T}})
&=\Theta_{n,2}(F\circ\mathrm{pr}_n)\cB_n^\iota\\
&=(x,y)\widetilde{\cB}_n^\iota\cB_n^\iota\\
&=(x,y)\nu_n.
\end{align}
For $a=\sum_{j=0}^{2^n-1}a_j\gamma^j\in R_n$,
we have $\gamma^j\nu_n=\nu_n$, and hence
\begin{align}
\begin{aligned}
a\nu_n=\left(\sum_{j=0}^{2^n-1}a_j\right)\nu_n
=\aug_\Gamma(a)\nu_n,\quad P_n(a\nu_n)=\aug_\Gamma(a).
\end{aligned}
\end{align}
Since $G_0$ has exponent $2$, we have $w^\iota=w$. By \eqref{eq:3} and \eqref{eq:4}, we have
\begin{align}
\begin{aligned}
\Theta_{0,1}(s_n^*F)
&=\Theta_{0,1}(F\circ\mathrm{pr}_n\circ
 m_{\cB_n^{\mathsf T}}\circ j_n)\\
&=-wP_n(x\nu_n)+P_n(y\nu_n)\\
&=-w\aug_\Gamma(x)+\aug_\Gamma(y).
\end{aligned}
\end{align}
Thus, under $R_0^*\cong R_0$ via $\Theta_{0,1}$, the map
$s_n^*$ is given by
\begin{align}
q'_n(\overline{(x,y)})
=-w\aug_\Gamma(x)+\aug_\Gamma(y).
\end{align}
This formula is independent of the representative since
$\aug_\Gamma(D_n)=0$.
It is compatible with the natural transition maps of
Subsection \ref{subsec:Coleman}, since augmentation commutes
with the projections $R_n\to R_{n-1}$.
Taking the projective limit, as in the proof of
\cite[Theorem 4.16]{Kat21}, gives the exact sequence
\begin{align}\label{eq:exact_Col_prime}
\xymatrix{
0\ar[r]
&\varprojlim_n\what{E}(\fm_n)^*\ar^-{\col'}[r]
&R^{\oplus2}\ar^{q'}[r]
&R_0\ar[r]&0,
}
\end{align}
where
\begin{align}
q'(x,y)=-w\aug_\Gamma(x)+\aug_\Gamma(y).
\end{align}
Here the right-hand term $R_0$ is identified with $R_0^*$
by $\Theta_{0,1}$.
We put $q:=q'\circ m_{U^{-1}}$. Since $U^{-1}
=\begin{pmatrix}0&1\\-1&-a_2\end{pmatrix}$, we see that
\begin{align}
\begin{aligned}
q(x,y)
&=q'(-y,x-a_2y)\\
&=\aug_\Gamma(x)+(w-a_2)\aug_\Gamma(y)
\end{aligned}
\end{align}
and we obtain the exact sequence
\begin{align}\label{eq:exact_Col}
\xymatrix{
0\ar[r]
&\varprojlim_n\what{E}(\fm_n)^*\ar^-{\col}[r]
&R^{\oplus2}\ar^{q}[r]
&R_0\ar[r]&0.
}
\end{align}
For any $y\in R$, since
$\aug_\Gamma:R\to R_0$ is surjective, we can choose
$x\in R$ such that
\begin{align}
\aug_\Gamma(x)=-(w-a_2)\aug_\Gamma(y).
\end{align}
Then $q(x,y)=0$, hence $(x,y)\in\Ker(q)=\Imag(\col)$.
Hence the second projection of $\Imag(\col)$ is all of $R$.
By Definition \ref{Definition: Local_condition},
\begin{align}
(E^{\flat}_{\infty})^{\vee}
\cong\Imag(\col^{\flat})=R.
\end{align}
\end{proof}

\begin{prop}\label{Prop:sharp_condition}
Put $N_{\Delta}:=1+\delta$. There is an exact sequence
\begin{align}\label{seq:sharp_condition}
\xymatrix{
0\ar[r] & (E^{\sharp}_{\infty})^{\vee}\ar[r] & R \ar[r] & R_0/(w-a_2)\ar[r] & 0,
}
\end{align}
and
\begin{align}
R_0/(w-a_2)\cong \begin{cases}
R_0/(2N_{\Delta})&(a_2=0),\\
R_0/2R_0&(a_2=\pm 2).
\end{cases}
\end{align}
\end{prop}

\begin{proof}
By \eqref{eq:exact_Col}, an element $x\in R$ belongs to
$\Imag(\col^{\sharp})$ if and only if there exists $y\in R$ such that
\begin{align}
\aug_{\Gamma}(x)+(w-a_2)\aug_{\Gamma}(y)=0.
\end{align}
Since $\aug_{\Gamma}:R\to R_0$ is surjective, this is equivalent to
\begin{align}
\aug_{\Gamma}(x)\in (w-a_2) R_0.
\end{align}
Therefore
\begin{align}
\Imag(\col^{\sharp})
=
\Ker\left(
R\overset{\aug_{\Gamma}}{\longrightarrow}R_0
\longrightarrow R_0/(w-a_2)
\right).
\end{align}
The composite map is surjective, and Definition \ref{Definition: Local_condition} gives
$(E^{\sharp}_{\infty})^{\vee}\cong\Imag(\col^{\sharp})$. This proves
the exactness of \eqref{seq:sharp_condition}. We now determine $R_0/(w-a_2)$. Put
\begin{align}
u:=a_2-2\varphi,
\qquad
v:=a_2u-1.
\end{align}
Since $a_2\in2\Z_2$, we have $v\equiv-1\pmod{2R_0}$, and hence
$v\in R_0^{\times}$. By the definition of $w$, $w=uv^{-1}N_{\Delta}$. Thus
\begin{align}
(w-a_2)=(v^{-1}(uN_{\Delta}-a_2v))=(uN_{\Delta}-a_2v).
\end{align}
We put $r:=uN_{\Delta}-a_2v$. Suppose first that $a_2=0$. Then $u=-2\varphi$ and $v=-1$, and hence $r=-2\varphi N_{\Delta}$.
Since $\varphi\in R_0^{\times}$, we obtain
\begin{align}
R_0/(w-a_2)\cong R_0/(2N_{\Delta}).
\end{align}
Suppose next that $a_2=\pm2$, and put $\varepsilon:=a_2/2\in\{\pm1\}$.
Then
\begin{align}
\frac{r}{2}
=(\varepsilon-\varphi)N_{\Delta}-\varepsilon v.
\end{align}
Since $G_0$ is a finite $2$-group, $R_0=\Z_2[G_0]$ is local with
maximal ideal
\begin{align}
\fM_0:=(2,\sigma-1\mid \sigma\in G_0).
\end{align}
The congruences
\begin{align}
\varepsilon\equiv1,\quad \varphi\equiv1,\quad v\equiv1\pmod{\fM_0}
\end{align}
imply that
\begin{align}
\frac{r}{2}\equiv1\pmod{\fM_0}.
\end{align}
It follows that $r/2\in R_0^{\times}$, hence $(r)=(2)$. Therefore
\begin{align}
R_0/(w-a_2)\cong R_0/2R_0.
\end{align}
\end{proof}

\begin{cor}\label{cor:rank_local_condition}
Let $\star\in\{\sharp,\flat\}$. For every prime $\fp\mid2$ of $F$, we have
\begin{align}
\rank_{\Lambda}\parenth{E^\star_{\infty,\fp}}^\vee
=
[F_\fp(\zeta_4):\Q_2]
=
2[F_\fp:\Q_2].
\end{align}
In particular,
\begin{align}
\rank_{\Lambda}
\parenth{
\prod_{\fp\mid2}E^\star_{\infty,\fp}
}^\vee
=4.
\end{align}
\end{cor}

\begin{proof}
This follows immediately from Propositions \ref{Prop:flat_condition} and \ref{Prop:sharp_condition}.
\end{proof}

\section{A $2$-adic Kida-type formula for chromatic Selmer groups}\label{sec:Kida}

In this section, we prove a $2$-adic Kida-type formula for doubly-signed Selmer groups. There are several approaches for proving Kida-type formulae for elliptic curves. For proofs using congruences, see \cite{PW06} and \cite{HL19}, while for proofs using quasi-projectivity, see \cite{Gre11}, \cite{Lim22}, and \cite{FM25}, for example.  In this paper, we adopt the method of calculating the Herbrand quotient as in \cite{Iwa81}, \cite{HM99}, \cite{Mat08}, and \cite{RS24}.
We retain the notation $M_n$, $M_n^{\cyc}$, and
$M_\infty$ introduced in Subsection \ref{subsec:doubly_Selmer} for $M=\Q$ or $F$.

\subsection{Kida-type formulae and Herbrand quotients}\label{subsec:HM_method}
In this subsection, we extend the method of Hachimori--Matsuno \cite{HM99}
and Matsuno \cite{Mat08} to the supersingular setting. The main additional difficulty in the supersingular
case arises from the local conditions at $2$.
In the ordinary case, the local conditions are
given by the images of the local Kummer maps,
whose compatibility with restriction follows
from their functoriality; see
\cite[Lemma 3.2]{Mat08} for the corresponding result.
In contrast, our chromatic local conditions
are defined through the kernels of the Coleman
maps, and hence require a direct comparison
of these maps under corestriction.

For simplicity, we put
$G:=\gal(F^{\cyc}/\Q^{\cyc})$.
For a $G$-module $\cA$ such that
$H^1(G,\cA)$ and $H^2(G,\cA)$ are finite,
we define its Herbrand quotient by
\begin{align}
h_G(\cA):=\frac{\sharp H^2(G,\cA)}{\sharp H^1(G,\cA)}.
\end{align}
In particular, $h_G(\cA)=1$ if $\cA$ is finite. Let $K/\Q_2$ be the unramified quadratic
extension. For $M=\Q_2$ or $K$, we put
\begin{align}
R_{M,n}:=\Z_2[\gal(M_n/\Q_2)]
\quad(n\geq-2),
\qquad
R_M:=\Z_2\llbracket\gal(M_\infty/\Q_2)\rrbracket.
\end{align}
We write $d_{M,n}$, $\eta_{M,n}$, $\pi_{M,n}$,
and $f_M$ for the corresponding objects constructed
in Subsection \ref{subsec:local_points} over $M$. The following lemma compares the systems of local points over $K$ and $\Q_2$. It will be used to
establish the compatibility of the chromatic Coleman maps with corestriction when $2$ is inert in $F/\Q$.
\begin{lem}\label{trace_local_points}
For $n\geq-2$,
\begin{align}
\Tr_{K_n/\Q_{2,n}}(d_{K,n})=-d_{\Q_2,n}.
\end{align}
\end{lem}
\begin{proof}
By Remark \ref{Rem:phi_relation_eta_n}, we have
\begin{align}\label{eq:eta_trace_zeta}
 (3-a_2\varphi)\eta_{K,n}=2\zeta^{\varphi^{-n}}.
\end{align}
Applying $\Tr_{K/\Q_2}=1+\varphi$ to \eqref{eq:eta_trace_zeta}, we obtain
\begin{align}
 (3-a_2)\Tr_{K/\Q_2}(\eta_{K,n})=2\Tr_{K/\Q_2}\left(\zeta^{\varphi^{-n}}\right)=-2.
\end{align}
On the other hand, the corresponding relation
over $\Q_2$ is
\begin{align}
(3-a_2)\eta_{\Q_2,n}=2.
\end{align}
Since $3-a_2\in\Z_2^\times$, it follows that
\begin{align}\label{eq:trace_eta_K}
\Tr_{K/\Q_2}(\eta_{K,n})
=-\eta_{\Q_2,n}.
\end{align}
By Definition \ref{Def_power}, we have
\begin{align}
f_K(T)=\sum_{j\geq0}x_j
\bigl((T+\zeta)^{2^j}-\zeta^{2^j}\bigr),
\end{align}
whereas
\begin{align}
f_{\Q_2}(T)=\sum_{j\geq0}x_j
\bigl((T+1)^{2^j}-1\bigr).
\end{align}
Since $K/\Q_2$ is unramified and
$\Q_{2,n}/\Q_2$ is totally ramified, we have
\begin{align}
\gal(K_n/\Q_{2,n})
\cong \gal(K/\Q_2)=\{1,\varphi\}.
\end{align}
In particular, for $a\in K$ and $b\in\Q_{2,n}$,
\begin{align}
\Tr_{K_n/\Q_{2,n}}(ab)
=\Tr_{K/\Q_2}(a)b.
\end{align}
Therefore
\begin{align}
&\Tr_{K_n/\Q_{2,n}}
\left(
f_K^{\varphi^{-n}}(\pi_{K,n})
\right)\notag\\
&\quad=
\sum_{j\geq0}x_j
\Tr_{K/\Q_2}
\left(
\zeta^{\varphi^{-n}2^j}
\right)(\zeta_{2^{n+2}}^{2^j}-1)\notag\\
&\quad=
-\sum_{j\geq0}x_j(\zeta_{2^{n+2}}^{2^j}-1)\notag\\
&\quad=-f_{\Q_2}(\pi_{\Q_2,n}),
\label{eq:trace_f_K}
\end{align}
where we used the fact that
$\zeta^{2^j}$ is a primitive third root of unity
for every $j\geq0$. By the definition of $d_{M,n}$, we have
\begin{align}
\log_{\what E}(d_{K,n})
=\eta_{K,n}
+f_K^{\varphi^{-n}}(\pi_{K,n})
\end{align}
and
\begin{align}
\log_{\what E}(d_{\Q_2,n})
=\eta_{\Q_2,n}
+f_{\Q_2}(\pi_{\Q_2,n}).
\end{align}
Hence \eqref{eq:trace_eta_K} and
\eqref{eq:trace_f_K} imply
\begin{align}
\Tr_{K_n/\Q_{2,n}}
\left(\log_{\what E}(d_{K,n})\right)
=-\log_{\what E}(d_{\Q_2,n}).
\end{align}
Since $\log_{\what E}$ is injective
by Lemma \ref{injective_log}, we conclude that
\begin{align}
\Tr_{K_n/\Q_{2,n}}(d_{K,n})
=-d_{\Q_2,n}.
\end{align}
\end{proof}
We put $G_{\mathrm{ur}}:=\gal(K_{\infty}/\Q_{2,\infty})$. For $\star\in\{\sharp,\flat\}$, taking $G_{\mathrm{ur}}$-coinvariants of
\begin{align}
\xymatrix{
0\ar[r] & \Ker(\col^{\star}_K) \ar[r] & H^1_{\iw,K}(\cT)\ar[r] & \Imag(\col^{\star}_K)\ar[r] &0
},
\end{align}
we obtain
\begin{align}
\xymatrix{
0\ar[r] & \Ker(\col^{\star}_K)_{G_{\mathrm{ur}}}/A \ar[r]  & H^1_{\iw,K}(\cT)_{G_{\mathrm{ur}}}\ar[r] & \Imag(\col^{\star}_K)_{G_{\mathrm{ur}}}\ar[r]&0,
}
\end{align}
where $A:=\Imag(H_1(G_{\mathrm{ur}},\Imag(\col^{\star}_K))\to \Ker(\col^{\star}_K)_{G_{\mathrm{ur}}})$.

\begin{lem}\label{Lem:cor_aug}
For $\star\in\{\sharp,\flat\}$, the following diagram is commutative.
\begin{align}\label{diagram:cor_aug}
\xymatrix{
0\ar[r] & \Ker(\col^{\star}_K)_{G_{\mathrm{ur}}}/A \ar[r] \ar^{\core}[d] & H^1_{\iw,K}(\cT)_{G_{\mathrm{ur}}}\ar^{\core}[d]\ar^{\col^{\star}_K}[r] & \Imag(\col^{\star}_K)_{G_{\mathrm{ur}}}\ar[r]\ar^{-\augm}[d] &0\\
0\ar[r] & \Ker(\col^{\star}_{\Q_2}) \ar[r] & H^1_{\iw,\Q_2}(\cT)\ar^{\col^{\star}_{\Q_2}}[r] & \Imag(\col^{\star}_{\Q_2})\ar[r] &0
}
\end{align}
Here, two vertical maps labeled $\core$ are induced by the corestriction map
\begin{align}
H^1_{\iw,K}(\cT)\to H^1_{\iw,\Q_2}(\cT),
\end{align}
while $\augm$ denotes the augmentation map
\begin{align}
\Z_2\llbracket \gal(K_{\infty}/\Q_2)\rrbracket\to \Z_2\llbracket \gal(\Q_{2,\infty}/\Q_2)\rrbracket.
\end{align}
\end{lem}
\begin{proof}
Since the commutativity of the left square is clear,
we show that the right square is commutative.
Fix $n\geq0$ and $z_n\in H^1(K_n,\cT)$. We identify $G_{\mathrm{ur}}$ with $\gal(K_n/\Q_{2,n})$.
For each $\sigma\in\gal(\Q_{2,n}/\Q_2)$,
let $\widetilde{\sigma}\in\gal(K_n/\Q_2)$
be its unique lift acting trivially on $K$.
Let $\rest_{K_n/\Q_{2,n}}$ denote the natural
inclusion of formal group points induced by
$\Q_{2,n}\subset K_n$.

By Lemma \ref{trace_local_points}, for $j=n,n-1$, we have
\begin{align}
&\sum_{\tau\in G_{\mathrm{ur}}}
(\tau\widetilde{\sigma}(d_{K,j}),z_n)_{K_n}\notag\\
&=(\rest_{K_n/\Q_{2,n}}
(\sigma\Tr_{K_n/\Q_{2,n}}(d_{K,j})),z_n)_{K_n}
\notag\\
&=-(\rest_{K_n/\Q_{2,n}}
(\sigma(d_{\Q_2,j})),z_n)_{K_n}\notag\\
&=-(\sigma(d_{\Q_2,j}),
\core_{K_n/\Q_{2,n}}(z_n))_{\Q_{2,n}}.
\end{align}
Hence, summing over
$\sigma\in\gal(\Q_{2,n}/\Q_2)$, we obtain
\begin{align}\label{eq:aug_P}
&\augm_n\left(
\sum_{g\in\gal(K_n/\Q_2)}
(g(d_{K,j}),z_n)_{K_n}g
\right)\notag\\
&=-\sum_{\sigma\in\gal(\Q_{2,n}/\Q_2)}
(\sigma(d_{\Q_2,j}),
\core_{K_n/\Q_{2,n}}(z_n))_{\Q_{2,n}}\sigma
\qquad (j=n,n-1).
\end{align}
Let $\widetilde{\cB}_{M,n}$ denote the adjugate
matrix of $\cB_{M,n}$ for $M=K$ or $\Q_2$.
By \eqref{eq:aug_P} and
\cite[(4.7)]{Kat21}, we obtain
\begin{align}
\augm_n(\col'_{K,n}(z_n))
\widetilde{\cB}_{\Q_2,n}^{\iota}=\augm_n\left(\col'_{K,n}(z_n)\widetilde{\cB}_{K,n}^{\iota}\right)=-\col'_{\Q_2,n}(\core_{K_n/\Q_{2,n}}(z_n))\widetilde{\cB}_{\Q_2,n}^{\iota}.
\end{align}
Since right multiplication by
$\widetilde{\cB}_{\Q_2,n}^{\iota}$
induces an injective map 
\begin{align}
\frac{R_{\Q_2,n}^{\oplus2}}
{R_{\Q_2,n}^{\oplus2}
((\gamma-1)\cB_{\Q_2,n})^\iota}
\hookrightarrow R_{\Q_2,n}^{\oplus2}
\end{align}
by \cite[Lemma 4.9 (2)]{Kat21}, we have
\begin{align}
\augm_n(\col'_{K,n}(z_n))
=-\col'_{\Q_2,n}(\core_{K_n/\Q_{2,n}}(z_n)).
\end{align}
Taking the projective limit, multiplying by
\begin{align}
\begin{pmatrix}
-a_2&-1\\
1&0
\end{pmatrix},
\end{align}
and taking the $\star$-component, we have
\begin{align}\label{eq:Col_aug_compatible}
\col^\star_{\Q_2}\circ\core
=-\augm\circ\col^\star_K.
\end{align}
By \eqref{eq:Col_aug_compatible}, the augmentation map sends
$\Imag(\col_K^\star)$ into $\Imag(\col_{\Q_2}^\star)$.
Since $\augm(\tau)=1$ for every $\tau\in G_{\mathrm{ur}}$,
it induces a homomorphism
\begin{align}
-\augm:
\Imag(\col_K^\star)_{G_{\mathrm{ur}}}
\longrightarrow \Imag(\col_{\Q_2}^\star).
\end{align}
This completes the proof.
\end{proof}

\begin{lem}\label{Lem:ker_finite}
For $\star\in\{\sharp,\flat\}$ and $i\geq1$, we have
\begin{align}\label{eq:col_image_G2_vanishing}
H_i(G_{\mathrm{ur}},\Imag(\col_K^{\star}))=0
\end{align}
and the augmentation map induces an isomorphism
\begin{align}
-\augm:\Imag(\col_K^{\star})_{G_{\mathrm{ur}}}\overset{\sim}{\longrightarrow}\Imag(\col_{\Q_2}^{\star}).
\end{align}
\end{lem}

\begin{proof}
For $M=\Q_2$ or $K$, Propositions
\ref{Prop:flat_condition} and
\ref{Prop:sharp_condition} give
\begin{align}
\Imag(\col_M^\flat)&=R_M,\\
\Imag(\col_M^\sharp)
&=\left\{
x\in R_M\;\middle|\;
\aug_\Gamma(x)\in(w-a_2)R_{M,0}
\right\}.
\end{align}
By the computation in the proof of
Proposition \ref{Prop:sharp_condition}, we have
\begin{align}
(w-a_2)R_{M,0}
=
\begin{cases}
2N_\Delta R_{M,0} &(a_2=0),\\
2R_{M,0} &(a_2=\pm2).
\end{cases}
\end{align}
Since $\aug_\Gamma:R_M\to R_{M,0}$ is
surjective with kernel $(\gamma-1)R_M$,
we obtain
\begin{align}\label{eq:col_images_for_descent}
\Imag(\col_M^\flat)=R_M,\quad
\Imag(\col_M^\sharp)
=
\begin{cases}
(\gamma-1,2N_\Delta)R_M &(a_2=0),\\
(\gamma-1,2)R_M &(a_2=\pm2).
\end{cases}
\end{align}
Since $K/\Q_2$ is unramified and
$\Q_{2,\infty}/\Q_2$ is totally ramified,
we have
\begin{align}
R_K\cong
\Z_2[G_{\mathrm{ur}}]\otimes_{\Z_2}R_{\Q_2}.
\end{align}
As the generators in \eqref{eq:col_images_for_descent} belong to $R_{\Q_2}$, we have
\begin{align}
\Imag(\col_K^\star)
&=\Imag(\col_{\Q_2}^\star)R_K\\
&=\Imag(\col_{\Q_2}^\star)
\oplus\varphi\Imag(\col_{\Q_2}^\star).
\end{align}
Here, $\varphi$ interchanges the two summands.
Thus multiplication induces an isomorphism
\begin{align}\label{eq:col_image_induced}
\Z_2[G_{\mathrm{ur}}]\otimes_{\Z_2}
\Imag(\col_{\Q_2}^\star)
&\xrightarrow{\sim}\Imag(\col_K^\star)
\end{align}
of $\Z_2[G_{\mathrm{ur}}]$-modules,
compatibly with the $\Delta$-action.
By Shapiro's lemma,
\eqref{eq:col_image_induced} gives
\begin{align}
H_i(G_{\mathrm{ur}},\Imag(\col_K^\star))=0
\qquad(i\geq1).
\end{align}
Furthermore, taking $G_{\mathrm{ur}}$-coinvariants
in \eqref{eq:col_image_induced},
the augmentation map is identified with
the canonical isomorphism
\begin{align}
\left(
\Z_2[G_{\mathrm{ur}}]\otimes_{\Z_2}
\Imag(\col_{\Q_2}^\star)
\right)_{G_{\mathrm{ur}}}
\xrightarrow{\ \sim\ }
\Imag(\col_{\Q_2}^\star).
\end{align}
Since multiplication by $-1$ is an
automorphism, the map
\begin{align}
-\augm:
\Imag(\col_K^\star)_{G_{\mathrm{ur}}}
\xrightarrow{\ \sim\ }
\Imag(\col_{\Q_2}^\star)
\end{align}
is also an isomorphism.
\end{proof}

We first note that the natural restriction maps
on local Galois cohomology induce maps on the
quotients by the chromatic local conditions.
Indeed, when $2$ is inert in $F/\Q$,
\eqref{eq:Col_aug_compatible} implies that
\begin{align}
\core_{K_\infty/\Q_{2,\infty}}
(\Ker(\col_K^\star))
\subseteq \Ker(\col_{\Q_2}^\star).
\end{align}
By the adjointness of restriction and corestriction
under local Tate duality, we obtain
\begin{align}
\rest_{K_\infty/\Q_{2,\infty}}
(E_{\infty,2}^\star)
\subseteq E_{\infty,\fp}^\star.
\end{align}
When $2$ splits in $F/\Q$, the compatibility
is immediate. Hence the restriction maps induce
the map
\begin{align}\label{restriction_map}
 \frac{H^1(\Q_{2,\infty},E[2^\infty])}
      {E^\star_{\infty,2}}
 \longrightarrow
 \left(
 \prod_{\fp\mid 2}
 \frac{H^1(F_{\fp,\infty},E[2^\infty])}
      {E^\star_{\infty,\fp}}
 \right)^G.
\end{align}
\begin{prop}\label{lem:local-control}
Let $\star\in\{\sharp,\flat\}$. Then the natural restriction map
\eqref{restriction_map} is an isomorphism.
\end{prop}

\begin{proof}
If $2$ splits in $F/\Q$, the assertion is immediate.
Suppose that $2$ is inert in $F$, namely
$F_{\fp}\cong K$. Then we identify $G$ with $G_{\mathrm{ur}}$.
By Lemma \ref{lem:2-primary_torsion},
the inflation--restriction sequence and local Tate
duality show that corestriction induces an isomorphism
\[
H^1_{\iw,K}(\cT)_{G_{\mathrm{ur}}}
\xrightarrow{\sim}
H^1_{\iw,\Q_2}(\cT).
\]
Moreover, Lemma \ref{Lem:ker_finite}
gives $A=0$ and shows that the right vertical map
in \eqref{diagram:cor_aug} is an isomorphism.
Hence, applying the snake lemma to
\eqref{diagram:cor_aug}, we obtain an isomorphism
\[
\core:\Ker(\col_K^\star)_{G_{\mathrm{ur}}}
\xrightarrow{\sim}
\Ker(\col_{\Q_2}^\star).
\]
By local Tate duality and the adjointness of restriction
and corestriction, this isomorphism is dual to
\eqref{restriction_map}.
\end{proof}

For $M=\Q$ or $F$, fix $\vec{s}=(s_\fp)_{\fp\in\Sigma_2(M)}\in\prod_{\fp\in\Sigma_2(M)}\{\sharp,\flat\}$. If $w$ is a prime of $M_{\infty}$ above $\fp\mid 2$,  write $s_w=s_{\fp}$. For $\fw\in\Sigma(M_{\infty})$, define
\begin{align}
\cH_{\fw}(E/M_{\infty}):=\begin{cases}
\frac{H^1(M_{\infty,\fw},E[2^\infty])}
     {E^{s_\fw}_{\infty,\fw}}&(\fw\in\Sigma_2(M_\infty))\\
H^1(M_{\infty,\fw},E[2^\infty])
&(\fw\in\Sigma(M_\infty)\setminus\Sigma_2(M_\infty)).
\end{cases}
\end{align}
Here, the dependence of $\cH_{\fw}(E/M_\infty)$
on the choice of $\vec{s}$ is suppressed from
the notation. We put
\begin{align}
\cH^{\vec{s}}(M_\infty)
:=
\prod_{\fw\in\Sigma(M_\infty)}
\cH_{\fw}(E/M_\infty).
\end{align}
For $\star\in\{\sharp,\flat\}$, we also write
$\cH^\star(M_\infty):=\cH^{\vec{s}_\star}(M_\infty)$,
where $\vec{s}_\star:=(\star)_{\fp\in\Sigma_2(M)}$.

\begin{cor}\label{Cor:Selmer_control}
The kernel and cokernel of the restriction map
\begin{align}
\sel^{\star}(E/\Q^{\cyc})\to \sel^{\vec{s}_{\star}}(E/F^{\cyc})^G
\end{align}
are finite. 
\end{cor}

\begin{proof}
We adapt the arguments of Hachimori--Matsuno
\cite[Lemma 3.3]{HM99} and Matsuno
\cite[Lemma 3.2]{Mat08}, first comparing the
relaxed Selmer groups and then imposing the
archimedean local conditions of Definition
\ref{Definition_Selmer_chromatic_cyclotomic}.

Since $E(F^{\cyc})[2^\infty]=0$ by Lemma
\ref{lem:2-primary_torsion}, the inflation--restriction
sequence gives an isomorphism
\begin{align}
\alpha:
H^1(\Q_\Sigma/\Q^{\cyc},E[2^\infty])
\xrightarrow{\;\sim\;}
H^1(F_\Sigma/F^{\cyc},E[2^\infty])^G.
\end{align}
Let
\begin{align}
\beta:
\cH^{\star}(\Q_\infty)^\Delta
\longrightarrow
\left(
\cH^{\vec{s}_\star}(F_\infty)^\Delta
\right)^G
\end{align}
be the map induced by restriction.
We claim that $\Ker(\beta)$ is finite.
For this, we consider the local restriction maps at each prime.
At the primes above $2$, the corresponding component
of $\beta$ is an isomorphism by Proposition
\ref{lem:local-control}, even after taking
$\Delta$-invariants.

For an odd prime $\ell\in\Sigma_\Q$,
let $v$ be a prime of $\Q^{\cyc}$ above
$\ell$ and $w$ a prime of $F^{\cyc}$
above $v$.
Since $(\Q^{\cyc})_v$ contains $\zeta_4$,
the extension $M_\infty/M^{\cyc}$ is split
at primes above $\ell$ for
$M\in\{\Q,F\}$.
Hence, for each prime $u$ of $M^{\cyc}$
above $\ell$, we have
\begin{align}
\left(
\prod_{\fw\mid u}\cH_{\fw}(E/M_\infty)
\right)^\Delta
\cong
H^1((M^{\cyc})_u,E)[2^\infty],
\end{align}
where $\fw$ runs over the primes of $M_\infty$
above $u$. In particular, the component of
$\beta$ at $v$ is identified with
\begin{align}
H^1((\Q^{\cyc})_v,E)[2^\infty]
\longrightarrow
\left(
\prod_{w\mid v}
H^1((F^{\cyc})_w,E)[2^\infty]
\right)^G.
\end{align}
Put
$G_w:=\gal((F^{\cyc})_w/(\Q^{\cyc})_v)$.
The inflation--restriction sequence gives an isomorphism
\begin{align}
\Ker\left(
H^1((\Q^{\cyc})_v,E)[2^\infty]
\longrightarrow
\left(
\prod_{w\mid v}
H^1((F^{\cyc})_w,E)[2^\infty]
\right)^G
\right)
\overset{\sim}{\longrightarrow}
H^1(G_w,E((F^{\cyc})_w)).
\end{align}
Since $\ell\neq 2$, we have
\begin{align}
E((F^{\cyc})_w)\otimes_{\Z}\Q_2/\Z_2=0,\quad H^1(G_w,E((F^{\cyc})_w))
&\cong
H^1(G_w,E((F^{\cyc})_w)[2^\infty]).
\end{align}
As $E((F^{\cyc})_w)[2^\infty]$ is
cofinitely generated over $\Z_2$, the group
$H^1(G_w,E((F^{\cyc})_w)[2^\infty])$ is finite (cf. \cite[Proof of Lemma 3.3]{HM99}).
Thus the local restriction map at $v$ has
finite kernel.
Since $\Sigma_\Q$ is finite and there are
only finitely many primes above each odd
$\ell$ in $\Q^{\cyc}$, we conclude that
$\Ker(\beta)$ is finite.

Let $\rho_{\mathrm{rel}}$ denote the natural
restriction map
\begin{align}
\sel_{\mathrm{rel}}^{\star}
(E/\Q^{\cyc})
\longrightarrow
\sel_{\mathrm{rel}}^{\vec{s}_\star}
(E/F^{\cyc})^G.
\end{align}
The restriction maps induce the following
commutative diagram:
\begin{align}
\xymatrix@C=1.5em@R=2.5em{
0\ar[r]
&\sel_{\mathrm{rel}}^{\star}(E/\Q^{\cyc})
 \ar[r]\ar[d]^{\rho_{\mathrm{rel}}}
&H^1(\Q_\Sigma/\Q^{\cyc},E[2^\infty])
 \ar[r]\ar[d]^{\alpha}_{\simeq}
&\cH^{\star}(\Q_\infty)^\Delta
 \ar[d]^{\beta}\\
0\ar[r]
&\sel_{\mathrm{rel}}^{\vec{s}_\star}(E/F^{\cyc})^G
 \ar[r]
&H^1(F_\Sigma/F^{\cyc},E[2^\infty])^G
 \ar[r]
&\left(\cH^{\vec{s}_\star}(F_\infty)^\Delta\right)^G.
}
\end{align}
Replacing the rightmost terms by the
images of the horizontal maps and
applying the snake lemma, we obtain
\begin{align}
\Ker(\rho_{\mathrm{rel}})=0,\qquad
\Coker(\rho_{\mathrm{rel}})
\hookrightarrow\Ker(\beta).
\end{align}
Since $\Ker(\beta)$ is finite,
$\rho_{\mathrm{rel}}$ is injective
with finite cokernel.

We next consider the local conditions
at the archimedean places.
Since $F/\Q$ is real quadratic, every
real place of $\Q^{\cyc}$ splits
completely in $F^{\cyc}$.
Consequently, restriction induces a
natural isomorphism
\begin{align}\label{eq:archimedean_control}
\cH_\infty(\Q,E)
\xrightarrow{\;\sim\;}
\cH_\infty(F,E)^G.
\end{align}
Let $\rho$ denote the natural restriction map
\begin{align}
\sel^{\star}(E/\Q^{\cyc})
\longrightarrow
\sel^{\vec{s}_\star}(E/F^{\cyc})^G.
\end{align}
By \eqref{eq:archimedean_control},
we have the following commutative diagram:
\begin{align}
\xymatrix@C=2.5em@R=2.5em{
0\ar[r]
&\sel^{\star}(E/\Q^{\cyc})
 \ar[r]\ar[d]^{\rho}
&\sel_{\mathrm{rel}}^{\star}(E/\Q^{\cyc})
 \ar[r]^-{\loc_\infty}
 \ar[d]^{\rho_{\mathrm{rel}}}
&\cH_\infty(\Q,E)\ar[d]^{\simeq}\\
0\ar[r]
&\sel^{\vec{s}_\star}(E/F^{\cyc})^G
 \ar[r]
&\sel_{\mathrm{rel}}^{\vec{s}_\star}(E/F^{\cyc})^G
 \ar[r]_-{(\loc_\infty)^G}
&\cH_\infty(F,E)^G.
}
\end{align}
Replacing the rightmost terms by the
images of the horizontal maps and
applying the snake lemma, we obtain
\begin{align}
\Ker(\rho)=0,\qquad
\Coker(\rho)
\hookrightarrow\Coker(\rho_{\mathrm{rel}}).
\end{align}
Since $\Coker(\rho_{\mathrm{rel}})$ is finite,
the assertion follows.
\end{proof}

\begin{cor}\label{Cor:Kida_Herbrand}
Let $E$ be an elliptic curve defined over $\Q$ with good supersingular reduction at $2$. Let $\star\in\{\sharp,\flat\}$. Assume that $\sel^{\star}(E/\Q^{\cyc})$ is $\Lambda$-cotorsion. If $\mu^{\star}(E/\Q)=0$, then $\sel^{\vec{s}_{\star}}(E/F^{\cyc})$ is also $\Lambda$-cotorsion, $\mu^{\vec{s}_{\star}}(E/F)=0$, and we have
\begin{align}
\lambda^{\vec{s}_{\star}}(E/F)=2\lambda^{\star}(E/\Q)+v_2(h_G(\sel^{\vec{s}_{\star}}(E/F^{\cyc}))).
\end{align}
\end{cor}
\begin{proof}
By Corollary \ref{Cor:Selmer_control}, the map
\begin{align}
\sel^\star(E/\Q^{\cyc})\longrightarrow
\sel^{\vec{s}_\star}(E/F^{\cyc})^G
\end{align}
has finite kernel and cokernel. The hypothesis
$\mu^\star(E/\Q)=0$ implies that
$\sel^\star(E/\Q^{\cyc})$ is cofinitely generated over $\Z_2$.
It follows that $\cX^{\vec{s}_\star}(E/F^{\cyc})_G$ is finitely generated over $\Z_2$.
As $G$ is a group of order $2$, the topological Nakayama lemma
shows that $\cX^{\vec{s}_\star}(E/F^{\cyc})$ is finitely generated over $\Z_2[G]$, and hence
over $\Z_2$. Thus
\begin{align}
\mu^{\vec{s}_\star}(E/F)=0,
\qquad
\corank_{\Z_2}\sel^{\vec{s}_\star}(E/F^{\cyc})
=\lambda^{\vec{s}_\star}(E/F).
\end{align}
By the formula for cofinitely
generated $\Z_2$-modules (cf. \cite[\S 9, pp.~285--287]{Iwa81}
and \cite[p.~589]{HM99}),
we have
\begin{align}
v_2\parenth{
h_G\parenth{
\sel^{\vec{s}_\star}(E/F^{\cyc})
}}
&=
\corank_{\Z_2}\sel^{\vec{s}_\star}(E/F^{\cyc})
-2\corank_{\Z_2}\sel^{\vec{s}_\star}(E/F^{\cyc})^G\\
&=
\lambda^{\vec{s}_\star}(E/F)
-2\lambda^\star(E/\Q),
\end{align}
where the last equality follows from Corollary \ref{Cor:Selmer_control}.
\end{proof}

\subsection{A proof of the Kida-type formula}\label{subsec:proof_Kida}

We use the structure of the chromatic local conditions
established in \S \ref{subsec:Herbrand_chromatic}
to prove the Kida-type formula.
We first establish the necessary localization results,
taking account of the archimedean local conditions.

\begin{lem}\label{Lem:chromatic_descent_error}
Let $M=\Q$ or $F$, and let $v\in\Sigma_2(M^{\cyc})$.
For $\star\in\{\sharp,\flat\}$, the cokernel of the natural map
\begin{align}\label{eq:local_chromatic_descent}
H^1((M^{\cyc})_v,E[2^\infty])
\longrightarrow
\parenth{
\frac{H^1((M_\infty)_\fw,E[2^\infty])}
{E^\star_{\infty,\fw}}
}^\Delta
\end{align}
is finite, where $\fw$ is the prime of $M_\infty$ above $v$.
\end{lem}

\begin{proof}
Let $\fp$ be the prime of $M$ below $v$. By Lemma \ref{lem:2-primary_torsion},
the inflation-restriction sequence gives
\begin{align}\label{isom:cyclotomic_local}
H^1((M^{\cyc})_v,E[2^\infty])\cong H^1((M_\infty)_\fw,E[2^\infty])^\Delta.
\end{align}
The $\Delta$-cohomology of the exact sequence 
\[
0\longrightarrow E^\star_{\infty,\fw}
\longrightarrow H^1((M_\infty)_{\fw},E[2^\infty])
\longrightarrow
\frac{H^1((M_\infty)_{\fw},E[2^\infty])}
{E^\star_{\infty,\fw}}
\longrightarrow 0
\]
and \eqref{isom:cyclotomic_local} imply that the cokernel of \eqref{eq:local_chromatic_descent} is contained in $H^1(\Delta,E^\star_{\infty,\fw})$.
By local Tate duality and Definition
\ref{Definition: Local_condition},
$(E^\star_{\infty,\fw})^\vee
\cong\Imag\col^\star_{\fw}$.
Since $R_{M_\fp}$ is free over $\Z_2[\Delta]$,
the exact sequence
\begin{align}
0\to\Imag\col^\star_{\fw}\to R_{M_\fp}
\to R_{M_\fp}/\Imag\col^\star_{\fw}\to0
\end{align}
gives
\begin{align}
H^1(\Delta,E^\star_{\infty,\fw})^\vee\cong H_1(\Delta,\Imag\col^\star_{\fw})\cong H_2(\Delta,R_{M_\fp}/\Imag\col^\star_{\fw}).
\end{align}
By Propositions \ref{Prop:flat_condition} and
\ref{Prop:sharp_condition}, the quotient
$R_{M_\fp}/\Imag\col^\star_{\fw}$
is finitely generated over $\Z_2$.
The second homology group
$H_2(\Delta,R_{M_\fp}/\Imag\col^\star_{\fw})$
is therefore finitely generated over $\Z_2$
and annihilated by $\sharp\Delta=2$, hence finite.
\end{proof}

\begin{prop}\label{Prop:real_signed_localization}
Let $M=\Q$ or $F$, let
$\vec{s}=(s_\fp)_{\fp\in\Sigma_2(M)}$, and suppose that
$\cX^{\vec{s}}(E/M^{\cyc})$ is finitely generated over
$\Z_2$.
Put
\begin{align}
\cJ_M^{\vec{s}}:=
\cH^{\vec{s}}(M_\infty)^\Delta
\times\cH_\infty(M,E).
\end{align}
Then there is an exact sequence of $\Lambda$-modules
\begin{align}\label{eq:strict_signed_localization}
0\longrightarrow
\sel^{\vec{s}}(E/M^{\cyc})
\longrightarrow
H^1(M_\Sigma/M^{\cyc},E[2^\infty])
\longrightarrow
\cJ_M^{\vec{s}}
\longrightarrow C
\longrightarrow 0,
\end{align}
where $C$ is finite. Moreover, for every $j\geq2$,
localization induces an isomorphism
\begin{align}\label{eq:arch_higher_localization}
H^j(M_\Sigma/M^{\cyc},E[2^\infty])
\xrightarrow{\;\sim\;}
\varinjlim_n\bigoplus_{v\in\Sigma_\infty(M_n^{\cyc})}
H^j((M_n^{\cyc})_v,E[2^\infty]).
\end{align}

\end{prop}
\begin{proof}
We follow the strategy of the proof of
\cite[Proposition 2.8]{LL22}, taking into account
the contributions from real places and the descent
of the chromatic local conditions. We define
\begin{align}
\cJ_{M,\mathrm{loc}}^{\vec{s}}
:={}&
\prod_{v\in\Sigma_2(M^{\cyc})}
\Imag\parenth{
H^1((M^{\cyc})_v,E[2^\infty])
\longrightarrow
\cH_{\fw}(E/M_\infty)^\Delta
}\\
&\times
\prod_{\substack{v\in\Sigma(M^{\cyc})\\v\nmid2,\ v\nmid\infty}}
\frac{H^1((M^{\cyc})_v,E[2^\infty])}
{E((M^{\cyc})_v)\otimes_{\Z}\Q_2/\Z_2}
\times\cH_\infty(M,E),
\end{align}
where $\fw$ is the unique prime of $M_\infty$ above $v$
for each $v\mid2$.
By Lemma \ref{Lem:chromatic_descent_error},
we have
\begin{align}
\cJ_{M,\mathrm{loc}}^{\vec{s}}
\subseteq\cJ_M^{\vec{s}},
\qquad
\sharp\parenth{\cJ_M^{\vec{s}}/
\cJ_{M,\mathrm{loc}}^{\vec{s}}}<\infty.
\end{align}
We write
\begin{align}
H^1_{\mathrm{Iw}}(M^{\cyc}/M,\cT)
:=
\varprojlim_n H^1(M_\Sigma/M_n^{\cyc},\cT),
\end{align}
where the transition maps are corestrictions and
$M_n^{\cyc}$ is as defined in
Subsection \ref{subsec:doubly_Selmer}.
We define the compact Selmer group by
\begin{align}
\mathfrak S^{\vec{s}}(E/M^{\cyc})
:=
\Ker\parenth{
H^1_{\mathrm{Iw}}(M^{\cyc}/M,\cT)
\longrightarrow
\prod_{v\in\Sigma_2(M^{\cyc})}
\parenth{(E_{\infty,\fw}^{s_\fw})^\Delta}^\vee
},
\end{align}
where $\fw$ is the unique prime of $M_\infty$ above $v$.
The map is induced by localization and local Tate duality.
We apply the Poitou--Tate exact sequence
\cite[Theorem I.4.10]{Mil06} to $E[2^m]$ over
$M_n^{\cyc}$, using Tate cohomology at the real places.
For each prime $u\mid2$ of $M_n^{\cyc}$, we take
as the local condition the inverse image of
$(E_{\infty,\fw}^{s_\fw})^\Delta$ under the natural map
\begin{align}
H^1((M_n^{\cyc})_u,E[2^m])
\longrightarrow H^1((M^{\cyc})_v,E[2^\infty]),
\end{align}
via the identification \eqref{isom:cyclotomic_local},
where $v$ and $\fw$ are the primes of $M^{\cyc}$ and
$M_\infty$ above $u$, respectively.
At all other places, we take the zero local condition.
With these local conditions, the argument of
\cite[Proposition A.3.2]{Per95} gives an exact sequence
at each finite level. Taking the direct limit over
$n$ and $m$, we obtain
\begin{align}\label{eq:real_CPT}
\begin{tikzpicture}[baseline=(current bounding box.center)]
  \node[inner sep=0pt] (top) at (0,0) {$\displaystyle
    0\longrightarrow\sel^{\vec{s}}(E/M^{\cyc})
    \longrightarrow H^1(M_\Sigma/M^{\cyc},E[2^\infty])
    \longrightarrow \cJ_{M,\mathrm{loc}}^{\vec{s}}$};
  \node[inner sep=0pt] (bottom) at ($(top.south)+(0,-1.35cm)$) {$\displaystyle
    \fS^{\vec{s}}(E/M^{\cyc})^\vee
    \longrightarrow H^2(M_\Sigma/M^{\cyc},E[2^\infty])
    \longrightarrow \cR_M^2
    \longrightarrow 0,$};
  \coordinate (A) at ([xshift=.55em]top.east);
  \coordinate (E) at ([xshift=-.55em]bottom.west);
  \coordinate (M) at ($(A)!0.5!(E)$);
  \coordinate (R) at (A |- M);
  \coordinate (L) at (E |- M);
  \draw[->]
    (A)
      .. controls +(1.65cm,0) and +(1.65cm,0) .. (R)
      .. controls +(-2.0cm,0) and +(2.0cm,0) .. (L)
      .. controls +(-1.65cm,0) and +(-1.65cm,0) .. (E);
\end{tikzpicture}
\end{align}
where
\begin{align}
\cR_M^2:=\varinjlim_n
\bigoplus_{v\in\Sigma_\infty(M_n^{\cyc})}
H^2((M_n^{\cyc})_v,E[2^\infty]).
\end{align}
The compatibility of local Tate duality with restriction
and corestriction shows that the projective limit of
the dual Selmer groups is $\fS^{\vec{s}}(E/M^{\cyc})$.
By Lemma \ref{lem:2-primary_torsion}, we have
$E(M_n^{\cyc})[2^m]=\{O\}$.
Moreover, for every finite prime $u$ of $M_n^{\cyc}$,
the local Tate duality gives
\begin{align}
\varinjlim_m H^2((M_n^{\cyc})_u,E[2^m])=0,
\end{align}
since $E((M_n^{\cyc})_u)[2^\infty]$ is finite. Thus only the real places contribute to $\cR_M^2$. The global Euler characteristic formula
(cf. \cite[Proposition 3]{Gre89}) gives
\begin{align}\label{eq:real_global_coranks}
\corank_\Lambda H^1(M_\Sigma/M^{\cyc},E[2^\infty])
-\corank_\Lambda H^2(M_\Sigma/M^{\cyc},E[2^\infty])
=[M:\Q].
\end{align}
The chromatic local computations in
Propositions \ref{Prop:flat_condition} and
\ref{Prop:sharp_condition}, after taking the trivial
$\Delta$-component over $\Q_2$, give
\begin{align}\label{eq:real_local_coranks}
\corank_\Lambda\cJ_{M,\mathrm{loc}}^{\vec{s}}
=[M:\Q].
\end{align}
Here, the contributions from the real places and from finite
places not dividing $2$ have $\Lambda$-corank zero.
Since the cohomology groups at real places are annihilated by $2$,
the same holds for $\cR_M^2$. In particular,
$\corank_\Lambda\cR_M^2=0$.
Taking $\Lambda$-coranks in \eqref{eq:real_CPT},
using \eqref{eq:real_global_coranks} and
\eqref{eq:real_local_coranks}, we obtain
\begin{align}
\rank_\Lambda\mathfrak S^{\vec{s}}(E/M^{\cyc})
=\corank_\Lambda\sel^{\vec{s}}(E/M^{\cyc})=0.
\end{align}
By Lemma \ref{lem:2-primary_torsion}
and the five-term exact sequence of
\cite[Theorem 1]{Jan14}, we have
\begin{align}
H^1_{\mathrm{Iw}}(M^{\cyc}/M,\cT)\cong\Hom_\Lambda\parenth{H^1(M_\Sigma/M^{\cyc},E[2^\infty])^\vee,\Lambda}.
\end{align}
Hence $H^1_{\mathrm{Iw}}(M^{\cyc}/M,\cT)$ is $\Lambda$-torsion-free. Since $\mathfrak S^{\vec{s}}(E/M^{\cyc})$ is a $\Lambda$-torsion submodule of $H^1_{\mathrm{Iw}}(M^{\cyc}/M,\cT)$, it follows that $\mathfrak S^{\vec{s}}(E/M^{\cyc})=0$.
It follows from \eqref{eq:real_CPT} that the localization map
\begin{align}
H^1(M_\Sigma/M^{\cyc},E[2^\infty])\longrightarrow\cJ_{M,\mathrm{loc}}^{\vec{s}}
\end{align}
is surjective and that
\begin{align}\label{Isom:second_globa-local}
H^2(M_\Sigma/M^{\cyc},E[2^\infty])\cong\cR_M^2.
\end{align}
Since the localization map to $\cJ_M^{\vec{s}}$ factors through
$\cJ_{M,\mathrm{loc}}^{\vec{s}}$, we obtain an isomorphism
\begin{align}
\Coker\parenth{H^1(M_\Sigma/M^{\cyc},E[2^\infty])\longrightarrow\cJ_M^{\vec{s}}}\cong\cJ_M^{\vec{s}}/\cJ_{M,\mathrm{loc}}^{\vec{s}}.
\end{align}
Taking $C:=\cJ_M^{\vec{s}}/\cJ_{M,\mathrm{loc}}^{\vec{s}}$, which is finite by Lemma \ref{Lem:chromatic_descent_error},
we obtain the exact sequence \eqref{eq:strict_signed_localization}.
The second assertion follows from
\eqref{Isom:second_globa-local} for $j=2$
and from Poitou--Tate duality for $j\geq3$.
\end{proof}

\begin{cor}\label{Cor:archimedean_mu_correction}
Suppose that $\Delta_E$ is positive.
Let $M=\Q$ or $F$, and assume that
$\cX^{\vec{s}}(E/M^{\cyc})$ is $\Lambda$-torsion and $\mu(\cX^{\vec{s}}(E/M^{\cyc}))=0$. Then there is an exact sequence
\begin{align}\label{eq:archimedean_mu_correction}
0\longrightarrow
\parenth{\Lambda/(2)}^{\oplus[M:\Q]}
\longrightarrow
\cX_{\mathrm{rel}}^{\vec{s}}(E/M^{\cyc})
\longrightarrow
\cX^{\vec{s}}(E/M^{\cyc})
\longrightarrow0.
\end{align}
In particular,
\begin{align}
\mu\parenth{\cX_{\mathrm{rel}}^{\vec{s}}(E/M^{\cyc})}
=[M:\Q].
\end{align}
\end{cor}

\begin{proof}
By Definition \ref{Definition_Selmer_chromatic_cyclotomic},
the inverse image of
$\{0\}\times\{0\}\times\cH_\infty(M,E)$
under the localization map in \eqref{eq:real_CPT}
is precisely $\sel_{\mathrm{rel}}^{\vec{s}}(E/M^{\cyc})$.
Since this localization map is surjective by the proof of
Proposition \ref{Prop:real_signed_localization},
its restriction induces a surjection
\begin{align}
\loc_\infty:
\sel_{\mathrm{rel}}^{\vec{s}}(E/M^{\cyc})
\twoheadrightarrow\cH_\infty(M,E).
\end{align}
Since $\Delta_E>0$, we have $H^1(\R,E[2^\infty])\cong\Z/2\Z$.
For every $n\geq0$, each real place of $M$
splits completely in $M_n^{\cyc}/M$.
Consequently, after choosing a real place above each
real place of $M$, we have an isomorphism
of $\Z_2[\gal(M_n^{\cyc}/M)]$-modules
\begin{align}
\prod_{v\in\Sigma_\infty(M_n^{\cyc})}H^1((M_n^{\cyc})_v,E[2^\infty])
\cong\parenth{(\Z/2\Z)[\gal(M_n^{\cyc}/M)]}^{\oplus[M:\Q]}.
\end{align}
These isomorphisms are compatible with the restriction maps on the left-hand side and the natural projections on the right-hand side after taking Pontryagin duals. Hence
\begin{align}
\cH_\infty(M,E)^\vee
\cong\varprojlim_n\parenth{(\Z/2\Z)[\gal(M_n^{\cyc}/M)]}^{\oplus[M:\Q]}
\cong\parenth{\Lambda/(2)}^{\oplus[M:\Q]}.
\end{align}
Since $\loc_\infty$ is surjective, we have
\begin{align}\label{Isom:image_local_mod2_Lambda}
\Imag(\loc_\infty)^\vee
=\cH_\infty(M,E)^\vee
\cong\parenth{\Lambda/(2)}^{\oplus[M:\Q]}.
\end{align}
Substituting \eqref{Isom:image_local_mod2_Lambda} into \eqref{eq:archimedean_comparison}
gives \eqref{eq:archimedean_mu_correction}.
By the additivity of $\mu$ in
\eqref{eq:archimedean_mu_correction}, we obtain
\begin{align}
\mu\parenth{\cX_{\mathrm{rel}}^{\vec{s}}(E/M^{\cyc})}
=\mu\parenth{\parenth{\Lambda/(2)}^{\oplus[M:\Q]}}+\mu\parenth{\cX^{\vec{s}}(E/M^{\cyc})}
=[M:\Q].
\end{align}
\end{proof}

\begin{cor}\label{Cor:Sel_Herbrand}
Let $\star\in\{\sharp,\flat\}$ and assume that
$\sel^{\vec{s}_\star}(E/F^{\cyc})$ is $\Lambda$-cotorsion
with $\mu^{\vec{s}_\star}(E/F)=0$. Then
\begin{align}\label{eq:strict_Sel_Herbrand}
h_G\parenth{\sel^{\vec{s}_\star}(E/F^{\cyc})}
=
\frac{h_G\parenth{H^1(F_\Sigma/F^{\cyc},E[2^\infty])}}
{h_G\parenth{\cH^{\vec{s}_\star}(F_\infty)^\Delta}
\cdot h_G\parenth{\cH_\infty(F,E)}}.
\end{align}
\end{cor}
\begin{proof}
Applying the multiplicativity of Herbrand quotients
to \eqref{eq:strict_signed_localization} with $M=F$
and $\vec{s}=\vec{s}_\star$, and noting that $h_G(C)=1$
since $C$ is finite, we obtain \eqref{eq:strict_Sel_Herbrand}.
\end{proof}

\begin{prop}\label{Prop:global_Herbrand_trivial}
Let $\star\in\{\sharp,\flat\}$ and assume that $\sel^{\star}(E/\Q^{\cyc})$ is $\Lambda$-cotorsion and $\mu^\star(E/\Q)=0$. Then
\begin{align}
H^i\left(
G,H^1(F_\Sigma/F_\infty,E[2^\infty])^\Delta
\right)=0
\qquad(i=1,2).
\end{align}
In particular,
\begin{align}
h_G\left(
H^1(F_\Sigma/F_\infty,E[2^\infty])^\Delta
\right)=1.
\end{align}
\end{prop}

\begin{proof}
The assumptions involving $\star$ are used only to apply
Proposition~\ref{Prop:real_signed_localization}
over $M=\Q$ and $F$.
We adapt Matsuno's argument of \cite[Proof of Proposition 3.5]{Mat08}
to our setting.
By Corollary \ref{Cor:Kida_Herbrand},
$\mu^{\vec{s}_\star}(E/F)=0$.
Apply Proposition \ref{Prop:real_signed_localization}
over $M=\Q$ and $F$. For $j\geq2$, it identifies
$H^j(M_\Sigma/M^{\cyc},E[2^\infty])$ with
\begin{align}
\cR_M^j:=\varinjlim_n
\bigoplus_{v\in\Sigma_\infty(M_n^{\cyc})}
H^j((M_n^{\cyc})_v,E[2^\infty]).
\end{align}
Since every real place of $\Q^{\cyc}$ splits completely
in $F^{\cyc}$, we have
\begin{align}\label{eq:induced_real_cohomology_new}
\cR_F^j\cong\cR_\Q^j\otimes_{\Z_2}\Z_2[G]
\qquad(j\geq2).
\end{align}
In particular, $\cR_F^j$ is cohomologically trivial
as a $G$-module, and the restriction map
$\cR_\Q^j\to(\cR_F^j)^G$ is an isomorphism.
Consider the Hochschild--Serre spectral sequence
\begin{align}
H^p\parenth{G,H^q(F_\Sigma/F^{\cyc},E[2^\infty])}
\Longrightarrow H^{p+q}(\Q_\Sigma/\Q^{\cyc},E[2^\infty]).
\end{align}
Write $E_r^{p,q}$ for the $(p,q)$-term on the
$r$-th page of this spectral sequence. By Lemma \ref{lem:2-primary_torsion},
\eqref{eq:induced_real_cohomology_new}, and Shapiro's lemma,
we have
\begin{align}
E_2^{p,0}&=0 \qquad(p\geq0),\\
E_2^{p,q}&=0 \qquad(p>0,\ q\geq2).
\end{align}
Furthermore, put
$H^n:=H^n(\Q_\Sigma/\Q^{\cyc},E[2^\infty])$,
and let $F^pH^n$ denote the decreasing filtration
on $H^n$ associated with the spectral sequence.
By Proposition \ref{Prop:real_signed_localization}
and \eqref{eq:induced_real_cohomology_new},
the restriction map induces an isomorphism
\begin{align}
\rest:H^n\xrightarrow{\sim}
H^n(F_\Sigma/F^{\cyc},E[2^\infty])^G
=E_2^{0,n}
\qquad(n\geq2).
\end{align}
This is the edge homomorphism, which factors as
\begin{align}
H^n\twoheadrightarrow H^n/F^1H^n
\cong E_\infty^{0,n}
\hookrightarrow E_2^{0,n}.
\end{align}
The last map is injective since no differential
enters the column $p=0$ in this first-quadrant
spectral sequence.
Since the composite is an isomorphism for $n\geq2$,
its kernel $F^1H^n$ vanishes.
As the filtration is decreasing, we obtain
\begin{align}
F^pH^n=0
\qquad(p\geq1,\ n\geq2).
\end{align}
By the convergence of the spectral sequence,
it follows that
\begin{align}
E_\infty^{i,1}
\cong F^iH^{i+1}/F^{i+1}H^{i+1}
=0
\qquad(i=1,2).
\end{align}
We next show that $E_2^{i,1}=E_\infty^{i,1}$
for $i=1,2$.
Since the edge homomorphism
$H^2\to E_2^{0,2}$ is an isomorphism,
the natural inclusion
$E_\infty^{0,2}\hookrightarrow E_2^{0,2}$
is also an isomorphism.
Consequently, $d_2^{0,2}=0$.
On the other hand, the outgoing differentials
from $E_2^{i,1}$ have targets
$E_2^{i+2,0}=0$ by Lemma
\ref{lem:2-primary_torsion},
while those on subsequent pages have targets
with negative second index.
Thus no nonzero differential enters or leaves
$E_r^{i,1}$ for $i=1,2$ and $r\geq2$.
Therefore
\begin{align}\label{vanishing_cohomology}
H^i\parenth{G,H^1(F_\Sigma/F^{\cyc},E[2^\infty])}
=E_2^{i,1}
\cong E_\infty^{i,1}
=0
\qquad(i=1,2).
\end{align}
Since Lemma \ref{lem:2-primary_torsion}
gives $E(F_\infty)[2^\infty]=0$,
by the inflation--restriction sequence and \eqref{vanishing_cohomology}, we have
\begin{align}
H^i(G,H^1(F_\Sigma/F_\infty,E[2^\infty])^\Delta)\cong H^i(G,H^1(F_\Sigma/F^{\cyc},E[2^\infty]))
=0
\end{align}
for $i\in\{1,2\}$.
\end{proof}

\begin{rem}
In the good ordinary case, Matsuno
\cite[Proposition 3.5]{Mat08} proved the corresponding
isomorphisms, which, in our notation, take the form
\begin{align}
H^i(G,H^1(F_\Sigma/F^{\cyc},E[2^\infty]))
\cong H^i(G,E(F^{\cyc})[2^\infty])
\qquad(i=1,2).
\end{align}
In our supersingular setting, Lemma
\ref{lem:2-primary_torsion} gives
$E(F^{\cyc})[2^\infty]=\{O\}$, and we obtain
\begin{align}
H^i(G,H^1(F_\Sigma/F^{\cyc},E[2^\infty]))=0
\qquad(i=1,2).
\end{align}
Moreover, the Hochschild--Serre spectral sequence
degenerates at $E_2$, rather than at $E_3$
as in Matsuno's proof.
\end{rem}

\begin{lem}\label{lem:local_Herbrand_at_two}
Let $F/\Q$ be a real quadratic extension in which $2$ is
unramified, and let $\star\in\{\sharp,\flat\}$. Then
\begin{align}
h_G\left(
\left(
\prod_{\fw\in\Sigma_2(F_\infty)}
\cH_{\fw}(E/F_\infty)
\right)^\Delta
\right)=1.
\end{align}
\end{lem}

\begin{proof}
If $2$ splits in $F/\Q$, the two local factors are
interchanged by $G$. Their product, also after taking
$\Delta$-invariants, is an induced $G$-module; hence its
Herbrand quotient is $1$ by Shapiro's lemma. Suppose that $2$ is inert in $F/\Q$, and put
$K:=F_{\fp}$. Since $K/\Q_2$ is unramified
and $\Q_{2,\infty}/\Q_2$ is totally ramified,
restriction induces a natural isomorphism
\begin{align}
G_{\mathrm{ur}}=\gal(K_\infty/\Q_{2,\infty})
\cong\gal(K/\Q_2)\cong G.
\end{align}
We use these identifications throughout the proof.
We first show that $H^1_{\iw,K}(\cT)$ is free of rank $2$ over $R_K$.
Lemma \ref{lem:2-primary_torsion} and local Tate
duality give $H^0(K_n,\cT)=H^2(K_n,\cT)=0$. The Hochschild--Serre spectral sequence
and the above vanishing give
\begin{align}
H^i(G_n,H^1(K_n,\cT))
\simeq H^{i+1}(\Q_2,\cT)=0
\qquad(i=1,2).
\end{align}
Since $G_n$ is a finite $2$-group,
the cohomological triviality criterion
(cf. \cite[Chapter VI, Theorem 8.7]{Bro94}),
together with Nakayama's lemma, shows
that $H^1(K_n,\cT)$ is free over
$\Z_2[G_n]$.
Its rank is $2$ by the local Euler
characteristic formula. Put $H_n:=\gal(K_{n+1}/K_n)$.
Since $E(K_{n+1})[2^\infty]=\{O\}$,
the inflation-restriction sequence gives
\begin{align}\label{res_isom}
H^1(K_n,E[2^\infty])\overset{\sim}{\longrightarrow}H^1(K_{n+1},E[2^\infty])^{H_n}.
\end{align}
By the local Tate duality, the dual of \eqref{res_isom} is the
corestriction isomorphism
\begin{align}
H^1(K_{n+1},\cT)_{H_n}
\overset{\sim}{\longrightarrow}
H^1(K_n,\cT).
\end{align}
Together with
$\Z_2[G_{n+1}]_{H_n}\cong\Z_2[G_n]$, Nakayama's lemma allows us to choose free bases compatible with corestriction at every finite level. Taking the inverse limit, we obtain
\begin{align}
H^1_{\iw,K}(\cT)\cong\varprojlim_n\Z_2[G_n]^{\oplus2}\cong R_K^{\oplus2}.
\end{align}
By the definition of the $\star$-Coleman map, we have the exact sequence of $R_K$-modules
\begin{align}\label{Col_def_ext_seq}
0\longrightarrow \Ker(\col_K^\star)\longrightarrow H^1_{\iw,K}(\cT)
\longrightarrow \Imag(\col_K^\star)\longrightarrow0.
\end{align}
By Lemma \ref{Lem:ker_finite},
$\Imag(\col_K^\star)$ is induced as a
$G_{\mathrm{ur}}$-module, compatibly with $\Delta$.
Since $H^1_{\iw,K}(\cT)\cong R_K^{\oplus2}$
and
$R_K\cong\Z_2[G_{\mathrm{ur}}]\otimes_{\Z_2}R_{\Q_2}$,
the same holds for $H^1_{\iw,K}(\cT)$,
which is also free over $\Z_2[\Delta]$.
Taking $\Delta$-homology of \eqref{Col_def_ext_seq} yields
\begin{align}\label{seq:homology_def_Coleman}
0\longrightarrow H_1(\Delta,\Imag(\col_K^\star))
\longrightarrow \Ker(\col_K^\star)_\Delta
\longrightarrow H^1_{\iw,K}(\cT)_\Delta
\longrightarrow \Imag(\col_K^\star)_\Delta
\longrightarrow0.
\end{align}
The terms $H_1(\Delta,\Imag(\col_K^\star))$, $H^1_{\iw,K}(\cT)_{\Delta}$, and $\Imag(\col_K^\star)_\Delta$ in \eqref{seq:homology_def_Coleman} are induced
$G_{\mathrm{ur}}$-modules, since the $G_{\mathrm{ur}}$-factor is independent
of the $\Delta$-action.
The long exact Tate cohomology sequences and
Shapiro's lemma therefore show that $\Ker(\col_K^\star)_\Delta$
is cohomologically trivial for $G_{\mathrm{ur}}$.
Thus, the local Tate duality gives
\begin{align}
\left(\cH_{\fp}(E/F_\infty)^\Delta\right)^\vee
\cong \Ker(\col_K^\star)_\Delta.
\end{align}
By Pontryagin duality, the required Herbrand quotient
is $1$.
\end{proof}

\begin{thm}\label{Thm:Matsuno-type_formula}
Let $E$ be an elliptic curve defined over $\Q$ with good supersingular reduction at $2$. Let $\star\in\{\sharp,\flat\}$ and let $F$ be a real quadratic field in which $2$ is unramified. For each finite prime $\fw\nmid 2$ of $F^{\cyc}$, put
\begin{align}\label{eq:modified_ramification_index}
e'_E(\fw):=
\begin{cases}
0&\text{if $E$ has additive reduction at $v$},\\
e_{F^{\cyc}/\Q^{\cyc}}(\fw)&\text{otherwise},
\end{cases}
\end{align}
where $e_{F^{\cyc}/\Q^{\cyc}}(\fw)$ is the ramification index of $\fw$ in $F^{\cyc}/\Q^{\cyc}$ and $v$ is the prime of $\Q^{\cyc}$ below $\fw$. If $\sel^{\star}(E/\Q^{\cyc})$ is $\Lambda$-cotorsion and $\mu^{\star}(E/\Q)=0$, then $\sel^{\vec{s}_{\star}}(E/F^{\cyc})$ is also $\Lambda$-cotorsion and $\mu^{\vec{s}_{\star}}(E/F)=0$. Moreover, we have
\begin{align}
\lambda^{\vec{s}_\star}(E/F)=2\lambda^{\star}(E/\Q)+\sum_{\fw\in P_1}(e'_E(\fw)-1)+2\sum_{\fw\in P_2}(e'_E(\fw)-1),
\end{align}
where $P_1$ and $P_2$ are the sets of primes of $F^{\cyc}$ defined by
\begin{align}
&P_1:=\{\fw\ |\ \fw\nmid 2,\ E\text{ has split multiplicative reduction at }\fw \},\\
&P_2:=\{\fw\ |\ \fw\nmid 2,\ E\text{ has good reduction at }\fw,\ E[2]((F^{\cyc})_{\fw})\neq\{O\}\}.
\end{align}
\end{thm}

\begin{proof}
Since $F/\Q$ is real quadratic,
$\cH_\infty(F,E)$ is induced as a $G$-module.
Thus Shapiro's lemma gives $h_G\parenth{\cH_\infty(F,E)}=1$.
Together with Lemma \ref{lem:local_Herbrand_at_two}, Corollaries \ref{Cor:Kida_Herbrand} and
\ref{Cor:Sel_Herbrand}, and Proposition
\ref{Prop:global_Herbrand_trivial}, we obtain
\begin{align}\label{eq:lambda_local_Herbrand}
&\lambda^{\vec{s}_\star}(E/F)-2\lambda^\star(E/\Q)
\notag\\
&\quad=-v_2\left(
h_G\left(
\left(\prod_{\fw\in\Sigma(F_\infty)\backslash \Sigma_2(F_{\infty})}
\cH_{\fw}(E/F_\infty)\right)^\Delta
\right)\right).
\end{align}
For every finite prime $v\in \Sigma(\Q^{\cyc})$ with $v\nmid2$,
choose a prime $u$ of $F^{\cyc}$ above $v$, and set
\begin{align}
G_u:=\gal((F^{\cyc})_u/(\Q^{\cyc})_v),\qquad
B_u:=\left(
\prod_{\substack{\fw\in\Sigma(F_\infty)\\\fw\mid u}}
\cH_{\fw}(E/F_\infty)
\right)^\Delta.
\end{align}
Let $\ell$ be the odd rational prime below $v$. The extension $(\Q^{\cyc})_v/\Q_\ell$ is the
unramified $\Z_2$-extension, while $\Q_\ell(i)/\Q_\ell$ is either the trivial or the unramified quadratic extension. Hence $i\in(\Q^{\cyc})_v\subseteq(F^{\cyc})_u$, and $u$ splits completely in $F_\infty/F^{\cyc}$. As $\Delta$ interchanges the two primes
above $u$, taking $\Delta$-invariants
identifies $B_u$ with the usual local term
\begin{align}
B_u\cong
H^1((F^{\cyc})_u,E[2^\infty])
\cong H^1((F^{\cyc})_u,E)[2^\infty].
\end{align}
Grouping the primes of $F^{\cyc}$
into $G$-orbits, we obtain
\begin{align}
\left(
\prod_{\substack{\fw\in\Sigma(F_\infty)\\\fw\nmid2\infty}}
\cH_{\fw}(E/F_\infty)\right)^\Delta\cong\prod_{\substack{v\in\Sigma(\Q^{\cyc})\\v\nmid2\infty}}
\operatorname{Ind}_{G_u}^G B_u.
\end{align}
Therefore, by Shapiro's lemma and multiplicativity,
\begin{align}\label{eq:local_orbit_decomposition}
v_2\left(
h_G\left(
\left(
\prod_{\substack{\fw\in\Sigma(F_\infty)\\\fw\nmid2\infty}}
\cH_{\fw}(E/F_\infty)\right)^\Delta\right)\right)=\sum_{\substack{v\in\Sigma(\Q^{\cyc})\\v\nmid2\infty}}
v_2\left(h_{G_u}(B_u)\right).
\end{align}
If $G_u$ is trivial, then $h_{G_u}(B_u)=1$.
Otherwise $G_u=G$ and the local extension
$(F^{\cyc})_u/(\Q^{\cyc})_v$ is totally ramified
of degree $2$.
As in \cite[Lemma 4.2]{HM99}, the vanishing of
$H^j((\Q^{\cyc})_v,E[2^\infty])$ and $H^j((F^{\cyc})_u,E[2^\infty])$
for $j\geq2$, together with the Hochschild--Serre
spectral sequence, gives
\begin{align}\label{eq:local_HS_Herbrand}
H^i(G_u,B_u)
\cong H^i\bigl(G_u,E((F^{\cyc})_u)[2^\infty]\bigr)
\qquad(i=1,2).
\end{align}
Here $E((F^{\cyc})_u)\otimes_{\Z}\Q_2/\Z_2=0$, since $(F^{\cyc})_u$ is
an unramified $\Z_2$-extension of a finite extension of $\Q_\ell$ with $\ell\neq2$. 
The local calculation follows from
\cite[Proposition 5.1 and Corollary 5.2]{HM99},
with the correction for changes in reduction
type as in \cite[Theorem 3.1]{Mat08}.
More precisely, if $G_u=G$, then
\begin{align}
v_2(h_{G_u}(B_u))=
\begin{cases}
-(e'_E(u)-1) & \text{if }u\in P_1,\\
-2(e'_E(u)-1) & \text{if }u\in P_2,\\
0 & \text{otherwise}.
\end{cases}
\end{align}
Since $e_{F^{\cyc}/\Q^{\cyc}}(u)=2$ precisely
in the nontrivial local case, this can be written
uniformly as
\begin{align}
v_2\left(h_{G_u}(B_u)\right)
=-\mathbf 1_{P_1}(u)\bigl(e'_E(u)-1\bigr)
-2\cdot \mathbf 1_{P_2}(u)\bigl(e'_E(u)-1\bigr).
\end{align}
Here, $\mathbf{1}_{P_i}$ denotes the indicator
function of $P_i$ for $i=1,2$. Substituting into \eqref{eq:local_orbit_decomposition}
and then into \eqref{eq:lambda_local_Herbrand},
and reindexing by primes of $F^{\cyc}$, yields
\begin{align}
\lambda^{\vec{s}_\star}(E/F)-2\lambda^\star(E/\Q)
=\sum_{\fw\in P_1}\bigl(e'_E(\fw)-1\bigr)
+2\sum_{\fw\in P_2}\bigl(e'_E(\fw)-1\bigr),
\end{align}
as required.
\end{proof}

\subsection{A Matsuno-type difference formula in the supersingular case}\label{subsec:Matsuno_supersingular}
Using the $2$-adic Kida-type formula obtained in \S\ref{sec:Kida}, we derive a supersingular analogue of Matsuno's formula \cite[Theorem 5.1]{Mat08}. This is an algebraic counterpart of the formula \cite[Theorem 1.3]{Ada26+} for the chromatic $2$-adic $L$-function.

\begin{thm}\label{chromatic_Matsuno-type_formula}
Let $E$ be an elliptic curve defined over $\Q$ with good supersingular reduction at $2$ and let $D>0$ be a square-free integer such that $D\equiv 1\pmod 4$. Fix $\star\in\{\sharp,\flat\}$. Assume that the conductor $N_E$ of $E$ is square-free and that $\sel^{\star}(E/\Q^{\cyc})$ is $\Lambda$-cotorsion. If $\mu^{\star}(E/\Q)=0$, then $\mu^{\star}(E^D/\Q)=0$ and we have
\begin{align}\label{eq:Matsuno-type_formula}
\lambda^{\star}(E^D/\Q)=\lambda^{\star}(E/\Q)+\sum_{\substack{\ell\mid D\\ \ell\mid N_E}}2^{n_{\ell}}+\sum_{\substack{\ell\mid D\\ \ell\nmid N_E, 2\mid \sharp \widetilde{E}_{\ell}(\F_{\ell})}}2^{n_{\ell}+1}.
\end{align}
\end{thm}

\begin{proof}
The case $D=1$ is trivial. Hence, we assume
$D>1$ and put $F:=\Q(\sqrt{D})$, and fix an isomorphism
$\psi:E\to E^D$ over $F$.
We first establish the compatibility of the chromatic local
conditions under $\psi$. 
If $D\equiv1\pmod8$, then $2$ splits in $F/\Q$.
At either prime $\fp\mid2$, we have $F_\fp=\Q_2$ and
$a_2(E^D)=a_2(E)$. The two systems of local points constructed in
Subsection \ref{subsec:local_points} therefore have the same
formal logarithms. Since $\psi$ extends to an isomorphism of
the formal groups over $\Z_2$, it identifies these points up
to a fixed scalar in $\Z_2^\times$. By the construction in
Subsection \ref{subsec:Coleman}, the corresponding Coleman
maps differ by multiplication by the same unit. Hence, for
each $\star\in\{\sharp,\flat\}$, we obtain
\begin{align}
\psi_*\bigl(\Ker\col_E^\star\bigr)
=\Ker\col_{E^D}^\star.
\end{align}

Suppose now that $D\equiv5\pmod8$, so that $2$ is inert.
Put $K:=F_\fp$ and let
$\varphi\in\gal(K/\Q_2)$ be the nontrivial Frobenius.
We put $t:=\zeta-\zeta^2$.
Then 
\begin{align}
t^2=-3,\quad t^\varphi=-t.
\end{align}
Since $D/(-3)\equiv1\pmod8$, choose $b\in\Z_2^\times$ such that
$\sqrt D=bt$. Write $a:=a_2(E)$, so that $a_2(E^D)=-a$.
For $n\geq-2$, denote by $d_{a,n}$ and $d_{-a,n}$ the local
points of $E$ and $E^D$, respectively, constructed over
$K_n=K(\zeta_{2^{n+2}})$ using the same $\zeta$.

Let $x_j(*)$ be the sequence of Definition \ref{Def:sequence} corresponding to $*\in\{a,-a\}$.
Its recurrence gives
\begin{align}\label{eq:twist_xj_sign}
x_j(-a)=(-1)^j x_j(a)\qquad(j\geq0).
\end{align}
Put $u_n:=\zeta_{2^{n+2}}$ and
$\alpha_{n,j}:=(\zeta^{\varphi^{-n}})^{2^j}$.
As $\varphi(\zeta)=\zeta^2$, a direct computation yields
\begin{align}\label{eq:twist_alpha_h}
t\alpha_{n,j}=(-1)^{n+j}(1+2\varphi)(\alpha_{n,j})
\qquad(n\geq-2,\ j\geq0).
\end{align}
Let $\eta_{a,n}$ and $\eta_{-a,n}$ denote the elements in
Definition \ref{Def:eta_system} for $E$ and $E^D$.
By Remark \ref{Rem:phi_relation_eta_n}, these elements satisfy
\begin{align}
(3-a\varphi)\eta_{a,n}=2\alpha_{n,0},\qquad
(3+a\varphi)\eta_{-a,n}=2\alpha_{n,0}.
\end{align}
Using $t^\varphi=-t$ and \eqref{eq:twist_alpha_h}, we obtain
\begin{align}
(3-a\varphi)(t\eta_{-a,n})
=2t\alpha_{n,0}
=(3-a\varphi)\bigl((-1)^n(1+2\varphi)\eta_{a,n}\bigr).
\end{align}
The factor $3-a\varphi$ is a unit in $\Z_2[\varphi]$, hence
\begin{align}\label{eq:twist_eta_h}
t\eta_{-a,n}=(-1)^n(1+2\varphi)\eta_{a,n}.
\end{align}
The formula for the formal logarithms in
\eqref{eq:log_dn}, together with
\eqref{eq:twist_xj_sign}, \eqref{eq:twist_alpha_h},
and \eqref{eq:twist_eta_h}, now gives
\begin{align}\label{eq:twist_log_dn}
t\log_{\what{E^D}}(d_{-a,n})
=(-1)^n(1+2\varphi)\log_{\what E}(d_{a,n}).
\end{align}
Here $\varphi$ acts on the coefficients in $K$ and fixes
$u_n$, so that \eqref{eq:twist_alpha_h} applies termwise to
the series in \eqref{eq:log_dn}.
Let $\omega_E$ and $\omega_{E^D}$ be the invariant
differentials associated with the chosen minimal models.
Since the quadratic twist is unramified at $2$, there exists
$c\in\Z_2^\times$ with
$\psi^*\omega_{E^D}=c\omega_E/\sqrt D$.
The compatibility of formal logarithms with $\psi$ implies
\begin{align}
\log_{\what E}(\psi^{-1}(d_{-a,n}))
=\frac{\sqrt D}{c}\log_{\what{E^D}}(d_{-a,n}).
\end{align}
Putting $\varepsilon:=b/c\in\Z_2^\times$ and applying
\eqref{eq:twist_log_dn}, we obtain
\begin{align}
\log_{\what E}(\psi^{-1}(d_{-a,n}))
=\varepsilon(-1)^n(1+2\varphi)\log_{\what E}(d_{a,n}).
\end{align}
By Lemma \ref{injective_log}, the formal logarithm is
injective. Therefore
\begin{align}\label{eq:twist_local_points}
\psi^{-1}(d_{-a,n})
=\varepsilon(-1)^n(1+2\varphi)d_{a,n}
\qquad(n\geq-2).
\end{align}

We next compare the Coleman maps. Let $\chi$ be the
unramified quadratic character of
$\gal(K_\infty/\Q_2)$, and define the ring automorphism $\operatorname{Tw}_\chi:R_K\to R_K$ by $\operatorname{Tw}_{\chi}(g):=\chi(g)g$ for each $g\in \gal(K_{\infty}/\Q_2)$.
In particular, $\chi(\varphi)=-1$ and
$\chi(\gamma)=1$, so
$\operatorname{Tw}_\chi(1+2\varphi)=1-2\varphi$.
For $n\geq0$, let
\begin{align}
P_{a,n}(z_n):=
\left(\sum_{g\in G_n}(gd_{a,n},z_n)_n g,\ 
\sum_{g\in G_n}(gd_{a,n-1},z_n)_n g\right),
\end{align}
and define $P_{-a,n}$ similarly for $E^D$.
Let $g_E$ and $g_{E^D}$ denote the natural actions of $g\in G_n$ on $E(K_n)$ and $E^D(K_n)$, respectively. Then we have
\begin{align}
g(\psi)=[\chi(g)]\circ\psi,
\qquad
\chi(g)=\frac{g(\sqrt D)}{\sqrt D}.
\end{align}
Consequently, for $g\in G_n$,
\begin{align}
\psi^{-1}\circ g_{E^D}\circ\psi
=[\chi(g)]\circ g_E.
\end{align}
Moreover, $(1+2\varphi)^\iota=1+2\varphi$, and right multiplication by
$(1+2\varphi)^\iota$ describes the effect of acting by $1+2\varphi$ on
the local points in the pairing vector. The functoriality
of the local Tate pairing and
\eqref{eq:twist_local_points} therefore give
\begin{align}\label{eq:twist_pairing_vector}
P_{-a,n}(\psi_*z_n)
=\varepsilon(-1)^n(1-2\varphi)
\operatorname{Tw}_\chi(P_{a,n}(z_n))S,
\end{align}
where $S:=\begin{pmatrix}1&0\\0&-1\end{pmatrix}$.
For $*\in\{a,-a\}$, let $A_n(*)$ and $\cB_n(*)$ be the matrices of
Subsection \ref{subsec:Coleman},
and let $\widetilde{\cB}_n(*)$ denote
the adjugate matrix of $\cB_n(*)$.
Since $A_n(-a)=-SA_n(a)S$, we have $\cB_n(-a)=(-1)^nS\cB_n(a)S$. Thus, the corresponding adjugate matrices satisfy
\begin{align}\label{eq:twist_adjugate}
\widetilde{\cB}_n(-a)
=(-1)^nS\widetilde{\cB}_n(a)S.
\end{align}
Recall from \cite[(4.7)]{Kat21} and the construction of
$\col'_n$ in Subsection \ref{subsec:Coleman} that
$P_{a,n}(z_n)=\col'_{a,n}(z_n)
\widetilde{\cB}_n(a)^\iota$, where this identity is read
via the injective map on the finite-level quotient of
\cite[Lemma 4.9 (2)]{Kat21}.
Since $\operatorname{Tw}_\chi$ fixes $\gamma$, it fixes
the entries of $\widetilde{\cB}_n(a)^\iota$.
Combining \eqref{eq:twist_pairing_vector} with
\eqref{eq:twist_adjugate} and using the same injection,
we obtain
\begin{align}\label{eq:finite-layer_twist}
\col'_{-a,n}(\psi_*z_n)
=\varepsilon(1-2\varphi)
\operatorname{Tw}_\chi(\col'_{a,n}(z_n))S.
\end{align}

For $*\in\{a,-a\}$, put
\begin{align}
U(*):=\begin{pmatrix}
-*& -1\\
1 & 0\end{pmatrix}.
\end{align}
Then
\begin{align}
SU(-a)=U(a)
\begin{pmatrix}
-1 & 0\\
0 & 1
\end{pmatrix}.
\end{align}
Taking projective limits of \eqref{eq:finite-layer_twist} and using
Definition \ref{Definition: Coleman_map}, we obtain

\begin{align}\label{eq:twist_col_compatibility}
\bigl(\col_{E^D}^\sharp(\psi_*z),
\col_{E^D}^\flat(\psi_*z)\bigr)
=\varepsilon(1-2\varphi)
\operatorname{Tw}_\chi\!\left(
-\col_E^\sharp(z),\col_E^\flat(z)\right).
\end{align}
The factor $\varepsilon(1-2\varphi)$ is a unit in $R_K$,
and $\operatorname{Tw}_\chi$ is an automorphism.
It follows that
$\psi_*(\Ker\col_E^\star)=\Ker\col_{E^D}^\star$.
By the functoriality of the local Tate pairing and
Definition \ref{Definition: Local_condition}, we have
\begin{align}\label{eq:twist_local_condition_compatibility}
\psi_*\bigl(E^\star_{\infty,\fp}\bigr)=E^{D,\star}_{\infty,\fp}.
\end{align}

At primes not lying above $2$, the isomorphism $\psi$ induces the isomorphisms on the corresponding local cohomology groups. 
It follows that $\psi$ induces a $\gal(F_\infty/F)$-equivariant isomorphism
\begin{align}
\sel^{\vec{s}_{\star}}(E/F_\infty)\cong\sel^{\vec{s}_{\star}}(E^D/F_\infty).
\end{align}

For every real place $v$ of $F^{\cyc}$, the isomorphism
$\psi:E\to E^D$ over $F$ induces an isomorphism
\begin{align}
\psi_{*,v}:
H^1((F^{\cyc})_v,E[2^\infty])
\xrightarrow{\sim}
H^1((F^{\cyc})_v,E^D[2^\infty]),
\end{align}
compatible with the archimedean localization maps.
Consequently, taking $\Delta$-invariants of the
preceding isomorphism over $F_\infty$ and imposing
the real local conditions yields
\begin{align}\label{eq:twist_signed_selmer_compatibility}
\sel^{\vec{s}_{\star}}(E/F^{\cyc})
\cong\sel^{\vec{s}_{\star}}(E^D/F^{\cyc}).
\end{align}
Next, we compute the $\mu$-invariants. Put $G:=\gal(F^{\cyc}/\Q^{\cyc})$.
By Theorem \ref{Thm:Matsuno-type_formula} applied to $E$,
$\sel^{\vec{s}_\star}(E/F^{\cyc})$ is
$\Lambda$-cotorsion with $\mu^{\vec{s}_\star}(E/F)=0$.
The isomorphism \eqref{eq:twist_signed_selmer_compatibility}
therefore gives
\begin{align}\label{eq:twist_F_invariants}
\mu^{\vec{s}_\star}(E^D/F)=0,\qquad
\lambda^{\vec{s}_\star}(E^D/F)
=\lambda^{\vec{s}_\star}(E/F).
\end{align}
By Corollary \ref{Cor:Selmer_control} applied to $E^D$,
the restriction map
\begin{align}
\sel^\star(E^D/\Q^{\cyc})\longrightarrow
\sel^{\vec{s}_\star}(E^D/F^{\cyc})^G
\end{align}
has finite kernel and cokernel. Since $(\sel^{\vec{s}_\star}(E^D/F^{\cyc})^G)^{\vee}$ is a quotient of the finitely generated
$\Z_2$-module
$\cX^{\vec{s}_\star}(E^D/F^{\cyc})$, $\cX^{\star}(E^D/\Q^{\cyc})$ is also finitely
generated over $\Z_2$. Hence it is $\Lambda$-torsion and
$\mu^\star(E^D/\Q)=0$.

Then we can apply Theorem \ref{Thm:Matsuno-type_formula}
to both $E$ and $E^D$.
Since $E$ and $E^D$ are isomorphic over $F$, the corresponding sets $P_1$ and $P_2$ coincide.
Subtracting the two Kida-type formulae and using
\eqref{eq:twist_F_invariants}, we obtain
\begin{align}\label{eq:twist_kida_subtract}
&2\lambda^\star(E^D/\Q)-2\lambda^\star(E/\Q)\notag\\
&=\sum_{\fw\in P_1}\bigl(e'_E(\fw)-e'_{E^D}(\fw)\bigr)
+2\sum_{\fw\in P_2}\bigl(e'_E(\fw)-e'_{E^D}(\fw)\bigr).
\end{align}
Since $e'_E(\fw)=e'_{E^D}(\fw)$ for $\ell\nmid D$, only the primes $\ell\mid D$ contribute to the right-hand side. For every such prime, $F^{\cyc}/\Q^{\cyc}$ is totally
ramified of degree $2$ at the primes above $\ell$,
and the number of these primes is $2^{n_\ell}$
by \cite[Remark 3.21]{Ada26+}.

If $\ell\mid N_E$, then $E$ has multiplicative reduction
and $E^D$ has additive reduction at $\ell$.
Let $v$ be the prime of $\Q^{\cyc}$ below $\fw$.
Since $(\Q^{\cyc})_v/\Q_\ell$ is the unramified
$\Z_2$-extension, $(\Q^{\cyc})_v$ contains the unramified quadratic
extension of $\Q_\ell$. Hence $E$ has split multiplicative
reduction at every prime $\fw$ above $\ell$
(cf. \cite[Lemma 5.3 (ii)]{Mat08}).
Thus these primes belong to $P_1$, and each contributes
$e'_E(\fw)-e'_{E^D}(\fw)=2$.

If $\ell\nmid N_E$, then $E$ has good reduction
and $E^D$ has additive reduction at $\ell$.
Since $\gal((F^{\cyc})_{\fw}/\Q_\ell)$ is a pro-$2$ group
and $E[2]$ is unramified at $\ell$,
the corresponding primes belong to $P_2$ precisely when
$2\mid\sharp\widetilde E_\ell(\F_\ell)$,
by \cite[Lemma 5.3 (i)]{Mat08}.
In that case, each contributes
$2(e'_E(\fw)-e'_{E^D}(\fw))=4$.
These local calculations coincide with
\cite[Lemmas 5.2 and 5.3]{Mat08}.
Consequently, \eqref{eq:twist_kida_subtract} becomes
\begin{align}
2\lambda^\star(E^D/\Q)-2\lambda^\star(E/\Q)
=\sum_{\substack{\ell\mid D\\\ell\mid N_E}}2^{n_\ell+1}
+\sum_{\substack{\ell\mid D\\
\ell\nmid N_E,2\mid\sharp\widetilde E_\ell(\F_\ell)}}
2^{n_\ell+2}.
\end{align}
Dividing by $2$ proves \eqref{eq:Matsuno-type_formula}.
\end{proof}

We derive several corollaries of Theorem
\ref{chromatic_Matsuno-type_formula} by elementary arguments.

\begin{cor}\label{cor:parallel_chromatic_lambda}
Let $E$ be an elliptic curve over $\Q$ with good
supersingular reduction at $2$ and square-free
conductor. Assume that
$\sel^{\sharp}(E/\Q^{\cyc})$, $\sel^{\flat}(E/\Q^{\cyc})$ are $\Lambda$-cotorsion and $\mu^\sharp(E/\Q)=\mu^\flat(E/\Q)=0$.
Then, for every positive square-free integer
$D\equiv1\pmod4$, we have
\begin{align}\label{eq:parallel_chromatic_lambda}
\lambda^\sharp(E^D/\Q)-\lambda^\flat(E^D/\Q)
=
\lambda^\sharp(E/\Q)-\lambda^\flat(E/\Q).
\end{align}
\end{cor}

\begin{proof}
Subtracting the instances of \eqref{eq:Matsuno-type_formula}
for $\star=\sharp$ and $\star=\flat$ proves the assertion, since the local correction terms are independent
of $\star$ and hence cancel.
\end{proof}

\begin{rem}
Since $E[2]\cong E^D[2]$, Corollary \ref{cor:parallel_chromatic_lambda} is a $2$-adic analogue of the following result. Using a Greenberg--Vatsal-type congruence argument
for non-primitive plus/minus Selmer groups,
Kim \cite[Remark 3.3]{Kim09} showed that
for elliptic curves $E_1,E_2$ defined over $\Q$ with
good supersingular reduction at an odd prime $p$
satisfying $a_p(E_1)=a_p(E_2)=0$ and
$E_1[p]\simeq E_2[p]$ as $G_\Q$-modules,
\begin{align}
\lambda^+(E_1/\Q)-\lambda^-(E_1/\Q)=\lambda^+(E_2/\Q)-\lambda^-(E_2/\Q).
\end{align}
\end{rem}

\begin{cor}\label{cor:parallel_analytic_algebraic_lambda}
Let $E$ be an elliptic curve over $\Q$ with good
supersingular reduction at $2$ and square-free
conductor. Fix $\star\in\{\sharp,\flat\}$. Assume
that $\sel^{\star}(E/\Q^{\cyc})$ is $\Lambda$-cotorsion
and that
\begin{align}
\mu_{\an}^\star(E/\Q)=\mu_{\alg}^\star(E/\Q)=0.
\end{align}
Then, for every positive square-free integer
$D\equiv1\pmod4$, we have
\begin{align}\label{eq:parallel_analytic_algebraic_lambda}
\lambda_{\an}^\star(E^D/\Q)
-\lambda_{\alg}^\star(E^D/\Q)
=
\lambda_{\an}^\star(E/\Q)
-\lambda_{\alg}^\star(E/\Q).
\end{align}
In particular, if
$\lambda_{\an}^\star(E/\Q)=\lambda_{\alg}^\star(E/\Q)$,
then
$\lambda_{\an}^\star(E^D/\Q)=\lambda_{\alg}^\star(E^D/\Q)$
for every such $D$.

Moreover, if the above assumptions hold for both
$\star=\sharp$ and $\star=\flat$, then
\begin{align}\label{eq:parallel_four_lambda}
\begin{pmatrix}
\lambda_{\an}^\sharp(E^D/\Q)\\
\lambda_{\an}^\flat(E^D/\Q)\\
\lambda_{\alg}^\sharp(E^D/\Q)\\
\lambda_{\alg}^\flat(E^D/\Q)
\end{pmatrix}
=
\begin{pmatrix}
\lambda_{\an}^\sharp(E/\Q)\\
\lambda_{\an}^\flat(E/\Q)\\
\lambda_{\alg}^\sharp(E/\Q)\\
\lambda_{\alg}^\flat(E/\Q)
\end{pmatrix}
+
C_E(D)
\begin{pmatrix}
1\\1\\1\\1
\end{pmatrix},
\end{align}
where
\begin{align}
C_E(D):=
\sum_{\substack{\ell\mid D\\\ell\mid N_E}}
2^{n_\ell}
+
\sum_{\substack{\ell\mid D\\
\ell\nmid N_E,\,
2\mid\sharp\widetilde E_\ell(\F_\ell)}}
2^{n_\ell+1}.
\end{align}
\end{cor}

\begin{proof}
By Theorem \ref{chromatic_Matsuno-type_formula}
and \cite[Theorem 1.3]{Ada26+}, both
$\lambda_{\alg}^\star$ and
$\lambda_{\an}^\star$ increase by $C_E(D)$
under the quadratic twist by $D$.
Subtracting these formulas gives
\eqref{eq:parallel_analytic_algebraic_lambda},
while applying them for both
$\star=\sharp$ and $\star=\flat$ gives
\eqref{eq:parallel_four_lambda}.
\end{proof}

\subsection{Distributional results of chromatic $\lambda$-invariants}\label{subsec:distribution}

Applying the same argument as in \cite[\S\S 3,4]{Ada26+} to the difference formula \eqref{eq:Matsuno-type_formula}, we first obtain the following refinement of Matsuno's result \cite[Proposition 6.1]{Mat08} in the supersingular setting, determining the exact natural densities.

\begin{cor}[cf. {\cite[Corollary 1.4]{Ada26+}}]\label{Cor:Adachi26+_1.4}
Let $E$ be an elliptic curve defined over $\Q$ with good supersingular reduction at $2$ and square-free conductor, and let $\star\in \{\sharp,\flat\}$. Assume that $\sel^{\star}(E/\Q^{\cyc})$ is $\Lambda$-cotorsion and $\mu^{\star}(E/\Q)=0$. Then the set
\begin{align}
\{\lambda^{\star}(E^{\ell}/\Q)\ |\ \ell\text{ is a rational prime and }\ell\equiv 1\pmod 4\}
\end{align}
is infinite. More precisely, for each nonnegative integer $k$, the natural density (among all primes) of
\begin{align}
\{\ell\text{ is a rational prime and }\ell\equiv 1\pmod 4\ |\ \lambda^{\star}(E^{\ell}/\Q)=\lambda^{\star}(E/\Q)+2^{k+1}\}
\end{align}
is $\frac{1}{3\cdot 2^{k+1}}$.
\end{cor}

Hatley--Ray \cite[Theorem 4.1]{HR24+} established a lower bound for the number of quadratic twists with prescribed Iwasawa $\lambda$-invariants in the good ordinary case at $2$. Using the difference formula \eqref{eq:Matsuno-type_formula}, we obtain the following supersingular analogue for chromatic $\lambda$-invariants.

\begin{thm}[cf. {\cite[Theorem 1.5]{Ada26+}}]\label{Thm:Ada26+_1.5}
Let $E$ be an elliptic curve defined over $\Q$ with good supersingular reduction at $2$ and square-free conductor, and let $\star\in \{\sharp,\flat\}$. Assume that $\sel^{\star}(E/\Q^{\cyc})$ is $\Lambda$-cotorsion and $\mu^{\star}(E/\Q)=0$. Let $N^{\star}$ be a nonnegative integer such that $N^{\star}\geq\lambda^{\star}(E/\Q)$ and $\lambda^{\star}(E/\Q)\equiv N^{\star}\pmod 2$. Let
\begin{align}
m^{\star}_{E,N^{\star}}(X):=\sharp\{D\leq X\ |\ D\in\Z_{>0}\text{ is square-free, }D\equiv1\pmod 4,\text{ and }\lambda^{\star}(E^D/\Q)=N^{\star}\}.
\end{align}
Then we have
\begin{align}
m^{\star}_{E,N^{\star}}(X)\gg \frac{X}{(\log X)^{\frac{2}{3}}}.
\end{align}
\end{thm}

\begin{proof}
By Theorem \ref{chromatic_Matsuno-type_formula},
the algebraic chromatic $\lambda$-invariants
satisfy the same difference formula as the
analytic invariants in \cite[Theorem 1.3]{Ada26+}.
Thus the counting argument in
\cite[\S 4, proof of Theorem 1.5]{Ada26+}
applies without change.
\end{proof}

\begin{rem}
Related lower bounds for the number of quadratic twists
with prescribed $2$-Selmer ranks have been studied
by Chang \cite{Cha10} and Mazur--Rubin \cite{MR10}. Chang \cite[Theorem 1.1]{Cha10} proved that,
for any elliptic curve $E/\Q$ with $E(\Q)[2]=\{O\}$,
there exists $0<\alpha<1$ such that
\begin{align}
\sharp\left\{
D\in\Z\ \middle|
\begin{array}{l}
0<|D|\le X,\quad D\text{ ; square-free},\\
\dim_{\F_2}\sel_2(E^D/\Q)\le 1
\end{array}
\right\}
\gg \frac{X}{(\log X)^\alpha}.
\end{align}
Mazur--Rubin \cite[Theorem 1.4]{MR10} obtained
more precise lower bounds over general number fields for the number
$N_r(E,X)$ of quadratic twists with prescribed
$2$-Selmer rank $r$, counted by the norm of the conductor of the
corresponding quadratic extension. In the case $K=\Q$, their result gives
\begin{align}
N_r(E,X)\gg
\begin{cases}
\dfrac{X}{(\log X)^{2/3}}
& \text{if }\gal(\Q(E[2])/\Q)\simeq S_3,\\[6pt]
\dfrac{X}{(\log X)^{1/3}}
& \text{if }\gal(\Q(E[2])/\Q)\simeq \Z/3\Z,
\end{cases}
\end{align}
provided that such a twist exists;
see also \cite[Theorem 1.5]{MR10}.
In our setting, only the $S_3$ case occurs
(cf. \cite[Proposition 2.1]{Ada26+}),
and its lower bound has the same order as ours,
although the invariants differ.
\end{rem}
\section{$\Z_2$-corank for quadratic twists}\label{sec:MW-rank_twists}
For an elliptic curve $A$ defined over $\Q$, let $r_{\mathrm{alg}}(A):=\rank_{\Z}A(\Q)$, and let $\Sha(A/\Q)$ be the Tate--Shafarevich group of $A$.
In this section, we study the $\Z_2$-corank of
$2^\infty$-Selmer groups over $\Q$ for quadratic twists of elliptic curves based on
the method of Hatley--Ray \cite{HR24+}. Recall that there is an exact sequence
\begin{align}\label{seq:ext_M-W_Selmer_Sha}
0\longrightarrow A(\Q)\otimes\Q_2/\Z_2
\longrightarrow \sel_{2^\infty}(A/\Q)
\longrightarrow \Sha(A/\Q)[2^\infty]
\longrightarrow 0.
\end{align}
In particular,
\begin{align}
r_{\alg}(A)\leq
\corank_{\Z_2}\sel_{2^\infty}(A/\Q),
\end{align}
with equality if $\Sha(A/\Q)[2^\infty]$ is finite. Thus, a bound for the
$\Z_2$-corank of the $2^\infty$-Selmer group also gives a bound for the
Mordell--Weil rank.

\subsection{Chromatic $\lambda$-invariants and $\Z_2$-coranks}\label{subsec:Z_2-corank_chromatic_lambda}
We first give the chromatic analogue of the bound
on the classical Selmer corank.
At odd supersingular primes with $a_p(E)=0$,
related signed control results were established
in \cite[Lemma 9.1 and Theorem 9.3]{Kob03}
and \cite[Lemma 3.9]{Kim13}.
A comparison with the classical Selmer group
also appears in the proof of
\cite[Proposition 3.2]{Kim09},
while the local comparison below is analogous to
\cite[Proposition 3.16]{Kim07}.

\begin{lem}\label{lem:chromatic_local_bottom}
Let $A/\Q$ be an elliptic curve with good supersingular reduction at $2$, and let $\star\in\{\sharp,\flat\}$. Then
\begin{align}\label{eq:bottom_local_inclusion}
\rest_{\Q_{2,\infty}/\Q_2}\left(A(\Q_2)\otimes \Q_2/\Z_2\right)
\subset A^\star_{\infty,2}.
\end{align}
\end{lem}

\begin{proof}
Apply the construction in \S\ref{subsec:Coleman} to $A$, and write $a_2:=a_2(A)$. Since $A$ is supersingular at $2$,
\begin{align}
\sharp\widetilde A(\F_2)=3-a_2\in\{1,3,5\}.
\end{align}
Hence we have
\begin{align}\label{eq:formal_kummer_bottom}
A(\Q_2)\otimes \Q_2/\Z_2
=\what A(2\Z_2)\otimes \Q_2/\Z_2.
\end{align}
By Proposition \ref{Prop:generate_formal},
\begin{align}\label{eq:bottom_formal_generator}
\what A(2\Z_2)=\langle d_{-2}\rangle_{\Z_2}.
\end{align}
Moreover, Proposition \ref{Prop:trace_relation}, with $F_{\fp}=\Q_2$, gives
\begin{align}\label{eq:bottom_trace_d0}
\tr_{k_0/\Q_2}(d_0)
=\bigl(a_2(a_2-2)-1\bigr)d_{-2}.
\end{align}
The coefficient in \eqref{eq:bottom_trace_d0} belongs to $\Z_2^\times$. Write $\col'=(C_1,C_2)$. By Definition \ref{Definition: Coleman_map},
\begin{align}\label{eq:chromatic_components_bottom}
\col^\sharp=-a_2C_1+C_2,
\qquad
\col^\flat=-C_1.
\end{align}
Let $z=(z_n)_n\in\Ker\col^\star$, put $C_{1,0}:=\aug_\Gamma(C_1(z))\in R_0$, and let $\aug_\Delta:R_0\to\Z_2$ be the augmentation map for $\Delta$. By the construction of $\col'$ (cf.~\cite[(4.7)]{Kat21}),
\begin{align}
\aug_\Delta(C_{1,0})
&=\sum_{\sigma\in\Delta}(\sigma d_0,z_0)_{k_0}\notag\\
&=(\tr_{k_0/\Q_2}(d_0),z_0)_{k_0}\notag\\
&=\bigl(a_2(a_2-2)-1\bigr)(d_{-2},\core_{k_0/\Q_2}(z_0))_{\Q_2}.
\label{eq:bottom_pairing_U}
\end{align}
If $\star=\flat$, then $C_1(z)=0$ by \eqref{eq:chromatic_components_bottom}. Suppose that $\star=\sharp$. Since $\Imag(\col)=\Ker(q)$ in \eqref{eq:exact_Col}, we have $(w-a_2) C_{1,0}=0$. If $a_2=\pm2$, Proposition \ref{Prop:sharp_condition} gives $w-a_2=2u$ for some $u\in R_0^\times$, and hence $C_{1,0}=0$. If $a_2=0$, then $F_{\fp}=\Q_2$ gives $\varphi=1$ and $w-a_2=2N_\Delta$.
Thus $N_\Delta C_{1,0}=0$, since $R_0$ is $\Z_2$-torsion-free, and therefore
\begin{align}
2\aug_\Delta(C_{1,0})=\aug_\Delta(N_\Delta C_{1,0})=0.
\end{align}
Hence $\aug_\Delta(C_{1,0})=0$ in every case. By \eqref{eq:bottom_pairing_U} and \eqref{eq:bottom_formal_generator},
\begin{align}\label{eq:annihilate_bottom_formal}
(\what A(2\Z_2),\core_{k_0/\Q_2}(z_0))_{\Q_2}=0.
\end{align}
The compatibility of the local Tate pairing with restriction and corestriction gives, for $P\in \what A(2\Z_2)\otimes\Q_2/\Z_2$,
\begin{align}
\langle \rest_{\Q_{2,\infty}/\Q_2}(P),z\rangle
=
\langle P,\core_{k_0/\Q_2}(z_0)\rangle=0.
\end{align}
Together with \eqref{eq:formal_kummer_bottom} and Definition \ref{Definition: Local_condition}, this proves \eqref{eq:bottom_local_inclusion}.
\end{proof}

\begin{prop}\label{prop:chromatic_rank_bound}
Let $A/\Q$ be an elliptic curve with good supersingular reduction at $2$, and let $\star\in\{\sharp,\flat\}$. Suppose that $\cX^\star(A/\Q^{\cyc})$ is $\Lambda$-torsion. Then
\begin{align}\label{eq:chromatic_rank_bound}
r_{\alg}(A)
\leq \corank_{\Z_2}\sel_{2^\infty}(A/\Q)
\leq \lambda^\star(A/\Q).
\end{align}
In particular, if both chromatic Selmer groups are $\Lambda$-cotorsion, then
\begin{align}
\corank_{\Z_2}\sel_{2^\infty}(A/\Q)
\leq
\min\{\lambda^\sharp(A/\Q),\lambda^\flat(A/\Q)\}.
\end{align}
\end{prop}

\begin{proof}
We adapt the argument in the proof of
\cite[Proposition 3.2]{Kim09} to the chromatic
setting at $2$, using Lemma \ref{lem:chromatic_local_bottom} for the
required local comparison.

Put $G_\infty:=\gal(\Q_\infty/\Q)=\Delta\times\Gamma$. By the same argument as in Lemma \ref{lem:2-primary_torsion},
applied to $A$, we have
\begin{align}\label{eq:global_torsion_zero}
A(\Q_\infty)[2^\infty]=0.
\end{align}
Let $\Sigma$ contain $2$, the bad primes of $A$, and the
archimedean place. By \eqref{eq:global_torsion_zero},
the inflation--restriction sequence yields
\begin{align}\label{eq:global_inf_res_rank}
H^1(\Q_\Sigma/\Q,A[2^\infty])
\overset{\sim}{\longrightarrow}
H^1(\Q_\Sigma/\Q_\infty,A[2^\infty])^{G_\infty}.
\end{align}
Let
\begin{align}
x\in \sel_{2^\infty}(A/\Q),
\qquad
x_\infty:=\rest_{\Q_\infty/\Q}(x).
\end{align}
By \eqref{eq:global_inf_res_rank}, we have
\begin{align}
x_\infty\in
H^1(\Q_\Sigma/\Q_\infty,A[2^{\infty}])^{G_\infty}.
\end{align}
We show that $x_\infty\in \sel^\star(A/\Q_\infty)$. For $M=\Q$ or $\Q_\infty$ and a place $v$ of $M$, we write $\loc_v$ for the localization map at $v$. Since $\loc_2(x)\in A(\Q_2)\otimes\Q_2/\Z_2$, Lemma \ref{lem:chromatic_local_bottom} gives
\begin{align}\label{eq:local_condition_at_2}
\loc_2(x_\infty)
&=
\rest_{\Q_{2,\infty}/\Q_2}\bigl(\loc_2(x)\bigr) \notag\\
&\in
\rest_{\Q_{2,\infty}/\Q_2}
\left(A(\Q_2)\otimes\Q_2/\Z_2\right)
\subset
A^\star_{\infty,2}.
\end{align}
Let $\ell\neq2$ be a finite prime in $\Sigma$ and let
$w\mid\ell$ be a prime of $\Q_\infty$.  By the definition of the
classical Selmer group,
\begin{align}
\loc_\ell(x)
\in
A(\Q_\ell)\otimes\Q_2/\Z_2
\subset H^1(\Q_\ell,A[2^{\infty}]).
\end{align}
The Kummer maps are compatible with restriction, and hence
\begin{align}
\rest_{\Q_{\infty,w}/\Q_\ell}
\left(A(\Q_\ell)\otimes\Q_2/\Z_2\right)
\subset
A(\Q_{\infty,w})\otimes\Q_2/\Z_2.
\end{align}
Therefore
\begin{align}\label{eq:local_condition_away_2}
\loc_w(x_\infty)
&=
\rest_{\Q_{\infty,w}/\Q_\ell}
\bigl(\loc_\ell(x)\bigr) \notag\\
&\in
A(\Q_{\infty,w})\otimes\Q_2/\Z_2=0.
\end{align}
If $w$ is an archimedean place of $\Q_\infty$, then
$\Q_{\infty,w}=\C$, since $\Q_\infty$ is totally imaginary. Hence
\begin{align}\label{eq:archimedean_local_condition}
H^1(\Q_{\infty,w},A[2^{\infty}])
=
H^1(\C,A[2^{\infty}])
=
0.
\end{align}
Thus \eqref{eq:archimedean_local_condition} together with
\eqref{eq:local_condition_at_2}
and \eqref{eq:local_condition_away_2}
shows that $x_\infty\in\sel^\star(A/\Q_\infty)$.
By \eqref{eq:global_torsion_zero},
the inflation--restriction sequence identifies
$\rest_{\Q^{\cyc}/\Q}(x)$ with
$x_\infty\in\sel^\star(A/\Q_\infty)^\Delta$.
Moreover, since
$A(\R)\otimes_\Z\Q_2/\Z_2=0$, we have
$\loc_\infty(x)=0$.
By the compatibility of localization with restriction,
\begin{align}
\loc_v\parenth{\rest_{\Q^{\cyc}/\Q}(x)}
=
\rest_{(\Q^{\cyc})_v/\R}\parenth{\loc_\infty(x)}
=0
\end{align}
for every real place $v$ of $\Q^{\cyc}$.
Hence, by Definition
\ref{Definition_Selmer_chromatic_cyclotomic},
$\rest_{\Q^{\cyc}/\Q}(x)$ belongs to
$\sel^\star(A/\Q^{\cyc})^\Gamma$,
and restriction induces the injection
\begin{align}\label{eq:classical_to_chromatic}
\sel_{2^\infty}(A/\Q)
\lhook\joinrel\longrightarrow
\sel^\star(A/\Q^{\cyc})^\Gamma.
\end{align}
Taking Pontryagin duals in \eqref{eq:classical_to_chromatic}, we obtain a surjection
\begin{align}
\cX^\star(A/\Q^{\cyc})_\Gamma
\twoheadrightarrow
\sel_{2^\infty}(A/\Q)^\vee.
\end{align}
Hence
\begin{align}\label{eq:corank_coinvariant}
\corank_{\Z_2}\sel_{2^\infty}(A/\Q)
\leq
\rank_{\Z_2}\cX^\star(A/\Q^{\cyc})_\Gamma.
\end{align}
For a finitely generated torsion $\Lambda\cong\Z_2\llbracket T\rrbracket$-module $X$, the structure theorem gives
\begin{align}\label{eq:coinvariant_lambda_bound}
\rank_{\Z_2}X/TX\leq\lambda(X).
\end{align}
Applying \eqref{eq:coinvariant_lambda_bound} to $X=\cX^\star(A/\Q^{\cyc})$ and using \eqref{eq:corank_coinvariant} gives the second inequality in \eqref{eq:chromatic_rank_bound}. The first inequality in \eqref{eq:chromatic_rank_bound} follows from the exact sequence \eqref{seq:ext_M-W_Selmer_Sha}.
\end{proof}

Applying Proposition \ref{prop:chromatic_rank_bound} to $A=E^D$ gives the chromatic counterpart of \cite[(3.2)]{HR24+}:
\begin{align}\label{eq:chromatic_HR_bound}
r_{\alg}(E^D)
\leq
\corank_{\Z_2}\sel_{2^\infty}(E^D/\Q)
\leq
\lambda^\star(E^D/\Q)
\qquad (\star\in\{\sharp,\flat\}).
\end{align}

\subsection{Quadratic twists with Selmer corank one}\label{subsec:Z_2-corank_twists}

We combine the difference formula \eqref{eq:Matsuno-type_formula} with
Proposition \ref{prop:chromatic_rank_bound} to construct prime quadratic twists
whose $2^\infty$-Selmer groups have $\Z_2$-corank one based on the strategy
of Hatley--Ray \cite[\S3]{HR24+}. 
We also extend the first prime-twist family to
square-free products of its primes.
\begin{rem}
In a series of works, Smith
\cite{Smi26a,Smi26b,Smi25+} established that,
for an elliptic curve $E/\Q$, the $2^\infty$-Selmer coranks $0$ and $1$ each
occur in $50\%$ of its quadratic twists. These results concern the full family of
quadratic twists, whereas our results concern
explicit positive-density families of prime twists.
\end{rem}
Throughout this subsection, let $E/\Q$ have good supersingular reduction at $2$ and square-free conductor $N_E$. Put
\begin{align}
L_E:=\Q(E[2]),\qquad
k_E:=\Q(\sqrt{\Delta_E}),\qquad
K_E:=\Q(i,\sqrt{-N_E}).
\notag\end{align}
The natural densities of the prime-twist families considered below are determined by whether $k_E\subseteq K_E$. We denote the root
number of an elliptic curve $A/\Q$ by $\omega(A)$. The first result is the chromatic counterpart of
\cite[Theorem 3.13]{HR24+}.

\begin{thm}\label{thm:prime_twists_corank_one}
Let $\star\in\{\sharp,\flat\}$, and assume that
\begin{align}
\omega(E)=-1,\qquad
\mu^\star(E/\Q)=0,\qquad
\lambda^\star(E/\Q)\leq2.
\notag\end{align}
Let $\cM_E^-$ be the set of rational primes $\ell$ satisfying
\begin{itemize}
\item $\ell\nmid 2N_E$,
\item $\widetilde E_{\ell}(\F_\ell)[2]=\{\widetilde{O}\}$,
\item $\ell$ splits completely in $K_E$.
\end{itemize}
Then, for every $\ell\in\cM_E^-$, we have
\begin{align}
\corank_{\Z_2}\sel_{2^\infty}(E^\ell/\Q)=1.
\notag\end{align}
Moreover, $\cM_E^-$ has natural density
\begin{align}\label{eq:density_prime_twists}
\operatorname{dens}(\cM_E^-)=
\begin{cases}
\frac{1}{6}&\text{if }k_E\subseteq K_E,\\[4pt]
\frac{1}{12}&\text{otherwise}.
\end{cases}
\end{align}
\end{thm}

\begin{proof}
We first prove the Selmer-corank assertion. Let
$\ell\in\cM_E^-$. Since $\ell$ splits completely in $K_E$, we have
$\ell\equiv1\pmod4$ and
\begin{align}
\left(\frac{-N_E}{\ell}\right)=1.
\notag\end{align}
Moreover, $\widetilde E_{\ell}(\F_\ell)[2]=\{\widetilde{O}\}$ implies that
$2\nmid\sharp\widetilde E_{\ell}(\F_\ell)$. Thus, the local correction terms in
\eqref{eq:Matsuno-type_formula} vanish for $D=\ell$, and hence
\begin{align}\label{eq:lambda_prime_twist_unchanged}
\mu^\star(E^\ell/\Q)=0,\qquad
\lambda^\star(E^\ell/\Q)=\lambda^\star(E/\Q)\leq2.
\end{align}
By \cite[p. 338, Corollary]{Roh96}, we see that
\begin{align}\label{eq:root_number_prime_twist}
\omega(E^\ell)
=\chi_\ell(-N_E)\omega(E)
=\left(\frac{-N_E}{\ell}\right)\omega(E),
\end{align}
where $\chi_\ell$ is the quadratic character associated with
$\Q(\sqrt\ell)/\Q$; see \cite[Proposition 3.11]{HR24+}.
It follows that $\omega(E^\ell)=-1$.
By the $2$-parity theorem of Dokchitser--Dokchitser
\cite[Theorem 1.4]{DD10},
\begin{align}\label{eq:parity_prime_twist}
\corank_{\Z_2}\sel_{2^\infty}(E^\ell/\Q)\equiv1\pmod2.
\end{align}
On the other hand, Proposition \ref{prop:chromatic_rank_bound} and
\eqref{eq:lambda_prime_twist_unchanged} yield
\begin{align}
\corank_{\Z_2}\sel_{2^\infty}(E^\ell/\Q)
\leq\lambda^\star(E^\ell/\Q)\leq2.
\notag\end{align}
Together with \eqref{eq:parity_prime_twist}, this proves the Selmer-corank assertion.

We next calculate the density. Recall that
$\gal(L_E/\Q)\cong S_3$, where $S_3$ denotes
the symmetric group of degree three (cf. \cite[Proposition 2.1]{Ada26+}). Moreover, $k_E$ is the unique quadratic subfield of $L_E$. Since $N_E>1$ is odd and square-free, $K_E/\Q$ is biquadratic. As $K_E/\Q$ is abelian and the only non-trivial proper Galois subfield
of $L_E/\Q$ is $k_E$, we have
\begin{align}
L_E\cap K_E=\Q\quad\text{or}\quad L_E\cap K_E=k_E,
\end{align}
and consequently
\begin{align}\label{composite_field_degree}
[K_EL_E:\Q]=
\begin{cases}
12&\text{if }k_E\subseteq K_E,\\
24&\text{otherwise}.
\end{cases}
\end{align}
Put $F_E:=K_EL_E$. Since $L_E/k_E$ is cyclic of degree $3$ and
$[K_Ek_E:k_E]$ is a power of $2$, we have
\begin{align}
\gal(F_E/K_Ek_E)\cong\gal(L_E/k_E)\cong\Z/3\Z.
\notag\end{align}
Let $S$ consist of the two non-identity elements of
$\gal(F_E/K_Ek_E)$. In $\gal(F_E/\Q)$, this set is stable
under conjugation. For any unramified prime $\ell\nmid2N_E$, the
condition $\widetilde E_{\ell}(\F_\ell)[2]=\{\widetilde{O}\}$ is equivalent to
$\Frob_\ell|_{L_E}$ being a $3$-cycle in $S_3$. Thus, $\ell$ belongs
to $\cM_E^-$ if and only if $\Frob_\ell\in S$. By Chebotarev's density theorem, we obtain
\begin{align}\label{dens_M_E}
\operatorname{dens}(\cM_E^-)
=\frac{\sharp S}{[F_E:\Q]}
=\frac2{[K_EL_E:\Q]}.
\end{align}
The assertion follows from \eqref{composite_field_degree} and \eqref{dens_M_E}.
\end{proof}

The algebraic chromatic $\lambda$-invariants can be difficult to compute
directly. The next corollary gives a criterion in terms of the analytic
invariants, which can be evaluated using Pollack's SageMath code
\cite{PolCode} (see also \cite[\S5]{Ada26+}).

\begin{cor}\label{cor:prime_twists_analytic}
Let $\star\in\{\sharp,\flat\}$, and assume that
\begin{align}
\omega(E)=-1,\qquad
\mu_{\an}^\star(E/\Q)=0,\qquad
\lambda_{\an}^\star(E/\Q)\leq2.
\notag\end{align}
Then every $\ell\in\cM_E^-$ satisfies
\begin{align}
\corank_{\Z_2}\sel_{2^\infty}(E^\ell/\Q)=1.
\notag\end{align}
\end{cor}

\begin{proof}
For $\ell\in\cM_E^-$, the analytic difference formula
\cite[Theorem 1.3]{Ada26+} has no local correction term. Hence
\begin{align}
\mu_{\an}^\star(E^{\ell}/\Q)=0,\qquad
\lambda_{\an}^\star(E^{\ell}/\Q)
=\lambda_{\an}^\star(E/\Q)\leq2.
\notag\end{align}
By Propositions \ref{prop:alg_an_lambda} and
\ref{prop:chromatic_rank_bound},
\begin{align}
\corank_{\Z_2}\sel_{2^\infty}(E^{\ell}/\Q)
\leq\lambda_{\alg}^\star(E^{\ell}/\Q)
\leq\lambda_{\an}^\star(E^{\ell}/\Q)\leq2.
\notag\end{align}
The root number computation in the proof of Theorem
\ref{thm:prime_twists_corank_one} shows that $\omega(E^\ell)=-1$.
The $2$-parity theorem \cite[Theorem 1.4]{DD10} therefore forces
\begin{align}
\corank_{\Z_2}\sel_{2^\infty}(E^\ell/\Q)=1.
\notag\end{align}
\end{proof}

\begin{cor}\label{cor:composite_twists_corank_one}
Assume the hypotheses of Theorem
\ref{thm:prime_twists_corank_one}.
Let $\cD_E$ be the set of square-free integers
$D>1$ whose prime factors all belong to $\cM_E^-$.
Then
\begin{align}
\corank_{\Z_2}\sel_{2^\infty}(E^D/\Q)=1
\qquad(D\in\cD_E).
\end{align}
Moreover,
\begin{align}
\sharp\{D\in\cD_E\mid D\le X\}&\gg \frac{X}{(\log X)^{1-\mathrm{dens}(\cM_E^-)}}\\
&=
\begin{cases}
\dfrac{X}{(\log X)^{5/6}} & \text{if } k_E\subseteq K_E,\\[8pt]
\dfrac{X}{(\log X)^{11/12}} & \text{otherwise}.
\end{cases}
\end{align}
\end{cor}

\begin{proof}
For $D\in\cD_E$, we have $D\equiv1\pmod4$,
and all local correction terms in
\eqref{eq:Matsuno-type_formula} vanish. Hence
\begin{align}
\mu^\star(E^D/\Q)=0,\qquad
\lambda^\star(E^D/\Q)=\lambda^\star(E/\Q)\le2.
\end{align}
Moreover,
\begin{align}
\omega(E^D)
=\omega(E)\prod_{\ell\mid D}
\left(\frac{-N_E}{\ell}\right)=-1.
\end{align}
Thus the corank assertion follows from Proposition
\ref{prop:chromatic_rank_bound} and the $2$-parity theorem.
The counting bound follows by the Selberg--Delange method
for the Chebotarev set $\cM_E^-$,
as in \cite[Theorem 3.12]{HR24+}.
\end{proof}

We next consider the complementary case $\omega(E)=+1$.
Unlike Theorem \ref{thm:prime_twists_corank_one},
where the chromatic $\lambda$-invariant does not change
under quadratic twists, we choose primes for which the local
correction term in \eqref{eq:Matsuno-type_formula}
is exactly $2$.

\begin{thm}\label{thm:prime_twists_corank_one_plus}
Let $\star\in\{\sharp,\flat\}$, and assume that
\begin{align}
\omega(E)=+1,\qquad
\mu^\star(E/\Q)=0,\qquad
\lambda^\star(E/\Q)=0.
\notag\end{align}
Let $\cM_E^+$ be the set of rational primes $\ell$ satisfying
\begin{itemize}
\item $\ell\nmid2N_E$,
\item $\ell\equiv5\pmod8$,
\item $\widetilde E_{\ell}(\F_\ell)[2]\neq\{\widetilde{O}\}$,
\item $\displaystyle\left(\frac{-N_E}{\ell}\right)=-1$.
\end{itemize}
Then, for every $\ell\in\cM_E^+$, we have
\begin{align}\label{eq:prime_twists_plus_conclusion}
\mu^\star(E^\ell/\Q)=0,\qquad
\lambda^\star(E^\ell/\Q)=2,\qquad
\corank_{\Z_2}\sel_{2^\infty}(E^\ell/\Q)=1.
\end{align}
Moreover, $\cM_E^+$ has natural density
\begin{align}\label{eq:density_prime_twists_plus}
\operatorname{dens}(\cM_E^+)=
\begin{cases}
\frac{1}{8}&\text{if }k_E\subseteq K_E,\\[4pt]
\frac{1}{12}&\text{otherwise}.
\end{cases}
\end{align}
\end{thm}

\begin{proof}
We first prove \eqref{eq:prime_twists_plus_conclusion}. For
$\ell\in\cM_E^+$, we have $\ell\equiv1\pmod4$ and
\begin{align}
n_\ell=v_2\parenth{\frac{\ell^2-1}{8}}=0.
\notag\end{align}
Since $\ell\nmid N_E$ and $\widetilde E_{\ell}(\F_\ell)[2]\neq\{\widetilde{O}\}$,
\eqref{eq:Matsuno-type_formula} gives
\begin{align}
\mu^\star(E^\ell/\Q)=0,\qquad
\lambda^\star(E^\ell/\Q)
=\lambda^\star(E/\Q)+2^{n_\ell+1}=2.
\notag\end{align}
By \eqref{eq:root_number_prime_twist},
\begin{align}
\omega(E^\ell)
=\left(\frac{-N_E}{\ell}\right)\omega(E)
=(-1)(+1)=-1.
\notag\end{align}
The $2$-parity theorem \cite[Theorem 1.4]{DD10}, combined with
Proposition \ref{prop:chromatic_rank_bound}, yields
\begin{align}
\corank_{\Z_2}\sel_{2^\infty}(E^\ell/\Q)\equiv1\pmod2,
\qquad
\corank_{\Z_2}\sel_{2^\infty}(E^\ell/\Q)\leq2.
\notag\end{align}
Thus we obtain $\corank_{\Z_2}\sel_{2^\infty}(E^\ell/\Q)=1$. 
It remains to compute the density. Put
\begin{align}
C:=\Q(\zeta_8),\qquad K:=\Q(\sqrt{-N_E}).
\notag\end{align}
Recall that $L_E\cap C=\Q$
(cf. \cite[Lemma 3.23]{Ada26+}).
Together with $\gal(L_E/\Q)\cong S_3$, this gives
\begin{align}
\gal(L_EC/\Q)
\cong S_3\times(\Z/2\Z)^{\oplus2}.
\notag\end{align}
The quadratic subfields of $L_EC$ are precisely those
of $k_EC$. Since $N_E>1$ is odd and square-free,
and $2\nmid\Delta_E$, both $K$ and $k_E$ have odd
square-free radicands. As the quadratic subfields
of $C$ are $\Q(i)$, $\Q(\sqrt{2})$, and
$\Q(\sqrt{-2})$, we obtain
\begin{align}\label{eq:common_density_equivalence}
K\subseteq L_EC
&\iff K=k_E\ \text{or}\
K=\Q(\sqrt{-\Delta_E})\notag\\
&\iff k_E\subseteq K_E.
\end{align}
For $\ell\nmid2N_E$, the condition $\ell\equiv5\pmod8$ fixes
$\Frob_\ell|_C$. The condition
$\widetilde E_{\ell}(\F_\ell)[2]\neq\{\widetilde{O}\}$ is equivalent to
$\Frob_\ell|_{L_E}$ being the identity or a transposition in $S_3$.
If $k_E\not\subseteq K_E$, then
\eqref{eq:common_density_equivalence} gives $K\not\subseteq L_EC$.
Consequently,
\begin{align}
\gal(L_ECK/\Q)\cong S_3\times(\Z/2\Z)^{\oplus3},
\notag\end{align}
and the three Frobenius conditions are independent. By the Chebotarev density theorem,
\begin{align}
\operatorname{dens}(\cM_E^+)
=\frac14\cdot\frac46\cdot\frac12=\frac1{12}.
\notag\end{align}
If $k_E\subseteq K_E$, we have $K=k_E$ or $K=\Q(\sqrt{-\Delta_E})$.
For $\ell\equiv5\pmod8$ we have
$\left(\frac{-1}{\ell}\right)=1$, so in either case
\begin{align}
\left(\frac{-N_E}{\ell}\right)
=\sgn(\Frob_\ell|_{L_E}).
\notag\end{align}
Thus, the condition on the Legendre symbol in the definition of
$\cM_E^+$ forces $\Frob_\ell|_{L_E}$ to be a transposition.
There are three transpositions in $S_3$; all fix a non-zero vector
of $E[2]$. Hence
\begin{align}
\operatorname{dens}(\cM_E^+)
=\frac14\cdot\frac36=\frac18.
\notag\end{align}
\end{proof}

\begin{rem}
The condition $k_E\subseteq K_E$, which governs
the density formulae \eqref{eq:density_prime_twists}
and \eqref{eq:density_prime_twists_plus},
admits a simple criterion in terms of the
minimal discriminant:
\begin{align}
k_E\subseteq K_E
\quad\Longleftrightarrow\quad
v_\ell(\Delta_E)\equiv1\pmod2
\quad\text{for every prime }\ell\mid N_E.
\end{align}
\end{rem}

As in Corollary \ref{cor:prime_twists_analytic},
we obtain the following version of Theorem
\ref{thm:prime_twists_corank_one_plus}
with hypotheses on the analytic Iwasawa invariants.

\begin{cor}\label{cor:prime_twists_plus_analytic}
Assume that $\omega(E)=+1$, and let
$\star\in\{\sharp,\flat\}$. If
\begin{align}
\mu_{\an}^\star(E/\Q)=0,\qquad
\lambda_{\an}^\star(E/\Q)=0,
\notag\end{align}
then every $\ell\in\cM_E^+$ satisfies
\begin{align}
\corank_{\Z_2}\sel_{2^\infty}(E^\ell/\Q)=1.
\notag\end{align}
Moreover, suppose that $L(E,1)\neq0$ and the
$2$-primary part of the rank-zero
Birch--Swinnerton-Dyer formula holds for $E$.
Then $\mu_{\an}^\sharp(E/\Q)$ and
$\lambda_{\an}^\sharp(E/\Q)$ cannot both vanish, and the condition
\begin{align}
\mu_{\an}^\flat(E/\Q)=\lambda_{\an}^\flat(E/\Q)=0
\end{align}
is equivalent to
\begin{align}
\Sha(E/\Q)[2^\infty]=0
\text{ and } k_E\subseteq K_E.
\end{align}
In particular, whenever these equivalent conditions hold,
\begin{align}
\operatorname{dens}(\cM_E^+)=\frac{1}{8}.
\notag\end{align}
\end{cor}

\begin{proof}
The first assertion follows from Theorem
\ref{thm:prime_twists_corank_one_plus}
by the same argument as in the proof of
Corollary \ref{cor:prime_twists_analytic}.
Suppose now that the $2$-primary part of the
rank-zero Birch--Swinnerton-Dyer formula
holds for $E$.
By \cite[Convention 2.1 and Proposition 3.3]{Ada26+},
the Shimura period $\Omega_{f_E}^+$ and the
N\'eron period $\Omega_E^+$ satisfy
\begin{align}\label{eq:period_unit_cor48}
\frac{\Omega_{f_E}^+}{\Omega_E^+}
\in\Z_2^\times.
\end{align}
This comparison also applies to non-optimal elliptic
curves and incorporates the factor $2$ arising when
the real locus is disconnected.
The formulae in \cite[Table 2]{Spr17}, expressed in the Shimura-period
normalization used here, give the following values at the trivial character:
\begin{align}\label{eq:special_values_cor48}
L_2^\star(E,0)=b_\star(a_2)\frac{L(E,1)}{\Omega_{f_E}^+},
\end{align}
where
\begin{align}
b_\sharp(X)=-X^3+2X^2+3X-4,\quad b_\flat(X)=-X^2+2X+1.\notag
\end{align}
Since $a_2\in\{-2,0,2\}$, these factors satisfy
\begin{align}\label{eq:parity_special_values_cor48}
\bigl(b_\sharp(-2),b_\sharp(0),b_\sharp(2)\bigr)=(6,-4,2),\quad \bigl(b_\flat(-2),b_\flat(0),b_\flat(2)\bigr)=(-7,1,1).
\end{align}
In particular, $b_\sharp(a_2)\in2\Z_2$ and
$b_\flat(a_2)\in\Z_2^\times$.
Let $c_\ell(E)$ be the local Tamagawa number of $E$ at $\ell$.
Since $E(\Q)[2]=\{O\}$,
the $2$-primary part of the rank-zero Birch--Swinnerton-Dyer formula gives
\begin{align}\label{eq:BSD_valuation_cor48}
v_2\parenth{\frac{L(E,1)}{\Omega_E^+}}&=v_2\parenth{\sharp\Sha(E/\Q)[2^\infty]}+\sum_{\ell\mid N_E}v_2(c_\ell(E))
\geq0.
\end{align}
It follows from \eqref{eq:period_unit_cor48},
\eqref{eq:special_values_cor48}, \eqref{eq:parity_special_values_cor48}, and \eqref{eq:BSD_valuation_cor48} that
$L_2^\sharp(E,0)\in2\Z_2$.
Consequently, $\mu_{\an}^{\sharp}(E/\Q)$ and
$\lambda_{\an}^{\sharp}(E/\Q)$ cannot both vanish.
For the flat sign, $b_\flat(a_2)$ is a $2$-adic unit.
Therefore, \eqref{eq:BSD_valuation_cor48} and
\eqref{eq:period_unit_cor48} yield the exact valuation formula
\begin{align}\label{eq:flat_special_value_valuation_cor48}
v_2\parenth{L_2^\flat(E,0)}
&=v_2\parenth{\frac{L(E,1)}{\Omega_E^+}}\notag\\
&=v_2\parenth{\sharp\Sha(E/\Q)[2^\infty]}
 +\sum_{\ell\mid N_E}v_2(c_\ell(E)).
\end{align}
Since $L_2^\flat(E,T)\in\Z_2\llbracket T\rrbracket$
by \cite[Theorem 2.14]{Spr17}, its analytic $\mu$- and
$\lambda$-invariants both vanish if and only if
$L_2^\flat(E,0)$ is a $2$-adic unit.

As $N_E$ is square-free, the reduction is
multiplicative at every prime $\ell\mid N_E$.
For the Kodaira type $I_n$, where $n=v_\ell(\Delta_E)$,
the Tamagawa number is $n$ in the split case,
and is $1$ or $2$ according as $n$ is odd or even
in the nonsplit case. Then we have
\begin{align}
c_\ell(E)\in\Z_2^\times
\quad\Longleftrightarrow\quad
v_\ell(\Delta_E)\equiv1\pmod2.
\notag\end{align}
By the criterion in the preceding remark, all
$c_\ell(E)$ are odd exactly when $k_E\subseteq K_E$.
It follows from \eqref{eq:flat_special_value_valuation_cor48}
that
\begin{align}\label{eq:flat_BSD_criterion_cor48}
\mu_{\an}^\flat(E/\Q)
=\lambda_{\an}^\flat(E/\Q)=0
\quad\Longleftrightarrow\quad
\Sha(E/\Q)[2^\infty]=0
\ \text{and}\ k_E\subseteq K_E.
\end{align}
Therefore \eqref{eq:density_prime_twists_plus} gives
the stated density.
\end{proof}

\begin{rem}\label{Rem:relation_previous_works}
We mention some previous works on $2$-Selmer ranks
in quadratic twists.
\begin{enumerate}
\item Barrera Salazar--Pacetti--Tornar\'{i}a
\cite[Theorems 3.3 and 3.7]{BPT21}
studied prime twists using $2$-descent and
a modified class group $\Cl_*(A_\Q,E)$
attached to the cubic field
$A_\Q:=\Q[x]/(F(x))$, where
$E:y^2=F(x)$ is a Weierstrass equation
with $F(x)\in \Z[x]$ monic and irreducible of degree $3$;
see also \cite[Theorem 2.16]{BPT21}.
Under their local hypotheses
\cite[Hypotheses 2.1]{BPT21}, their results imply
\begin{align}
\sel_2(E^\ell/\Q)=\sel_2(E/\Q)
\qquad(\ell\in\cM_E^-).
\end{align}
Indeed, every $\ell\in\cM_E^-$ is inert in $A_\Q$
and satisfies
\begin{align}
\ell\equiv1\pmod4,\qquad
\left(\frac{\Delta_E/N_E}{\ell}\right)
=
\left(\frac{-\Delta_E/N_E}{\ell}\right)=1.
\end{align}
When $\Delta_E>0$, these conditions identify
$\cM_E^-$, up to finitely many primes, with
their set $\cC_{+,\square}$, and hence \cite[Theorem 3.7]{BPT21} implies that the $2$-Selmer group is unchanged under the quadratic twist by $\ell$. When $\Delta_E<0$, it follows instead from
\cite[Theorem 3.3]{BPT21}.
Moreover, in this case, we give the exact density of $\cM_E^-$ by \eqref{eq:density_prime_twists}. In particular, if $\omega(E)=-1$ and
$\dim_{\F_2}\sel_2(E/\Q)=1$, the $2$-parity theorem
gives
\begin{align}
\corank_{\Z_2}\sel_{2^\infty}(E^\ell/\Q)=1
\qquad(\ell\in\cM_E^-).
\end{align}
By contrast, Theorem \ref{thm:prime_twists_corank_one}
obtains the same conclusion using chromatic Iwasawa
invariants, without assuming $\dim_{\F_2}\sel_2(E/\Q)=1$
or imposing the local conditions of
\cite[Hypotheses 2.1]{BPT21}.

The family $\cM_E^+$ of Theorem
\ref{thm:prime_twists_corank_one_plus} is disjoint
from their inert-prime families, since
$\widetilde{E}_{\ell}(\F_\ell)[2]\neq\{\widetilde{O}\}$ for
$\ell\in\cM_E^+$. Moreover, for every $\ell\in\cM_E^+$,
Tate's algorithm (see \cite[Chapter IV, Section 9, Step 6]{Sil94}) gives
\begin{align}
c_\ell(E^\ell)
=\sharp\widetilde E_{\ell}(\F_\ell)[2]\in\{2,4\}.
\end{align}
Thus, $E^\ell$ fails the local condition
of \cite[Definition 1.6 and Remark 1.7]{BPT21}
at $\ell$.

\item Chang \cite[Theorem 1.2 and Corollary 1.4]{Cha06}
proved that, if $E(\Q)[2]=\{O\}$ and
$\sel_2(E/\Q)=0$, then there exists a set
$\mathcal P$ of primes of positive Dirichlet
density such that
\begin{align}
\sel_2(E^\ell/\Q)=0
\quad(\ell\in\mathcal P).
\end{align}
In particular,
\begin{align}
\corank_{\Z_2}\sel_{2^\infty}(E^\ell/\Q)=0
\quad(\ell\in\mathcal P).
\end{align}

\item Mazur--Rubin \cite[Theorems 1.4 and 1.5]{MR10}
showed that, in our $S_3$ setting, quadratic twists
with $2$-Selmer rank one, and hence $2^\infty$-Selmer
corank one, satisfy a lower bound
$\gg X/(\log X)^{2/3}$; see also
Chang \cite[Theorem 1.1]{Cha10}.
Under additional local hypotheses,
Barrera Salazar--Pacetti--Tornar\'{i}a
\cite[Corollary 3.4]{BPT21}
obtained related counting bounds.
Although these results give stronger counting bounds
for the full family of quadratic twists,
they do not imply the conclusion of
Corollary~\ref{cor:composite_twists_corank_one}
for every square-free $D$ whose prime factors
all belong to $\cM_E^-$.
Indeed, the primes in $\cM_E^-$ are of good reduction,
satisfy $\widetilde E_{\ell}(\F_\ell)[2]=\{\widetilde{O}\}$,
and split completely in $K_E$.
\end{enumerate}
\end{rem}

We give an example of Corollary
\ref{cor:prime_twists_analytic} which is not covered by the results of
Barrera Salazar--Pacetti--Tornar\'{i}a
\cite{BPT21}.

\begin{eg}
Consider the elliptic curve defined over $\Q$ with LMFDB label $185.b1$ \cite{lmfdb},
\begin{align}
E:\quad y^2+y=x^3-x^2-5x+6.
\notag\end{align}
It has $N_E=185=5\cdot37$, $\Delta_E=925=5^2\cdot37>0$,
$a_2(E)=0$, $\omega(E)=-1$, and $c_5(E)=2$.
According to the Iwasawa invariants recorded in
LMFDB \cite{lmfdb} based on Pollack's SageMath code \cite{PolCode}, we have
\begin{align}
\mu_{\an}^{\sharp}(E/\Q)=\mu_{\an}^{\flat}(E/\Q)=0,
\qquad
\parenth{\lambda_{\an}^{\sharp}(E/\Q),
\lambda_{\an}^{\flat}(E/\Q)}=(2,1).
\notag\end{align}
Thus, Corollary \ref{cor:prime_twists_analytic} applies with
$\star=\flat$. Since $k_E=\Q(\sqrt{37})$ is not contained in
$K_E=\Q(i,\sqrt{-185})$, the set $\cM_E^-$ has density $1/12$.
On the other hand, since $c_5(E)=2$, the condition $(\dagger.\mathrm{iii})$ in \cite[Definition 1.6]{BPT21} fails
at $5$; by \cite[Remark 1.7]{BPT21}, neither
$(\dagger.\mathrm{i})$ nor $(\dagger.\mathrm{ii})$ holds there.
Hence \cite[Hypotheses 2.1]{BPT21} do not apply to this example.
\end{eg}

\begin{rem}
A related Iwasawa-theoretic approach at the good
ordinary prime $p=3$ was developed by
Burungale--Skinner \cite[Theorems 2.8 and 2.9]{BS23}.
Using the congruence results of Greenberg--Vatsal
\cite{GV00}, they proved that, for certain elliptic
curves with a rational cyclic subgroup of order $3$,
a positive proportion of quadratic twists have
$3$-adic $\lambda$-invariant one.
By Perrin-Riou's $3$-adic Gross--Zagier formula
and Kolyvagin's theorem, they further showed
that these twists have Mordell--Weil rank one.
\end{rem}

\section{Examples}\label{sec:examples}

In this section, we give examples of Corollaries
\ref{cor:prime_twists_analytic} and
\ref{cor:prime_twists_plus_analytic}. 
We use the Iwasawa invariants recorded in LMFDB
\cite{lmfdb}, based on Pollack's SageMath code
\cite{PolCode}; see also \cite[\S 5]{Ada26+}.
The displayed analytic $2$-adic Iwasawa invariants
are obtained from the $p=2$ counterparts of the
plus/minus invariants introduced by Perrin-Riou
\cite{Per03}.
As explained in \cite[Remark 5.1]{Ada26+}, if both
plus/minus $\mu$-invariants vanish, then the corresponding $\mu$- and $\lambda$-invariants coincide with Sprung's chromatic invariants, with $+=\flat$ and $-=\sharp$.
In all the examples below,
\begin{align}
\mu_{\an}^{\sharp}(E/\Q)=\mu_{\an}^{\flat}(E/\Q)=0.
\notag\end{align}
The prime sets in the following tables are determined by testing the
splitting and $2$-torsion conditions in the respective corollaries.
For a minimal Weierstrass model
$E:y^2+\alpha_1xy+\alpha_3y
=x^3+\alpha_2x^2+\alpha_4x+\alpha_6$
and an odd prime $\ell$ of good reduction, the non-zero elements
of $\widetilde E_{\ell}(\F_\ell)[2]$ correspond to roots in $\F_\ell$ of
\begin{align}\label{eq:examples_2torsion_poly}
4(x^3+\alpha_2x^2+\alpha_4x+\alpha_6)
+(\alpha_1x+\alpha_3)^2.
\end{align}
Thus the listed primes can also be checked directly without
computing Selmer groups of individual twists.

\subsection{Examples with root number $-1$}

Consider the following five elliptic curves defined over $\Q$:
\begin{align}
&E_{37}:\quad y^2+y=x^3-x,\\
&E_{43}:\quad y^2+y=x^3+x^2,\\
&E_{91}:\quad y^2+y=x^3+x,\\
&E_{101}:\quad y^2+y=x^3+x^2-x-1,\\
&E_{185}:\quad y^2+y=x^3-x^2-5x+6.
\end{align}
Their LMFDB labels are, respectively, $37.a1$, $43.a1$,
$91.a1$, $101.a1$, and $185.b1$.
Each elliptic curve has square-free conductor, good supersingular reduction
at $2$, and root number $-1$.
The values of the $\lambda$-invariants and the densities of
$\cM_E^-$ are given in Table \ref{tab:examples_minus}.
The final column lists the first five primes in $\cM_E^-$.

\begin{table}[ht]
\centering
\small
\begin{tabular}{@{}lrrccl@{}}
\toprule
$E$ & $N_E$ & $\Delta_E$ &
$(\lambda_{\an}^{\sharp}(E/\Q),\lambda_{\an}^{\flat}(E/\Q))$ &
$\operatorname{dens}(\cM_E^-)$ & First five primes in $\cM_E^-$\\
\midrule
$37.a1$  & $37$  & $37$   & $(1,2)$ & $1/6$  & $41,53,73,101,149$\\
$43.a1$  & $43$  & $-43$  & $(1,2)$ & $1/6$  & $13,17,41,53,97$\\
$91.a1$  & $91$  & $-91$  & $(1,16)$& $1/6$  & $5,29,53,73,89$\\
$101.a1$ & $101$ & $101$  & $(1,4)$ & $1/6$  & $5,13,17,137,157$\\
$185.b1$ & $185$ & $925$  & $(2,1)$ & $1/12$ & $41,101,149,181,229$\\
\bottomrule
\end{tabular}
\caption{Examples of Corollary \ref{cor:prime_twists_analytic}.}
\label{tab:examples_minus}
\end{table}

For $E_{37},E_{43},E_{91}$, and $E_{101}$,
the square class of $\Delta_E$ is that of $N_E$ or $-N_E$;
thus $\Q(\sqrt{\Delta_E})\subseteq\Q(i,\sqrt{-N_E})$.
By \eqref{eq:density_prime_twists}, their sets $\cM_E^-$
have density $1/6$.
For $E_{185}$, however,
$\Q(\sqrt{\Delta_E})=\Q(\sqrt{37})$ is not contained in
$\Q(i,\sqrt{-185})$, giving density $1/12$.
For each prime in the last column of Table
\ref{tab:examples_minus}, we have
\begin{align}
\ell\equiv1\pmod4,\qquad
\left(\frac{-N_E}{\ell}\right)=1,\qquad
\widetilde E_{\ell}(\F_\ell)[2]=\{\widetilde{O}\}.
\notag\end{align}
We take $\star=\sharp$ for $E_{37},E_{43},E_{91}$, and $E_{101}$, and $\star=\flat$ for $E_{185}$. In each case, $\lambda^{\star}_{\an}(E/\Q)=1$. Therefore, Corollary
\ref{cor:prime_twists_analytic} gives
\begin{align}
\corank_{\Z_2}\sel_{2^\infty}(E^\ell/\Q)=1
\qquad\text{for every }\ell\in\cM_E,
\notag\end{align}
not merely for the primes displayed in the table.

\subsection{Examples with root number $+1$}

Next consider the three elliptic curves defined over $\Q$:
\begin{align}
E_{11}&:\quad y^2+y=x^3-x^2,\\
E_{19}&:\quad y^2+y=x^3+x^2+x,\\
E_{67}&:\quad y^2+y=x^3+x^2-12x-21.
\end{align}
Their LMFDB labels are $11.a3$, $19.a3$, and $67.a1$,
respectively. Each has good supersingular reduction at $2$,
square-free conductor, and root number $+1$.
For each of them we have
\begin{align}
\mu_{\an}^{\flat}(E/\Q)=0,\qquad
\lambda_{\an}^{\flat}(E/\Q)=0.
\notag\end{align}
In particular, Corollary \ref{cor:prime_twists_plus_analytic}
applies with $\star=\flat$.
Table \ref{tab:examples_plus} gives the invariants, densities,
and the first five primes in $\cM_E^+$.

\begin{table}[ht]
\centering
\small
\begin{tabular}{@{}lrrccl@{}}
\toprule
$E$ & $N_E$ & $\Delta_E$ &
$(\lambda_{\an}^{\sharp}(E/\Q),\lambda_{\an}^{\flat}(E/\Q))$ &
$\operatorname{dens}(\cM_E^+)$ & First five primes in $\cM_E^+$\\
\midrule
$11.a3$ & $11$ & $-11$ & $(1,0)$ & $1/8$ & $13,29,61,101,109$\\
$19.a3$ & $19$ & $-19$ & $(3,0)$ & $1/8$ & $13,29,37,53,109$\\
$67.a1$ & $67$ & $-67$ & $(1,0)$ & $1/8$ & $5,13,53,61,101$\\
\bottomrule
\end{tabular}
\caption{Examples of Corollary \ref{cor:prime_twists_plus_analytic}.}
\label{tab:examples_plus}
\end{table}
In these three cases,
since $\Q(\sqrt{\Delta_E})=\Q(\sqrt{-N_E})$, \eqref{eq:density_prime_twists_plus} gives
$\operatorname{dens}(\cM_E^+)=1/8$.
For each prime $\ell$ in the last column of Table
\ref{tab:examples_plus}, one readily checks that
\begin{align}
\ell\equiv5\pmod8,\qquad
\left(\frac{-N_E}{\ell}\right)=-1,\qquad
\sharp\widetilde E_{\ell}(\F_\ell)[2]=2.
\notag\end{align}
The last equality follows by counting the roots of
\eqref{eq:examples_2torsion_poly} modulo $\ell$. The analytic difference formula \cite[Theorem 1.3]{Ada26+}
then gives, for every $\ell\in\cM_E^+$,
\begin{align}
\mu_{\an}^{\flat}(E^\ell/\Q)=0,\qquad
\lambda_{\an}^{\flat}(E^\ell/\Q)=2.
\notag\end{align}
Since the root number of $E^\ell$ is $-1$, Corollary
\ref{cor:prime_twists_plus_analytic} yields
\begin{align}
\corank_{\Z_2}\sel_{2^\infty}(E^\ell/\Q)=1
\qquad\text{for every }\ell\in\cM_E^+.
\notag\end{align}
These examples lie outside the inert-prime families of
\cite{BPT21}, since $\widetilde E_{\ell}(\F_\ell)[2]\neq\{\widetilde{O}\}$
for every $\ell\in\cM_E^+$.

{
\bibliographystyle{abbrv}
\bibliography{biblio}
}

\end{document}